\documentclass{article}
\usepackage[T1]{fontenc}
\usepackage{amsthm, fullpage}
\usepackage{amsmath}
\usepackage{amssymb}
\usepackage{amsfonts}
\usepackage{eucal}
\usepackage[all]{xy}
\usepackage[frenchb, english]{babel}
\usepackage{pslatex}
\usepackage[pagebackref]{hyperref}

\newtheorem{prop}{Proposition}[subsection]
\newtheorem{theo}[prop]{Théor\`eme}
\newtheorem*{theo**}{Théorème}

\newtheorem{conj*}{Conjecture}[section]
\newtheorem{lemm}[prop]{Lemme}

\theoremstyle{definition}

\newtheorem{vide}[prop]{}
\newtheorem{defi}[prop]{Définition}

\theoremstyle{remark}
\newtheorem{rema}[prop]{Remarques}

\newtheorem{nota}[prop]{Notations}

\numberwithin{equation}{prop}

\newcommand{\riso}{ \overset{\sim}{\longrightarrow}\, }
\newcommand{\liso}{ \overset{\sim}{\longleftarrow}\, }

\newcommand{\Spf}{\mathrm{Spf}\,}

\renewcommand{\sp}{\mathrm{sp}}

\renewcommand{\AA}{{\mathcal{A}}}

\newcommand{\FF}{{\mathcal{F}}}
\newcommand{\B}{{\mathcal{B}}}

\newcommand{\E}{{\mathcal{E}}}
\newcommand{\G}{{\mathcal{G}}}

\newcommand{\M}{{\mathcal{M}}}

\newcommand{\D}{{\mathcal{D}}}

\newcommand{\PP}{{\mathcal{P}}}

\renewcommand{\O}{{\mathcal{O}}}

\newcommand{\V}{\mathcal{V}}
\renewcommand{\S}{\mathcal{S}}

\newcommand{\Y}{\mathcal{Y}}
\newcommand{\ZZ}{\mathcal{Z}}
\newcommand{\X}{\mathfrak{X}}

\newcommand{\U}{\mathfrak{U}}

\newcommand{\DD}{\mathbb{D}}
\renewcommand{\L}{\mathbb{L}}
\newcommand{\R}{\mathbb{R}}
\newcommand{\Q}{\mathbb{Q}}
\newcommand{\Z}{\mathbb{Z}}
\newcommand{\N}{\mathbb{N}}

\newcommand{\hdag}{  \phantom{}{^{\dag} }    }

\begin{document}
\title{Stabilité par produit tensoriel de la surholonomie}
\author{Daniel Caro
\footnote{L'auteur a bénéficié du soutien du réseau européen TMR \textit{Arithmetic Algebraic Geometry}
(contrat numéro UE MRTN-CT-2003-504917)\newline
}
}

\date{}

\maketitle

\begin{abstract}
\selectlanguage{english}
Let $\V$ be a complete discrete valued ring of mixed characteristic $(0,p)$, 
$K$ its field of fractions,  
$k$ its residue field which is supposed to be perfect. 
Let $X$ be a separated $k$-scheme of finite type and
$Y$ be a smooth open of $X$. 
We check that the equivalence of categories $\sp _{(Y,X),+}$ (from the category of  overconvergent isocrystals on $(Y,X)/K$ 
to that of overcoherent isocrystals on $(Y,X)/K$)
commutes with tensor products. 
Next, in Berthelot's theory of arithmetic $\D$-modules, 
we prove the stability under tensor products of the devissability in overconvergent isocrystals.
With Frobenius structures,
we get the stability under tensor products of the overholonomicity.
\end{abstract}

\selectlanguage{frenchb}

\tableofcontents

\section*{Introduction}

Soient $\V$ un anneau de valuation discrète complet d'inégales caractéristiques,
de corps résiduel parfait $k$ de caractéristique $p$,
de corps des fractions $K$.
Afin d'obtenir une cohomologie $p$-adique sur les $k$-variétés algébriques 
(i.e. les $k$-schémas séparés de type fini)
satisfaisant aux propriétés analogues à celles de la cohomologie étale $l$-adique 
(avec $l\not =p$) sur les $k$-variétés algébriques construite par Grothendieck dans les années 60,
Berthelot a construit une version arithmétique de la théorie des $\D$-modules (voir son introduction \cite{Beintro2}).
Après avoir défini (dans le cadre de sa théorie des $\D$-modules arithmétiques) la notion d'holonomie,
il avait conjecturé sa stabilité par les six opérations cohomologiques de Grothendieck, i.e. image directe, image directe extraordinaire, 
image inverse, image inverse extraordinaire, dualité et produit tensoriel. 
Dans le cadre des $k$-variétés quasi-projectives,
ces conjectures ont été validées  (voir \cite{caro-stab-holo}), 
leur preuve utilisant cependant les travaux de  \cite{caro_surholonome} et surtout \cite{caro-Tsuzuki} décrits ci-dessous.
Pour contourner ces conjectures, nous avions introduit la notion de surholonomie (voir \cite{caro_surholonome}).
Dans un travail en commun avec Nobuo Tsuzuki (soumis pour publication en 2008 et publié récemment dans \cite{caro-Tsuzuki}), 
nous avions établi la surholonomie des $F$-isocristaux surcohérents.
Comme corollaire, 
nous y avions indiqué via une référence à ce papier ci-présent (la première version date de 2007)
la stabilité de la surholonomie par les six opérations de Grothendieck, 
en particulier par produit tensoriel. En fait, depuis 2008, on bénéficie dans la théorie des $\D$-modules arithmétiques de plusieurs progrès:
1) si $X$ est une $k$-variété et $Y$ est un ouvert lisse de $X$, on a défini 
dans \cite{caro-pleine-fidelite} la notion d'isocristaux surcohérents sur $(Y,X)/K$, ce qui généralise 
les constructions de \cite{caro_devissge_surcoh} (où on traitait le cas où $X$ est propre et avec 
une structure de Frobenius) ; 2) on dispose d'une notion de surcohérence pour les systèmes inductifs de $\D$-modules arithmétiques 
(voir \cite{caro-stab-sys-ind-surcoh}). 
Pour obtenir une version la plus forte possible, grâce au progrès numéroté 1) ci-dessus, 
nous étendons dans ce papier la notion de dévissabilité en isocristaux surconvergents
pour les systèmes inductifs de $\D$-modules arithmétiques, ce qui étend naturellement 
la construction de \cite{caro-2006-surcoh-surcv} (où on ne traitait que le cas où les $\V$-schémas formels sont propres 
et où les $\D$-modules sont munis d'une structure de Frobenius). Modulo une équivalence de catégories (pour être précis, voir \ref{equi-devcoh-devqc}),
cette notion de dévissabilité correspond d'ailleurs à celle définie à la fin de \cite{caro-pleine-fidelite}.
Le but de ce papier est de vérifier la stabilité par produit tensoriel de cette dévissabilité en isocristaux surconvergents 
(resp. de la surcohérence avec structure de Frobenius) pour les systèmes inductifs de $\D$-modules arithmétiques.
Remarquons que ces deux stabilités par produits tensoriels 
sont des propriétés plus fortes que celle de la dévissabilité en isocristaux surconvergents au sens de 
\cite{caro-pleine-fidelite} ou celle de la surcohérence au sens de \cite{caro_surcoherent} (voir la remarque \ref{rema-stab-tens-qcdev++}).

\bigskip 
Décrivons à présent le contenu de ce papier.
Soit $(\PP, T,X,Y)$ un cadre, i.e. 
$\PP$ est un $\V$-schéma formel séparé et lisse, 
$T$ est un diviseur de sa fibre spéciale notée $P$, $X$ est une sous-$k$-variété fermée de 
$P$ et $Y := X \setminus T$.
On désigne par
$\smash{\D} ^\dag _{\PP} (\hdag T ) _{\Q}$, l'{\it anneau des opérateurs différentiels sur $\PP$
de niveau fini à singularités surconvergentes le long de $T $} (voir \cite[4.2.5]{Be1}).
Lorsque $Y$ est lisse, 
on dispose de l'équivalence 
de la catégorie 
des isocristaux surconvergents sur 
$(\PP, T,X,Y)$ dans celle des isocristaux surcohérents sur 
$(\PP, T,X,Y)$
que l'on note 
$\sp _{X \hookrightarrow \PP, T,+} \colon 
\mathrm{Isoc} ^{\dag}( \PP, T , X /K)
\cong 
\mathrm{Isoc} ^{\dag \dag}( \PP, T , X /K)$
(voir \cite{caro-pleine-fidelite} ou bien ici \ref{vide-eq-cat-spxPT+}).
 La catégorie $\mathrm{Isoc} ^{\dag \dag}( \PP, T , X /K)$
 est une sous-catégorie strictement pleine de celle des $\D ^\dag _{\PP} (\hdag T) _\Q$-modules cohérents à support
dans $X$.
Entre parenthèse, cette équivalence de catégories se note aussi 
$\sp _{(Y,X),+}\colon
\mathrm{Isoc}^{ \dag} (Y,X/K) \cong 
\mathrm{Isoc}^{\dag \dag} (Y,X/K)$,
car cela ne dépend pas des choix de $\PP$ et 
$T$.

Dans le premier chapitre, nous précisons les notations
et le cadre dans lequel nous travaillerons, avec parfois quelques nouvelles définitions. 
Par exemple, lorsque $Y$ est lisse, 
on définit $\mathrm{Isoc} ^{(\bullet)}( \PP, T , X /K)$ comme étant la sous-catégorie strictement pleine de
celle des 
$\smash{\widehat{\D}} _{\PP} ^{(\bullet)}(T)$-modules cohérents à lim-ind-isogénies près (voir \cite[2]{caro-stab-sys-ind-surcoh}
pour cette notion de cohérence)
caractérisée par l'équivalence de catégories
$\underrightarrow{\lim} 
\colon
\mathrm{Isoc} ^{(\bullet)}( \PP, T , X /K)\cong
\mathrm{Isoc} ^{\dag \dag}( \PP, T , X /K)$.
Nous rappelons aussi
les résultats les plus récents (éparpillés dans la littérature) que nous utiliserons. 
Dans le second chapitre, nous explicitons les définitions de produits tensoriels (interne ou externe) dans les différentes 
(i.e. dans le contexte des systèmes inductifs ou pas) 
catégories apparaissant en théorie des $\D$-modules arithmétiques
et nous vérifions quelques propriétés les concernant, e.g. quelques commutations du produit tensoriel 
aux images inverses extraordinaires, images directes etc. 
Nous établissons en outre, lorsque $X=P$, 
la stabilité par produit tensoriel de $\mathrm{Isoc} ^{(\bullet)}( \PP, T , P /K)$
(voir \ref{stab-isoc-bullet-pas}).
Lorsque $Y$ est lisse, nous prouvons dans la troisième partie que 
$\mathrm{Isoc} ^{(\bullet)}( \PP, T , X /K)$ et par conséquent
$\mathrm{Isoc} ^{\dag \dag}( \PP, T , X /K)$
est stable par produit tensoriel. De plus, on vérifie que 
l'équivalence de catégories
$\sp _{X \hookrightarrow \PP, T,+} $ commutent aux produits tensoriels.
On donne aussi un résultat analogue dans le cadre des $\V$-schémas formels faiblement lisses.
Dans la dernière partie, on introduit la notion de dévissabilité en isocristaux surconvergents pour les systèmes inductifs
de $\D$-modules arithmétiques sur des cadres quelconques et vérifions qu'elle étend celle 
donnée dans \cite{caro-2006-surcoh-surcv}. On valide ensuite la stabilité de cette dévissabilité par produit tensoriel. 
Avec des structures de Frobenius, cela entraîne comme attendu la stabilité de la surholonomie ou de la surcohérence par produit tensoriel.  
\bigskip

\noindent \textbf{Notations et convention}
La lettre $\V$ désigne un anneau de valuation discrète complet,
de corps résiduel parfait $k$ de caractéristique $p>0$, de corps des
fractions $K$ de caractéristique $0$.
Soit $\pi$ une uniformisante de $\V$.
De plus, $s\geq 1$ est un entier fixé et $F$ la puissance
$s$ème de l'endomorphisme de Frobenius. Les modules sont par défaut à gauche.
Si $\E$ est un faisceau abélien, on pose $\E _\Q:=\E \otimes _\Z \Q$.
On note $m$ un entier positif.

En général,
les $\V$-schémas formels faibles
(voir par exemple \cite{meredith-weakformalschemes} ou \cite{caro_devissge_surcoh})
sont désignés par des lettres romanes surmontées du symbole {\og $\dag$\fg},
les $\V$-schémas formels par des lettres calligraphiques ou gothiques,
les $k$-schémas par des lettres romanes
e.g.,
$P ^\dag$, $\PP$, $P $.
De plus, si $P ^\dag$ est un $\V$-schéma formel faible,
$\PP$ indique le $\V$-schéma formel induit par complétion $p$-adique,
$P $ sa fibre spéciale
et pour tout entier $i $, on notera $P _i$ les réductions modulo $\pi ^{i +1}$ de $P ^\dag $ 
(de même pour des relations analogues).
Une variété ou une $k$-variété signifie un $k$-schéma séparé et de type fini. 

Si $f\colon P ^{\dag  \prime} \rightarrow P ^{\dag }$ est un morphisme de $\V$-schémas formels faibles lisses,
par abus de notations, $f \colon \PP' \rightarrow \PP$ et
$f _0$ ou $f\colon P' \rightarrow P$ seront les morphismes induits.
On note $d _P$ la dimension de $P$.

Soient $\AA$ un faisceau d'anneaux sur un espace topologique $X$.
Si $\sharp$ est l'un des symboles $+$, $-$, ou $\mathrm{b}$, alors $D ^{\sharp} ( \AA )$ désigne
la catégorie dérivée des complexes de $\AA$-modules (à gauche par défaut) vérifiant les conditions correspondantes d'annulation
des faisceaux de cohomologie. Lorsque l'on souhaite préciser entre droite et gauche, on précise alors comme suit
$D ^{\sharp} ( \overset{ ^\mathrm{g}}{}\AA )$ ou $D ^{\sharp} ( \AA \overset{ ^\mathrm{d}}{})$.
Si on ne veut pas faire de choix entre à droite ou à gauche, on écrit
$D ^{\sharp} ( \overset{ ^*}{}\AA )$.
On note $D ^{\mathrm{b}} _{\mathrm{coh}} ( \AA )$
la sous-catégorie pleine de $D  ( \AA )$
des complexes à cohomologie cohérente et bornée.

On suppose (sans nuire à la généralité) que tous les
$k$-schémas sont réduits.
Si $\PP$ est un $\V$-schéma formel lisse, $T $ est un diviseur de $P$, alors,
pour alléger les notations, on notera 
$\smash{\widehat{\D}} _{\PP} ^{(m)} (T):=
\widehat{\B} ^{(m)} _{\PP} ( T)  \smash{\widehat{\otimes}} _{\O _{\PP}} \smash{\widehat{\D}} _{\PP} ^{(m)}$, 
où $\widehat{\B} ^{(m)} _{\PP} ( T) $ désigne les faisceaux d'anneaux construits par Berthelot dans
\cite[4.2.3]{Be1} et 
$\smash{\D} _{\PP} ^{(m)}$ est le faisceau des opérateurs différentiels de niveau $m$ sur $\PP$
(voir \cite[2.2]{Be1}).
Enfin, si $f \colon \X \to \Y$ est un morphisme de $\V$-schémas formels lisses,
pour tout entier $i \in \N$,
on note $ f _{i} \colon X _i \to Y _i$ le morphisme induit modulo $\pi ^{i+1}$.

Lorsqu'un diviseur est vide, 
nous omettrons de l'indiquer dans toutes
les notations le faisant intervenir.

\section{Notations, définitions, rappels et compléments}

\subsection{Catégories de systèmes inductifs en théorie de $\D$-modules arithmétiques}

Soient $\PP$ un $\V$-schéma formel lisse, 
$T$ un diviseur de sa fibre spéciale.

\begin{vide}
\label{limqc}
Soit $\sharp \in \{\emptyset, +,-, \mathrm{b}\}$.
On rappelle les notations de \cite[1.1 et 1.2]{caro-stab-sys-ind-surcoh}:
on dispose du système inductif d'anneaux 
$\smash{\widehat{\D}} _{\PP} ^{(\bullet)}(T) : =
(\smash{\widehat{\D}} _{\PP} ^{(m)}(T) )_{m\in \N}$ (les morphismes de transition sont construits dans \cite{Be1}).
Nous disposons de
la catégorie $M (\smash{\widehat{\D}} _{\PP} ^{(\bullet)} (T))$ 
des $\smash{\widehat{\D}} _{\PP} ^{(\bullet)} (T)$-modules
et
de 
$D ^{\sharp}( \smash{\widehat{\D}} _{\PP} ^{(\bullet)}(T))$, 
la catégorie dérivée des complexes de $\smash{\widehat{\D}} _{\PP} ^{(\bullet)} (T)$-modules à
cohomologie bornée selon $\sharp$.
Les objets de $M ( \smash{\widehat{\D}} _{\PP} ^{(\bullet)}(T))$ ou de
$D ( \smash{\widehat{\D}} _{\PP} ^{(\bullet)}(T))$
seront notés 
 $\E ^{(\bullet)}= (\E ^{(m)} , \alpha ^{(m',m)})$, 
 où $m,m'$ parcourent les entiers positifs tels que $m' \geq m$,
 où $\E ^{(m)} $ est un complexe de $\smash{\widehat{\D}} _{\PP} ^{(m)}(T)$-modules
 et $\alpha ^{(m',m)} \colon \E ^{(m)}\to \E ^{(m')}$ sont les morphismes $\smash{\widehat{\D}} _{\PP} ^{(m)}(T)$-linéaires de transition.

La catégorie obtenue en localisant $D ^{\sharp}( \smash{\widehat{\D}} _{\PP} ^{(\bullet)}(T))$
(resp. $M( \smash{\widehat{\D}} _{\PP} ^{(\bullet)}(T))$)
par rapport aux ind-isogénies se note
$\smash{\underrightarrow{D}} _{\Q} ^{\sharp}( \smash{\widehat{\D}} _{\PP} ^{(\bullet)}(T))$
(resp. $\smash{\underrightarrow{M}} _{\Q} ( \smash{\widehat{\D}} _{\PP} ^{(\bullet)}(T))$).
La catégorie obtenue en localisant 
$\smash{\underrightarrow{D}} ^{\sharp} _{\Q}
( \smash{\widehat{\D}} _{\PP} ^{(\bullet)}(T))$
par rapport aux lim-isomorphismes 
sera notée
$\smash{\underrightarrow{LD}} ^{\sharp} _{\Q}
( \smash{\widehat{\D}} _{\PP} ^{(\bullet)}(T))$
(resp. $\smash{\underrightarrow{LM}}  _{\Q}
( \smash{\widehat{\D}} _{\PP} ^{(\bullet)}(T))$).

On note
 $\underrightarrow{LM} _{\Q, \mathrm{coh}} (\smash{\widehat{\D}} _{\PP} ^{(\bullet)} (T))$
la sous-catégorie pleine de
$\underrightarrow{LM} _{\Q} (\smash{\widehat{\D}} _{\PP} ^{(\bullet)} (T))$
des $\smash{\widehat{\D}} _{\PP} ^{(\bullet)} (T)$-modules cohérents à lim-ind-isogénie près
(voir \cite[2]{caro-stab-sys-ind-surcoh}).
On note $\smash{\underrightarrow{LD}} ^{\mathrm{b}} _{\Q ,\mathrm{coh}}
(\smash{\widehat{\D}} _{\PP} ^{(\bullet)}(T))$
(resp. $\smash{\underrightarrow{LD}} ^{\mathrm{b}} _{\Q ,\mathrm{qc}}
(\smash{\widehat{\D}} _{\PP} ^{(\bullet)}(T))$)
la sous-catégorie pleine triangulée de 
$\smash{\underrightarrow{LD}} ^{\mathrm{b}} _{\Q}
(\smash{\widehat{\D}} _{\PP} ^{(\bullet)}(T))$
des complexes cohérents (resp. quasi-cohérents).

En passant à la limite inductive sur le niveau puis en tensorisant par $\Q$, on obtient le foncteur
$\underrightarrow{\lim}\colon \smash{\underrightarrow{LD}} ^{\mathrm{b}} _{\Q}
(\smash{\widehat{\D}} _{\PP} ^{(\bullet)}(T))
\rightarrow
D  ^{\mathrm{b}} ( \smash{\D} ^\dag _{\PP} (\hdag T) _{\Q} )$
(on prendra garde au fait que la notation pourrait être trompeuse car on n'indique pas cette tensorisation par $\Q$). 
D'après \cite[2.4.4]{caro-stab-sys-ind-surcoh},
celui-ci induit les équivalences de catégories
\begin{gather}
\label{eq-coh-lim}
\underrightarrow{\lim}  \colon
\smash{\underrightarrow{LD}} ^{\mathrm{b}} _{\Q ,\mathrm{coh}}
(\smash{\widehat{\D}} _{\PP} ^{(\bullet)}(T))
\cong
D ^\mathrm{b} _\mathrm{coh} ( \smash{\D} ^\dag _{\PP} (\hdag T) _{\Q} ),
\\
\label{M-eq-coh-lim}
\underrightarrow{\lim} 
\colon
\smash{\underrightarrow{LM}}  _{\Q, \mathrm{coh}}
(\smash{\widehat{\D}} _{\PP} ^{(\bullet)}(T))
\cong
\mathrm{Coh} ( \smash{\D} ^\dag _{\PP} (\hdag T) _{\Q} ).
\end{gather}
où $\mathrm{Coh} ( \smash{\D} ^\dag _{\PP} (\hdag T) _{\Q} )$ 
désigne la catégorie des 
$\smash{\D} ^\dag _{\PP} (\hdag T) _{\Q}$-modules cohérents.
\end{vide}

\begin{vide}
\label{defi-M-L}
Les définitions et propriétés de \ref{limqc}
restent valables en substituant  {\og $\B$ \fg}
à {\og $\D$ \fg}:
on note $M (\smash{\widehat{\B}} _{\PP} ^{(\bullet)} (T))$ la catégorie
des $\smash{\widehat{\B}} _{\PP} ^{(\bullet)} (T)$-modules
et
$D ^{\sharp} (\smash{\widehat{\B}} _{\PP} ^{(\bullet)} (T))$ la catégorie dérivée des
complexes de 
$\smash{\widehat{\B}} _{\PP} ^{(\bullet)} (T)$-modules à cohomologie bornée selon 
$\sharp \in \{\emptyset, +,-, \mathrm{b}\}$.
On note
$\smash{\underrightarrow{D}} _{\Q}  ^{\sharp} ( \smash{\widehat{\B}} _{\PP} ^{(\bullet)}(T))$,
$\smash{\underrightarrow{M}} _{\Q}  ( \smash{\widehat{\B}} _{\PP} ^{(\bullet)}(T))$,
$\underrightarrow{LD} _{\Q} (\smash{\widehat{\B}} _{\PP} ^{(\bullet)} (T))$,
$\underrightarrow{LM} ^{\sharp} _{\Q} (\smash{\widehat{\B}} _{\PP} ^{(\bullet)} (T))$
les catégories localisées de manière analogue à \ref{limqc}. 
On dispose aussi de la notion de
$\smash{\widehat{\B}} _{\PP} ^{(\bullet)} (T)$-modules
cohérents à lim-ind-isogénies près et on note
$\underrightarrow{LM} _{\Q,\mathrm{coh}} (\smash{\widehat{\B}} _{\PP} ^{(\bullet)} (T))$
la sous-catégorie pleine 
de
$\underrightarrow{LM} _{\Q} (\smash{\widehat{\B}} _{\PP} ^{(\bullet)} (T))$
des 
$\smash{\widehat{\B}} _{\PP} ^{(\bullet)} (T)$-modules
cohérents à lim-ind-isogénies près.
On dispose du foncteur canonique pleinement fidèle
$\underrightarrow{LM}  _{\Q, \mathrm{coh}}
(\widehat{\B} _{\PP } ^{(\bullet)} (T))
\to 
\underrightarrow{LD}  ^\mathrm{b} _{\Q, \mathrm{coh}}
(\widehat{\B} _{\PP } ^{(\bullet)} (T))$
et de l'équivalence de catégories
\begin{equation}
\label{eq-coh-lim-B}
\underrightarrow{\lim} \colon
\smash{\underrightarrow{LD}} ^{\mathrm{b}} _{\Q ,\mathrm{coh}}
(\smash{\widehat{\B}} _{\PP} ^{(\bullet)}(T))
\cong
D ^\mathrm{b} _\mathrm{coh} ( \O _{\PP} (\hdag T) _{\Q} ).
\end{equation}

\end{vide}

\subsection{Foncteurs locaux, image directe, image inverse extraordinaire, dual}
Soient $f \colon \PP' \to \PP$ un morphisme de $\V$-schémas formels lisses, 
$T$ et $T'$ des diviseurs respectifs de $P$ et $P'$ tels que
$f ( P '\setminus T' ) \subset P \setminus T$.

\begin{vide}
Pour tous $\E \in D _{\mathrm{qc}} ^\mathrm{b}
(\overset{^\mathrm{g}}{} \smash{\widehat{\D}} _{\PP} ^{(m)} (T ))$
et $\M \in D _{\mathrm{qc}} ^\mathrm{b}
(\smash{\widehat{\D}} _{\PP} ^{(m)} (T ) \overset{^\mathrm{d}}{})$,
on pose

\begin{gather} \notag
\M _i := \M \otimes ^\L _{\smash{\widehat{\D}} _{\PP} ^{(m)} (T )} \smash{\D} _{P _i} ^{(m)} (T ),\
\E _i := \smash{\D} _{P _i} ^{(m)} (T ) \otimes ^\L _{\smash{\widehat{\D}} _{\PP} ^{(m)} ( T )} \E,\\
\M \smash{\widehat{\otimes}} ^\L _{\smash{\widehat{\D}} _{\PP} ^{(m)} (T )} \E :=
\R \underset{\underset{i}{\longleftarrow}}{\lim}\, ( \M _i \otimes ^\L  _{\smash{\D} _{P _i} ^{(m)} (T )} \E _i).
\end{gather}
\end{vide}

\begin{vide}
\label{hdagT-nota}
Pour tout second diviseur $D \subset T $ de $P$,
on dispose respectivement des foncteurs oublis et extensions
\begin{gather}
\mathrm{oub} _{D, T}\colon 
\smash{\underrightarrow{LD}} ^{\mathrm{b}} _{\Q,\mathrm{qc}} ( \smash{\widehat{\D}} _{\PP} ^{(\bullet)}(T))
\to
\smash{\underrightarrow{LD}} ^{\mathrm{b}} _{\Q,\mathrm{qc}} ( \smash{\widehat{\D}} _{\PP} ^{(\bullet)}(D)),
\
\mathrm{oub} _{D, T}\colon 
D ^\mathrm{b} (\D ^\dag _{\PP} (\hdag T) _\Q) \to 
D ^\mathrm{b} (\D ^\dag _{\PP} (\hdag D) _\Q),
\\
 (\hdag T ,\,D) :=\D ^\dag _{\PP} (\hdag T) _\Q
\smash{\overset{\L}{\otimes}}   ^{\dag}
_{\D ^\dag  _{\PP} (\hdag D) _\Q} -
:=
\widehat{\B} ^{(\bullet)} _{\PP} ( T )  \smash{\widehat{\otimes}} ^\L
_{\widehat{\B} ^{(\bullet)}  _{\PP} ( D) } -
\colon
\smash{\underrightarrow{LD}} ^{\mathrm{b}} _{\Q,\mathrm{qc}} ( \smash{\widehat{\D}} _{\PP} ^{(\bullet)}(D))
\to
\smash{\underrightarrow{LD}} ^{\mathrm{b}} _{\Q,\mathrm{qc}} ( \smash{\widehat{\D}} _{\PP} ^{(\bullet)}(T)),
\\
(\hdag T, D) := 
\D ^\dag _{\PP} (\hdag T) _\Q \otimes _{ \D ^\dag _{\PP} (\hdag D) _\Q} - 
\colon 
D ^\mathrm{b} _\mathrm{coh} ( \D ^\dag _{\PP} (\hdag D) _\Q)
\to 
D ^\mathrm{b} _\mathrm{coh} (\D ^\dag _{\PP} (\hdag T) _\Q).
\end{gather}
Pour tout
$\E ^{(\bullet)} \in \underrightarrow{LD}  ^\mathrm{b} _{\Q, \mathrm{qc}}
(\overset{^\mathrm{g}}{} \smash{\widehat{\D}} _{\PP} ^{(\bullet)} (D))$,
on écrit aussi
$ \E ^{(\bullet)} (\hdag D ,\,T) :=(\hdag T ,\,D) (\E ^{(\bullet)})$.
Il nous arrivera d'omettre d'indiquer $D$, ce qui ne court aucun risque (voir par exemple la remarque \ref{rema-fct-qcoh2coh-coh}).
\end{vide}

\begin{rema}
\label{rema-fct-qcoh2coh-coh}
Soit $D \subset T$ un diviseur.
On note 
$ \underrightarrow{LD}  ^\mathrm{b} _{\Q, \mathrm{coh}}
(\overset{^\mathrm{g}}{} \smash{\widehat{\D}} _{\PP} ^{(\bullet)} (D ))
\cap 
 \underrightarrow{LD}  ^\mathrm{b} _{\Q, \mathrm{coh}}
(\overset{^\mathrm{g}}{} \smash{\widehat{\D}} _{\PP} ^{(\bullet)} (T ))$
la sous-catégorie pleine de
$ \underrightarrow{LD}  ^\mathrm{b} _{\Q, \mathrm{coh}}
(\overset{^\mathrm{g}}{} \smash{\widehat{\D}} _{\PP} ^{(\bullet)} (T ))$
des complexes 
$\E ^{(\bullet)} \in  \underrightarrow{LD}  ^\mathrm{b} _{\Q, \mathrm{coh}}
(\overset{^\mathrm{g}}{} \smash{\widehat{\D}} _{\PP} ^{(\bullet)} (D ))$. 
Les foncteurs $oub _{D,T}$ et $(\hdag T, ~D)$ induisent des équivalences quasi-inverses
entre
$ \underrightarrow{LD}  ^\mathrm{b} _{\Q, \mathrm{coh}}
(\overset{^\mathrm{g}}{} \smash{\widehat{\D}} _{\PP} ^{(\bullet)} (D ))
\cap 
 \underrightarrow{LD}  ^\mathrm{b} _{\Q, \mathrm{coh}}
(\overset{^\mathrm{g}}{} \smash{\widehat{\D}} _{\PP} ^{(\bullet)} (T ))$
et la sous-catégorie pleine de 
$ \underrightarrow{LD}  ^\mathrm{b} _{\Q, \mathrm{coh}}
(\overset{^\mathrm{g}}{} \smash{\widehat{\D}} _{\PP} ^{(\bullet)} (D ))$
dont les objets sont aussi des objets de
$  \underrightarrow{LD}  ^\mathrm{b} _{\Q, \mathrm{coh}}
(\overset{^\mathrm{g}}{} \smash{\widehat{\D}} _{\PP} ^{(\bullet)} (T ))$. 
De même, 
on note 
$D ^\mathrm{b} _\mathrm{coh} (\overset{^\mathrm{g}}{} \smash{\D} ^\dag _{\PP} (\hdag D) _{\Q} ) \cap 
D ^\mathrm{b} _\mathrm{coh} (\overset{^\mathrm{g}}{} \smash{\D} ^\dag _{\PP} (\hdag T) _{\Q} )$
la sous-catégorie pleine de 
$D ^\mathrm{b} _\mathrm{coh} (\overset{^\mathrm{g}}{} \smash{\D} ^\dag _{\PP} (\hdag T) _{\Q} )$
des complexes 
$\E \in D ^\mathrm{b} _\mathrm{coh} (\overset{^\mathrm{g}}{} \smash{\D} ^\dag _{\PP} (\hdag D) _{\Q} )$
;
les foncteurs 
$oub _{D,T}$ et $(\hdag T, ~D)$ induisent des équivalences quasi-inverses
entre
$D ^\mathrm{b} _\mathrm{coh} (\overset{^\mathrm{g}}{} \smash{\D} ^\dag _{\PP} (\hdag D) _{\Q} ) \cap 
D ^\mathrm{b} _\mathrm{coh} (\overset{^\mathrm{g}}{} \smash{\D} ^\dag _{\PP} (\hdag T) _{\Q} )$
et la sous-catégorie pleine de 
$D ^\mathrm{b} _\mathrm{coh} (\overset{^\mathrm{g}}{} \smash{\D} ^\dag _{\PP} (\hdag D) _{\Q} )$
dont les objets sont des complexes de
$D ^\mathrm{b} _\mathrm{coh} (\overset{^\mathrm{g}}{} \smash{\D} ^\dag _{\PP} (\hdag T) _{\Q} )$.

Grâce au corollaire 
\cite[3.5.4]{caro-stab-sys-ind-surcoh},
on vérifie que
le foncteur $\underrightarrow{\lim} $ 
se factorise en l'équivalence de catégories
\begin{equation}
\label{QcohTD}
\underrightarrow{\lim} 
\colon 
 \underrightarrow{LD}  ^\mathrm{b} _{\Q, \mathrm{coh}}
(\overset{^\mathrm{g}}{} \smash{\widehat{\D}} _{\PP} ^{(\bullet)} (D ))
\cap 
\underrightarrow{LD}  ^\mathrm{b} _{\Q, \mathrm{coh}}
(\overset{^\mathrm{g}}{} \smash{\widehat{\D}} _{\PP} ^{(\bullet)} (T ))
\cong 
D ^\mathrm{b} _\mathrm{coh} (\overset{^\mathrm{g}}{} \smash{\D} ^\dag _{\PP} (\hdag D) _{\Q} ) \cap 
D ^\mathrm{b} _\mathrm{coh} (\overset{^\mathrm{g}}{} \smash{\D} ^\dag _{\PP} (\hdag T) _{\Q} ).
\end{equation}
\end{rema}

\begin{nota}
\label{fct-qcoh2coh}
Soient $\X$, $\Y$ deux $\V$-schémas formels lisses, 
$D _1 \subset T _1$ deux diviseurs de $X$, $D _2 \subset  T _2$ deux diviseurs de $Y$.
Soient 
$\phi ^{(\bullet)}
\colon 
\smash{\underrightarrow{LD}} ^{\mathrm{b}} _{\Q,\mathrm{qc}}
( \smash{\widehat{\D}} _{\X} ^{(\bullet)}(T _1))
\to 
\smash{\underrightarrow{LD}} ^{\mathrm{b}} _{\Q,\mathrm{qc}}
( \smash{\widehat{\D}} _{\PP} ^{(\bullet)}(T))$
un foncteur
et
$\psi ^{(\bullet)}
\colon 
\smash{\underrightarrow{LD}} ^{\mathrm{b}} _{\Q,\mathrm{qc}}
( \smash{\widehat{\D}} _{\X} ^{(\bullet)}(T _1))
\times
\smash{\underrightarrow{LD}} ^{\mathrm{b}} _{\Q,\mathrm{qc}}
( \smash{\widehat{\D}} _{\Y} ^{(\bullet)}(T _2))
\to 
\smash{\underrightarrow{LD}} ^{\mathrm{b}} _{\Q,\mathrm{qc}}
( \smash{\widehat{\D}} _{\PP} ^{(\bullet)}(T))$
un bifoncteur.
On en déduit un foncteur 
$\mathrm{Coh} _{T _1} ( \phi ^{(\bullet)})
\colon 
D ^\mathrm{b} _\mathrm{coh} ( \smash{\D} ^\dag _{\X} (\hdag T _1) _{\Q} )
\to 
D ^\mathrm{b}  ( \smash{\D} ^\dag _{\PP} (\hdag T) _{\Q} ) $
en posant 
$\mathrm{Coh} _{T _1} ( \phi ^{(\bullet)}) := \underrightarrow{\lim} \circ \phi ^{(\bullet)}\circ (\underrightarrow{\lim} _{T _1} ) ^{-1}$,
où $(\underrightarrow{\lim} _{T_1})  ^{-1}$ désigne un foncteur quasi-inverse de 
l'équivalence de catégories de 
\begin{equation}
\label{eq-catLDBer-LD-D}
\underrightarrow{\lim} _{T _1}
\colon 
\underrightarrow{LD} ^{\mathrm{b}}  _{\Q, \mathrm{coh}} (\smash{\widehat{\D}} _{\X} ^{(\bullet)} (T _1))
\cong
D ^{\mathrm{b}} _{\mathrm{coh}}( \smash{\D} ^\dag _{\X} (\hdag T _1) _{\Q} ).
\end{equation}
De même, on obtient le bifoncteur
$\mathrm{Coh} _{T _1, T _2} ( \psi ^{(\bullet)})
\colon 
D ^\mathrm{b} _\mathrm{coh} ( \smash{\D} ^\dag _{\X} (\hdag T _1) _{\Q} )
\times 
D ^\mathrm{b} _\mathrm{coh} ( \smash{\D} ^\dag _{\Y} (\hdag T _2) _{\Q} )
\to 
D ^\mathrm{b}  ( \smash{\D} ^\dag _{\PP} (\hdag T) _{\Q} ) $
en posant 
$\mathrm{Coh} _{T _1, T _2} ( \psi ^{(\bullet)}) := 
\underrightarrow{\lim} \circ \psi ^{(\bullet)}\circ (\underrightarrow{\lim} _{T _1} \times \underrightarrow{\lim} _{T _2}) ^{-1}$.

D'après \cite[3.5.6]{caro-stab-sys-ind-surcoh},
les foncteurs 
$\mathrm{Coh} _{T _1} ( \phi ^{(\bullet)}) $
et
$\mathrm{Coh} _{D _1} ( \phi ^{(\bullet)}) $
sont isomorphes sur 
$D ^\mathrm{b} _\mathrm{coh} ( \smash{\D} ^\dag _{\X} (\hdag D _1) _{\Q} ) 
\cap D ^\mathrm{b} _\mathrm{coh} ( \smash{\D} ^\dag _{\X} (\hdag T _1) _{\Q} )$.
Pour les mêmes raisons,
les bifoncteurs 
$\mathrm{Coh} _{T _1, T _2} ( \psi ^{(\bullet)}) $
et
$\mathrm{Coh} _{D _1, D _2} ( \psi ^{(\bullet)}) $
sont isomorphes sur 
$(D ^\mathrm{b} _\mathrm{coh} ( \smash{\D} ^\dag _{\X} (\hdag D _1) _{\Q} ) 
\cap D ^\mathrm{b} _\mathrm{coh} ( \smash{\D} ^\dag _{\X} (\hdag T _1) _{\Q} ))
\times 
(D ^\mathrm{b} _\mathrm{coh} ( \smash{\D} ^\dag _{\Y} (\hdag D _2) _{\Q} ) 
\cap D ^\mathrm{b} _\mathrm{coh} ( \smash{\D} ^\dag _{\Y} (\hdag T _2) _{\Q} )
)$.
\end{nota}

\begin{vide}
On construit canoniquement (voir \cite[3.4, 3.5, 4.3]{Beintro2} et \cite[1.1.6--7]{caro_courbe-nouveau}) 
les foncteurs images directes par $f$ à singularités surconvergentes le long de $T$ et $T'$:
\begin{gather}
f  ^{(\bullet)} _{T,T',+} \colon
\smash{\underrightarrow{LD}} ^{\mathrm{b}} _{\Q,\mathrm{qc}} ( \smash{\widehat{\D}} _{\PP'} ^{(\bullet)}(T'))
\to 
\smash{\underrightarrow{LD}} ^{\mathrm{b}} _{\Q,\mathrm{qc}} ( \smash{\widehat{\D}} _{\PP} ^{(\bullet)}(T)),
\\
f _{T,T',+} \colon
D ^\mathrm{b} _\mathrm{coh} ( \smash{\D} ^\dag _{\PP'} (\hdag T') _{\Q} )
\to 
 D ^\mathrm{b} ( \smash{\D} ^\dag _{\PP} (\hdag T) _{\Q} ).
\end{gather}
On définit de plus canoniquement les foncteurs images inverses extraordinaires 
par $f$ à singularités surconvergentes le long de $T$ et $T'$:
\begin{gather}
f ^{(\bullet)!} _{T',T} \colon
\smash{\underrightarrow{LD}} ^{\mathrm{b}} _{\Q,\mathrm{qc}} ( \smash{\widehat{\D}} _{\PP} ^{(\bullet)}(T))
\to
\smash{\underrightarrow{LD}} ^{\mathrm{b}} _{\Q,\mathrm{qc}} ( \smash{\widehat{\D}} _{\PP'} ^{(\bullet)}(T')),
\\
f ^! _{T',T} \colon
D ^\mathrm{b} _\mathrm{coh} ( \smash{\D} ^\dag _{\PP} (\hdag T) _{\Q} )
\to 
 D ^\mathrm{b} ( \smash{\D} ^\dag _{\PP'} (\hdag T') _{\Q} ).
\end{gather}
\end{vide}

\begin{nota}
\label{nota-f*TT'}
On posera 
$\L f ^{(\bullet) *} _{T',T}:= f ^{(\bullet)!} _{T',T} [-d _{P'/P}]$
et
$\L f ^{*} _{T',T} := f ^{!} _{T',T} [-d _{P'/P}]$.
Si $f$ est lisse, on omet d'indiquer le symbole $\L$.
\end{nota}

\begin{vide}
\label{coh-Qcoh}
Avec les notations \ref{fct-qcoh2coh},
on dispose des égalités (à isomorphismes canoniques près):
$\mathrm{Coh} _{T'} (f  ^{(\bullet)} _{T,T',+} )=f  _{T,T',+}$,
$\mathrm{Coh} _{T}(f ^{(\bullet) !} _{T',T}) = f ^{!} _{T',T}$, 
$\mathrm{Coh} _{D} ((\hdag T, D) ) =(\hdag T, D) $
(voir par exemple \cite[4.3.2.2 et 4.3.7.1]{Beintro2}).
\end{vide}

\begin{vide}
\label{f+!-indTT'}
Lorsque $T' = f ^{-1} (T)$, on omet d'indiquer $T'$ dans les notations faisant intervenir $f$.
De plus, comme, modulo les foncteurs oublis, les foncteurs
$f ^{(\bullet)!} _T$ et $f ^{(\bullet)!} $ 
(resp. $f ^{(\bullet)} _{T,+}$ et $f ^{(\bullet)} _{+} $, 
resp. 
$f ^{(\bullet)*} _T$ et $f ^{(\bullet)*} $ )
sont canoniquement isomorphes (voir \cite[1.1.10--11]{caro_courbe-nouveau}),
on pourra sans risque omettre d'indiquer le diviseur $T$. 
De même, 
si on ne souhaite pas préciser le diviseur, on pourra noter simplement
$f ^!$, $f ^* $ et $f _{+}$ à la place de $f ^! _T$, $f ^* _T$ et $f _{T,+}$ (d'après \ref{fct-qcoh2coh}, cela est anodin). 
\end{vide}

\begin{vide}
Si $X$ est un sous-schéma fermé de $P$, $\R \underline{\Gamma} ^\dag
_X $ désigne le foncteur cohomologique local à support strict dans
$X$ (au sens de \cite[2.2.6]{caro_surcoherent}) et $(\hdag X)$ le
foncteur de localisation en dehors de $X$.
(\cite[2.2.6]{caro_surcoherent}). 
On pourra aussi consulter \cite[4]{caro-stab-sys-ind-surcoh} pour plus de précisions à leur propos.

Le foncteur dual
$\smash{\D} ^\dag _{\PP} (\hdag T) _{\Q}$-linéaire (voir
\cite[I.3.2]{virrion} pour la définition des foncteurs duaux)
se note $\DD _{\PP,T}$ ou $\DD _T$.

\end{vide}

\begin{vide}
\noindent $\bullet$ Soit $Y$ un sous-schéma de $P$. 
Pour tous sous-schémas fermés $X$ et $T$ de $P$ tels que $Y= X \setminus T$, 
pour tout $\E ^{(\bullet)} \in (F\text{-})\smash{\underrightarrow{LD}} ^{\mathrm{b}} _{\Q ,\mathrm{qc}}
(\smash{\widehat{\D}} _{\PP} ^{(\bullet)})$,
on pose
\begin{equation}
  \label{xtx't'}
\R\underline{\Gamma} ^\dag _{Y} (\E ^{(\bullet)}) :=
\R\underline{\Gamma} ^\dag _{X} (\hdag T ) (\E ^{(\bullet)}).
\end{equation}
La notation est justifiée par le fait que, 
d'après \cite[3.2.1]{caro-2006-surcoh-surcv}, 
ce complexe ne dépend canoniquement que de $Y$.
On dispose des propriétés ci-dessous : 

\noindent $\bullet$ Si $Y $ et $Y'$ sont deux sous-schémas de $P$, pour tout
$\E ^{(\bullet)} \in (F\text{-})\smash{\underrightarrow{LD}} ^{\mathrm{b}} _{\Q ,\mathrm{qc}}
(\smash{\widehat{\D}} _{\PP} ^{(\bullet)})$,
on dispose de l'isomorphisme canonique :
\begin{equation}
  \label{gammayY'}
  \R\underline{\Gamma} ^\dag _{Y} \circ \R\underline{\Gamma} ^\dag _{Y'} (\E ^{(\bullet)})
  \riso
  \R\underline{\Gamma} ^\dag _{Y \cap Y'} (\E ^{(\bullet)}).
\end{equation}

\noindent $\bullet$ Soit $\E ^{(\bullet)} \in (F\text{-})\smash{\underrightarrow{LD}} ^{\mathrm{b}} _{\Q ,\mathrm{qc}}
(\smash{\widehat{\D}} _{\PP} ^{(\bullet)})$.
Si $Y '$ est un ouvert (resp. un fermé) de $Y$, on dispose du morphisme canonique
$\R\underline{\Gamma} ^\dag _{Y} (\E ^{(\bullet)}) \rightarrow \R\underline{\Gamma} ^\dag _{Y'} (\E ^{(\bullet)})$
(resp. $\R\underline{\Gamma} ^\dag _{Y'} (\E ^{(\bullet)}) \rightarrow \R\underline{\Gamma} ^\dag _{Y} (\E ^{(\bullet)})$).
Si $Y'$ est un fermé de $Y$, on bénéficie aussi du triangle distingué de localisation :
$\R\underline{\Gamma} ^\dag _{Y'} (\E ^{(\bullet)}) \rightarrow \R\underline{\Gamma} ^\dag _{Y} (\E ^{(\bullet)})
\rightarrow \R\underline{\Gamma} ^\dag _{Y \setminus Y'} (\E ^{(\bullet)}) \rightarrow +1$.
\end{vide}

\subsection{Surcohérence, surholonomie, holonomie sur des cadres ou (couples de) variétés}

Rappelons d'abord, avec quelques ajouts utiles dans la suite de ce papier, les définitions et notations de \cite[1.2 et 4.1]{caro-image-directe}:
\begin{defi}
\label{defi-cadre}
On définit la catégorie des cadres (resp. cadres lisses, resp. cadres lisses en dehors du diviseur)
 de la manière suivante:
 
\begin{enumerate}
 \item Un {\og cadre $(\PP, T,X,Y)$\fg} est la donnée d'un $\V$-schéma formel séparé et lisse $\PP$, 
d'un diviseur $T$ de $P$, d'un sous-schéma fermé $X$ de $P$ tels que
$Y=X \setminus T$.
 Un cadre $(\PP, T,X,Y)$
est {\og lisse\fg} 
(resp. {\og lisse en dehors du diviseur\fg})
si $X$ (resp. $Y$) est lisse.

\item Un morphisme de cadres (resp. cadres lisses, resp. cadres lisses en dehors du diviseur) 
$\theta \colon   (\PP', T',X',Y')\to (\PP, T,X,Y)$
est la donnée d'un morphisme $f\colon   \PP' \to \PP$ tel que $f (X') \subset X$ et $f (Y') \subset Y$.
Si $a\colon  X' \to X$, $b \colon   Y' \to Y$ sont les morphismes induits par $f$, le morphisme $\theta$ se note aussi $(f,a,b)$.
\end{enumerate}
\end{defi}

 \begin{defi}
\label{defi-(d)plong}
On définit la catégorie des couples de $k$-variétés $d$-plongeables de la manière suivante: 
\begin{itemize}
\item Soient $X$ une $k$-variété et $Y$ un ouvert de $X$.  
Le couple $(Y,X)$ est un {\og couple de $k$-variétés $d$-plongeables\fg} 
s'il existe un cadre de la forme $(\PP, T,X,Y)$ 
(voir les conventions de \ref{defi-cadre}).

\item Soient $(Y',X')$ et $(Y,X)$ deux couples de $k$-variétés $d$-plongeables.
Un morphisme $(b,a)\colon (Y',X')\to (Y,X)$ de couples de $k$-variétés $d$-plongeables est un morphisme de variétés
$a \colon X'\to X$ induisant la factorisation $b\colon Y' \to Y$.
On pourra noter abusivement $a$ pour $(b,a)$.

\end{itemize}
\end{defi}

\begin{nota}
\label{nota-surhol-T}

Soit $(\PP, T,X,Y)$ un cadre.

$\bullet$
Selon \cite{caro_surholonome},
on note $(F\text{-})D ^\mathrm{b} _\mathrm{hol} (\D ^\dag _{\PP,\Q})$
(resp. $(F\text{-})D ^\mathrm{b} _\mathrm{surhol} (\D ^\dag _{\PP,\Q})$) 
la sous-catégorie pleine de 
$(F\text{-})D ^\mathrm{b}  (\D ^\dag _{\PP,\Q})$ des 
 $(F\text{-})$complexes holonomes (resp. surholonomes).

$\bullet$ On note $(F\text{-})D ^\mathrm{b} _\mathrm{hol}  (\PP, T, X/K)$
(resp. $(F\text{-})D ^\mathrm{b} _\mathrm{surhol}  (\PP, T, X/K)$) 
 la sous-catégorie pleine de 
$(F\text{-})D ^\mathrm{b} _\mathrm{coh} (\D ^\dag _{\PP} (\hdag T) _\Q)$ des $(F\text{-})$complexes $\FF$ à support dans $X$ tels que 
  $oub _{T} (\FF)\in D ^\mathrm{b} _\mathrm{hol} (\D ^\dag _{\PP,\Q})$
  (resp. $oub _{T} (\FF)\in D ^\mathrm{b} _\mathrm{surhol} (\D ^\dag _{\PP,\Q})$).

 \end{nota}

\begin{nota}
\label{4.2.5-imagedirecte}
Soit $(Y,X)$ un couple de $k$-variétés proprement $d$-plongeables.
D'après \cite[4.2.5]{caro-image-directe},
la catégorie $(F\text{-})D ^\mathrm{b} _\mathrm{surhol}  (\PP, T, X/K)$
ne dépend pas, à isomorphisme canonique près, du choix
du $\V$-schéma formel $\mathcal{Q} $ propre et lisse,
de l'ouvert $\PP$ de $\mathcal{Q}$, 
du diviseur $T$ de $P$ et d'immersion fermée $X \hookrightarrow \PP$
vérifiant $Y = X \setminus T$.
Elle se notera alors sans ambiguïté 
$(F\text{-})D ^\mathrm{b} _\mathrm{surhol} (\D ^\dag _{(Y,X)/K}) $.
\end{nota}

\begin{nota}
\label{nota-holY}
$\bullet$
 Soit $(\PP, T,X,Y)$ un cadre tel que 
 $\PP$ soit un $\V$-schéma formel projectif et lisse.
Dans ce cas, il découle de 
\cite{caro-stab-holo} que l'on dispose de l'égalité:
  $$
  F\text{-}D ^\mathrm{b} _\mathrm{surhol}  (\PP, T, X/K)
=
  F\text{-}D ^\mathrm{b} _\mathrm{hol}  (\PP, T, X/K).$$

$\bullet$ Soit $Y$ une variété quasi-projective. 
De manière analogue au cas surholonome, 
nous avons vérifié dans \cite[2.2]{caro-stab-holo} 
que la catégorie
$  F\text{-}D ^\mathrm{b} _\mathrm{hol}  (\PP, T, X/K)$
ne dépend pas du choix du cadre
 $(\PP, T,X,Y)$ tel que 
 $\PP$ soit un $\V$-schéma formel projectif et lisse.
  On l'avait alors simplement noté
  $F\text{-}D ^\mathrm{b} _\mathrm{hol}  (Y/K)$
Pour rester cohérent avec le cas surholonome, 
nous la noterons aussi
$F \text{-}D ^\mathrm{b} _\mathrm{hol}  (\D ^\dag _{Y/K})$.
\end{nota}

\begin{nota}
\label{Qsurcoh}
Soit $(\PP, T,X,Y)$ un cadre.
D'après \cite{caro-stab-sys-ind-surcoh},
on note $\underrightarrow{LD}  ^\mathrm{b} _{\Q, \mathrm{surcoh}}
(\overset{^\mathrm{g}}{} \smash{\widehat{\D}} _{\PP} ^{(\bullet)} (T ))$ 
la sous-catégorie strictement pleine de 
$\underrightarrow{LD}  ^\mathrm{b} _{\Q, \mathrm{coh}}
(\overset{^\mathrm{g}}{} \smash{\widehat{\D}} _{\PP} ^{(\bullet)} (T ))$
des complexes $\E ^{(\bullet)}$ surcohérents, i.e., satisfaisant à la propriété suivante: 
pour tout morphisme lisse $f \colon \PP' \to \PP$, pour tout diviseur $T'$ de $P'$, on a alors
$(\hdag T') \circ f ^{(\bullet)*} (\E ^{(\bullet)} )
\in 
\underrightarrow{LD}  ^\mathrm{b} _{\Q, \mathrm{coh}}
(\overset{^\mathrm{g}}{} \smash{\widehat{\D}} _{\PP'} ^{(\bullet)} (f ^{-1}(T )))$.
Cette définition de surcohérence est analogue à celle des complexes surcohérents
de $D ^\mathrm{b} _\mathrm{surcoh} ( \smash{\D} ^\dag _{\PP} (\hdag T) _{\Q} )$
définis dans \cite{caro_surcoherent}.
D'aillleurs, on a vérifié dans \cite{caro-stab-sys-ind-surcoh} que l'on bénéficie du l'équivalence de catégories
\begin{equation}
\label{pl-fid-surcoh}
\underrightarrow{\lim} \colon 
(F\text{-})\smash{\underrightarrow{LD}} ^{\mathrm{b}} _{\Q ,\mathrm{surcoh}}
(\smash{\widehat{\D}} _{\PP} ^{(\bullet)}(T))
\cong 
(F\text{-})D ^\mathrm{b} _\mathrm{surcoh} ( \smash{\D} ^\dag _{\PP} (\hdag T) _{\Q} ).
\end{equation}
Conformément aux autres notations,
la sous-catégorie pleine de 
$(F\text{-})D ^\mathrm{b} _\mathrm{surcoh} ( \smash{\D} ^\dag _{\PP} (\hdag T) _{\Q} )$
des complexes à support dans $X$ se notera 
$(F\text{-})D ^\mathrm{b} _{\mathrm{surcoh}}  (\PP, T, X/K) $.
\end{nota}

\subsection{Isocristaux, dévissabilité en isocristaux surconvergents}

Soit $(\PP, T,X,Y)$ un cadre lisse en dehors du diviseur
(voir les conventions de \ref{defi-cadre}).

\begin{vide}
\label{nota-IsocDag}
On note $j \colon Y \subset X$ l'immersion ouverte canonique et
    $\mathrm{Isoc} ^{\dag} (\PP , T , X /K)$
la catégorie des $j ^\dag \O _{]X [ _{\PP}}$-modules cohérents munis d'une connexion surconvergente.
Cette catégorie est canoniquement isomorphe à celle des isocristaux sur $Y$ surconvergents le long de $T$ et notée
  $\mathrm{Isoc} ^{\dag} (Y ,X/K)$ (voir \cite[2.3.2 et 2.3.7]{Berig}).
\end{vide}

Rappelons les notations (voir \cite[3.1.2 et 3.5.10]{caro-pleine-fidelite}) et définitions
(voir \cite[5.4.5 et 5.4.7]{caro-pleine-fidelite}) suivantes:

\begin{defi}
  \label{defi-part-surcoh}

\begin{itemize}
\item On définit la catégorie $\mathrm{Isoc} ^{\dag \dag} (\PP, T, X/K)$ d'abord en considérant le cas où $X$ lisse qui est traité 
dans \cite{caro-construction} (voir un aperçu de la construction donné dans le paragraphe \ref{const-sp+lisse}). 
Dans le cas général, $\mathrm{Isoc} ^{\dag \dag} (\PP, T, X/K)$ est la catégorie
des $\D ^{\dag} _{\PP}(\hdag T) _{\Q}$-modules surcohérents $\E$ tels que, en posant $\U := \PP \setminus T$, 
on ait
$\E |\U \in \mathrm{Isoc} ^{\dag \dag} (\U, Y/K)$ (on n'indique pas le diviseur lorsqu'il est vide). 
Les objets de $\mathrm{Isoc} ^{\dag \dag} (\PP, T, X/K)$ sont {\og les isocristaux partiellement surcohérents sur $(\PP, T, X/K)$\fg}.
Lorsque $\PP$ est propre, on omet le qualificatif {\og partiellement\fg}. 

\item Une fois $(Y,X)/K$ fixé, la catégorie $(F\text{-})\mathrm{Isoc} ^{\dag \dag} (\PP, T, X/K)$ ne dépend, à équivalence canonique de catégories près, ni du choix de l'immersion fermée
$X \hookrightarrow \PP$ et ni de celui du diviseur $T$ de $P$ tel que $Y = X \setminus T$.
On note alors
sans ambiguïté 
$(F\text{-})\mathrm{Isoc}^{\dag \dag} (Y,X/K)$
à la place de $(F\text{-})\mathrm{Isoc} ^{\dag \dag} (\PP, T, X/K)$.
Ses objets sont les {\og $(F\text{-})$isocristaux partiellement surcohérents sur $(Y,X)/K$\fg} 
ou simplement {\og $(F\text{-})$isocristaux surcohérents sur $(Y,X)/K$\fg}.

\item Lorsque $X$ est propre, la catégorie 
$(F\text{-})\mathrm{Isoc}^{\dag \dag} (Y,X/K)$ ne dépend pas non plus du choix de la compactification propre $X$ de $Y$.
On la note alors
$(F\text{-})\mathrm{Isoc}^{\dag \dag} (Y/K)$.
Ses objets sont les {\og $(F\text{-})$isocristaux surcohérents sur $Y/K$\fg}.
\end{itemize}

\end{defi}

\begin{nota}
\label{nota-M-eq-isoc-lim}
On notera $\mathrm{Isoc} ^{(\bullet)}( \PP, T , X /K)$ la sous-catégorie strictement pleine de 
$\smash{\underrightarrow{LM}}  _{\Q, \mathrm{coh}}
(\smash{\widehat{\D}} _{\PP} ^{(\bullet)}(T))$
telle que l'équivalence de catégories de \ref{M-eq-coh-lim} se factorise en l'équivalence de catégories
\begin{equation}
\label{M-eq-isoc-lim}
\underrightarrow{\lim} 
\colon
\mathrm{Isoc} ^{(\bullet)}( \PP, T , X /K)\cong
\mathrm{Isoc} ^{\dag \dag}( \PP, T , X /K).
\end{equation}

\end{nota}

\begin{vide}
\label{vide-eq-cat-spxPT+}
D'après \cite[3.5.10 et 4.2.2]{caro-pleine-fidelite},
on dispose de l'équivalence canonique de catégories 
\begin{equation}
\label{eq-cat-spxPT+}
\sp _{X \hookrightarrow \PP, T,+} \colon 
\mathrm{Isoc} ^{\dag}( \PP, T , X /K)
\cong 
\mathrm{Isoc} ^{\dag \dag}( \PP, T , X /K).
\end{equation}
L'équivalence de catégories de \ref{eq-cat-spxPT+}
ne dépend canoniquement pas des choix faits et pourra simplement être notée 
$\sp _{(Y,X),+}\colon
\mathrm{Isoc}^{ \dag} (Y,X/K) \cong 
\mathrm{Isoc}^{\dag \dag} (Y,X/K)$
et, lorsque $X$ est propre,
$\sp _{Y,+}\colon 
\mathrm{Isoc}^{ \dag} (Y/K) \cong 
\mathrm{Isoc}^{\dag \dag} (Y/K)$.
\end{vide}

\begin{vide}
\label{corr-424525}
Soit 
$\theta = (f,a,b)\colon   (\PP', T',X',Y')\to (\PP, T,X,Y)$
un morphisme de cadres lisses en dehors du diviseur.
D'après \cite[3.1.7 et 3.5.10]{caro-pleine-fidelite},
on dispose des foncteurs exacts
\begin{align}
\label{stabIsoc*inv-i} 
(\hdag T') \circ \R \underline{\Gamma} ^\dag _{X'} \circ f ^! [-d _{Y'/Y}] \colon ,&(F\text{-})\mathrm{Isoc} ^{\dag\dag} (\PP, T, X/K) \to (F\text{-})\mathrm{Isoc} ^{\dag\dag} (\PP', T', X'/K),\\
\label{stabIsoc*inv-ii}
\DD _{T'}\circ \R \underline{\Gamma} ^\dag _{X'}\circ (\hdag T') \circ f ^! [-d _{Y'/Y}] \circ \DD _{T} \colon &
(F\text{-})\mathrm{Isoc} ^{\dag\dag} (\PP, T, X/K)\to (F\text{-})\mathrm{Isoc} ^{\dag\dag} (\PP', T', X'/K)).
\end{align}

$ \bullet$ Comme nous ne travaillions qu'avec des isocristaux, afin d'alléger les notations, 
cela nous avait conduit à poser
$\theta ^{!}\,:=\,(\hdag T') \circ \R \underline{\Gamma} ^\dag _{X'} \circ f ^! [-d _{Y'/Y}]$
et 
$\theta ^{+}\,:=\,\DD _{T'}\circ \R \underline{\Gamma} ^\dag _{X'}\circ (\hdag T') \circ f ^! [-d _{Y'/Y}] \circ \DD _{T}$.
De plus, par \cite[5.2.5]{caro-pleine-fidelite}, on avait vérifié que
les foncteurs $ \theta ^{+}$ et $ \theta ^{!} \colon 
\mathrm{Isoc} ^{\dag \dag} (\PP, T, X/K) \to \mathrm{Isoc} ^{\dag \dag} (\PP', T', X'/K) $
sont canoniquement isomorphes.
Par contre, d'après \cite[5.6]{Abe-Frob-Poincare-dual}, 
pour obtenir un isomorphisme compatible à Frobenius entre
$ \theta ^{+}$ et $ \theta ^{!}$, il faut rajouter un twist.  
En général, lorsque l'on travaille avec des $\D$-modules arithmétiques, 
on note plutôt 
$\theta ^{!}\,:=\,(\hdag T') \circ \R \underline{\Gamma} ^\dag _{X'} \circ f ^!$
et
$\theta ^{+}\,:=\,\DD _{T'}\circ \R \underline{\Gamma} ^\dag _{X'}\circ (\hdag T') \circ f ^!\circ \DD _{T}$.
Pour éviter toute ambiguïté, dans le cas des isocristaux, 
nous noterons 
$\theta ^{*}$ à la place de
$(\hdag T') \circ \R \underline{\Gamma} ^\dag _{X'} \circ f ^! [-d _{Y'/Y}]$ ou de 
$\DD _{T'}\circ \R \underline{\Gamma} ^\dag _{X'}\circ (\hdag T') \circ f ^! [-d _{Y'/Y}] \circ \DD _{T}$.
Si on tient compte des structures de Frobenius,
$\theta ^{*}$ désigne alors le foncteur
$(\hdag T') \circ \R \underline{\Gamma} ^\dag _{X'} \circ f ^! [-d _{Y'/Y}]$.

$\bullet$ Avec \ref{pl-fid-surcoh}, il vient 
$\mathrm{Isoc} ^{(\bullet)}( \PP, T , X /K) \subset
\smash{\underrightarrow{LD}} ^{\mathrm{b}} _{\Q ,\mathrm{surcoh}}
(\smash{\widehat{\D}} _{\PP} ^{(\bullet)}(T))$.
On obtient de plus le foncteur exact:
\begin{equation}
\label{stabIsoc*inv-i-bullet} 
\theta ^{*}:= (\hdag T') \circ \R \underline{\Gamma} ^\dag _{X'} \circ f ^! [-d _{Y'/Y}]
\colon 
\mathrm{Isoc} ^{(\bullet)} (\PP, T, X/K) 
\to 
\mathrm{Isoc} ^{(\bullet)} (\PP', T', X'/K).
\end{equation}

$\bullet$ On dispose du
foncteur canonique
$\theta ^*\colon  \mathrm{Isoc} ^{\dag}( \PP, T, X/K)\rightarrow
 \mathrm{Isoc} ^{\dag }( \PP', T', X'/K) $
 (voir \cite[2.3.2.2]{Berig}).

$ \bullet$  D'après le théorème \cite[4.2.4]{caro-pleine-fidelite}, 
avec ces notations,
on dispose alors de l'isomorphisme canonique
\begin{equation}
\label{iso-theta*com-sp}
\sp _{X '\hookrightarrow \PP', T',+} \circ \theta ^{*} 
\riso 
\theta ^{*}\circ  \sp _{X \hookrightarrow \PP, T,+}.
\end{equation}

\end{vide}

\begin{nota}
\label{rappel-dev-coh}
Soit $(\PP, T,X,Y)$ un cadre.
Lorsque $Y$ est lisse, 
grâce à \cite{caro-pleine-fidelite}, on a pu définir la catégorie des isocristaux surcohérents
$\mathrm{Isoc} ^{\dag \dag} ( \PP, T ,X /K)$ dans le cas où $\PP$ est seulement séparé et lisse, 
ce qui étend les travaux \cite{caro_devissge_surcoh} et \cite{caro-2006-surcoh-surcv}
du cas où $\PP$ est propre au cas où  $\PP$ est séparé et lisse.
On en avait déduit que la notion de 
 dévissabilité en isocristaux surconvergents de 
\cite[8.1.1]{caro-2006-surcoh-surcv} s'étendait 
du cas où $\PP$ est propre au cas où  $\PP$ est séparé et lisse.
Via la définition \cite[6.2.2]{caro-pleine-fidelite}, on avait alors 
noté  
$D ^\mathrm{b} _{\textrm{dév}}
(\smash{\D} ^\dag _{\PP} (\hdag T) _\Q)$
la sous-catégorie pleine de 
$D ^\mathrm{b} _{\textrm{coh}}
(\smash{\D} ^\dag _{\PP} (\hdag T) _\Q)$
des complexes dévissables en isocristaux surconvergents.
On note
$(F\text{-})D ^\mathrm{b} _{\text{dév}}  (\PP, T, X/K) $
la sous-catégorie pleine de 
$(F\text{-})D ^\mathrm{b} _{\text{dév}} (\smash{\D} ^\dag _{\PP} (\hdag T) _\Q)$
des complexes qui sont à support dans $X$.
\end{nota}

\begin{rema}
\label{rema-precision-in-isoc}
Soit $(\PP, T,X,Y)$ un cadre lisse en dehors du diviseur.
Soient $\E ^{(\bullet)}\in \underrightarrow{LD}  ^\mathrm{b} _{\Q, \mathrm{coh}}
(\overset{^\mathrm{g}}{} \smash{\widehat{\D}} _{\PP} ^{(\bullet)} (T ))$
et $\E := \underrightarrow{\lim} ~\E ^{(\bullet)} $.
Les notations sur l'hypothèse $2$ de la définition \cite[6.2.2]{caro-pleine-fidelite} 
(ou bien aussi de \cite[3.2.2]{caro-2006-surcoh-surcv} ) sont abusives: 
pour être précis et exact,
{\og les espaces de cohomologie de 
$\R\underline{\Gamma} ^\dag _{X} (\hdag T ) (\E )$ appartiennent
à $\mathrm{Isoc} ^{\dag \dag} ( \PP, T ,X /K)$\fg} signifie que
$\R\underline{\Gamma} ^\dag _{X} (\hdag T ) (\E ^{(\bullet)}) 
\in 
\underrightarrow{LD}  ^\mathrm{b} _{\Q, \mathrm{coh}}
(\overset{^\mathrm{g}}{} \smash{\widehat{\D}} _{\PP} ^{(\bullet)} (T ))$
et que, pour tout $j \in \Z$, on ait 
$\mathcal{H} ^j \R\underline{\Gamma} ^\dag _{X} (\hdag T ) (\E )
:=
\mathcal{H} ^j  \underrightarrow{\lim}
\R\underline{\Gamma} ^\dag _{X} (\hdag T ) (\E ^{(\bullet)})
\in 
\mathrm{Isoc} ^{\dag \dag} ( \PP, T ,X /K)$.

\end{rema}

\begin{vide}
Soit $(\PP, T,X,Y)$ un cadre.
On bénéficie des inclusions 
(voir \cite[6.2.3]{caro-pleine-fidelite} pour la dernière):
\begin{equation}
\label{subset-surholdev}
D ^\mathrm{b} _{\mathrm{surhol}} 
 (\PP, T, X/K)
\subset
D ^\mathrm{b} _{\mathrm{surcoh}} 
(\PP, T, X/K)
\subset
D ^\mathrm{b} _{\textrm{\rm dév}}
(\PP, T, X/K).
\end{equation}
D'après \cite[6.2.4]{caro-pleine-fidelite},
on dispose avec une structure de Frobenius
 des égalités
 \begin{equation}
\label{=Fsurholdev}
F\text{-}D ^\mathrm{b} _{\textrm{\rm dév}} (\PP, T, X/K)
=F\text{-}D ^\mathrm{b} _{\mathrm{surcoh}} 
(\PP, T, X/K)
=F\text{-}D ^\mathrm{b} _{\mathrm{surhol}} 
(\PP, T, X/K).
\end{equation}
Rappelons que ces égalités
sont des conséquences du théorème de réduction semi-stable de Kedlaya
(voir \cite{kedlaya-semistableI}, \cite{kedlaya-semistableII}, \cite{kedlaya-semistableIII},
\cite{kedlaya-semistableIV} ou
\cite{tsumono} pour la version concernant les $F$-isocristaux surconvergents unités)
précédemment conjecturée par Shiho (voir \cite{Shiho-log-isocI}, \cite{Shiho-log-isocII}).

\end{vide}

\section{Produits tensoriels en théorie des $\D$-modules arithmétiques}

\subsection{Produits tensoriels internes: définition et propriétés de commutation}
Soient $\PP$ un $\V$-schéma formel lisse, $T $ un diviseur de $P$.

\begin{vide}
Pour tous 
$\E\in D _{\mathrm{qc}} ^\mathrm{b}
(\overset{^\mathrm{g}}{} \smash{\widehat{\D}} _{\PP} ^{(m)} (T ))$
et $\M \in D _{\mathrm{qc}} ^\mathrm{b}
(\overset{^\mathrm{*}}{}\smash{\widehat{\D}} _{\PP} ^{(m)} (T ) )$,
on pose
\begin{gather} 
\notag
\M \smash{\widehat{\otimes}} ^\L _{\smash{\widehat{\B}} _{\PP} ^{(m)} (T )} \E :=
\R \underset{\underset{i}{\longleftarrow}}{\lim}\, ( \M _i \otimes ^\L  _{\smash{\B} _{P _i} ^{(m)} (T )} \E _i).
\end{gather}
\end{vide}

\begin{vide}
Pour tous 
$\E ^{(\bullet)} \in \underrightarrow{LD}  ^\mathrm{b} _{\Q, \mathrm{qc}}
(\overset{^\mathrm{g}}{} \smash{\widehat{\D}} _{\PP} ^{(\bullet)} (T ))$,
$\M ^{(\bullet)} \in \underrightarrow{LD}  ^\mathrm{b} _{\Q, \mathrm{qc}}
( \overset{^\mathrm{*}}{}\smash{\widehat{\D}} _{\PP} ^{(\bullet)} (T ) )$.
On définit les bifoncteurs produits tensoriels 
\begin{align}
\label{def-otimes-qc1}
 \smash{\overset{\L}{\otimes}}^{\dag} _{\O _{\PP } ( \hdag T ) _{\Q}}
  \colon
\underrightarrow{LD}  ^\mathrm{b} _{\Q, \mathrm{qc}}
( \overset{^*}{} \smash{\widehat{\D}} _{\PP} ^{(\bullet)} (T ) )
\times 
\underrightarrow{LD}  ^\mathrm{b} _{\Q, \mathrm{qc}}
(\overset{^\mathrm{g}}{} \smash{\widehat{\D}} _{\PP} ^{(\bullet)} (T ))
&
\to 
\underrightarrow{LD}  ^\mathrm{b} _{\Q, \mathrm{qc}}
( \overset{^*}{} \smash{\widehat{\D}} _{\PP} ^{(\bullet)} (T )  ).
\end{align}
en posant 
\begin{equation}
\notag
\M ^{(\bullet)}
\smash{\overset{\L}{\otimes}}   ^{\dag}
_{\O  _{\PP} (\hdag T) _\Q}\E ^{(\bullet)}
:=
(\M ^{(m)}  \smash{\widehat{\otimes}}
^\L _{\widehat{\B} ^{(m)}  _{\PP} ( T) }
\E ^{(m)}) _{m\in \N}
\end{equation}

\end{vide}

\begin{vide}
\label{def-otimes-coh}
Soient
$\E ^{(\bullet)}  
\in \underrightarrow{LD}  ^\mathrm{b} _{\Q, \mathrm{coh}}
(\overset{^\mathrm{g}}{} \smash{\widehat{\D}} _{\PP} ^{(\bullet)} (T ))$,
$\M ^{(\bullet)}
\in \underrightarrow{LD}  ^\mathrm{b} _{\Q, \mathrm{coh}}
(\overset{^\mathrm{*}}{} \smash{\widehat{\D}} _{\PP} ^{(\bullet)} (T ))$.
Posons 
$\E:= \underrightarrow{\lim}\, \E ^{(\bullet)}
\in D ^\mathrm{b} _\mathrm{coh} (\overset{^\mathrm{g}}{} \smash{\D} ^\dag _{\PP} (\hdag T) _{\Q} )$
et
$\M := \underrightarrow{\lim} \, \M ^{(\bullet)} 
\in D ^\mathrm{b} _\mathrm{coh} (\overset{^\mathrm{*}}{} \smash{\D} ^\dag _{\PP} (\hdag T) _{\Q} )$.
Avec les notations de \ref{fct-qcoh2coh}, 
on obtient le bifoncteur
$ \smash{\overset{\L}{\otimes}}^{\dag} _{\O _{\PP } ( \hdag T ) _{\Q}}
  := 
\mathrm{Coh} _{T , T }  
 ( \smash{\overset{\L}{\otimes}}^{\dag} _{\O _{\PP } ( \hdag T ) _{\Q}}
)$
 de la forme:
\begin{align}
\label{def-otimes-coh1}
 \smash{\overset{\L}{\otimes}}^{\dag} _{\O _{\PP } ( \hdag T ) _{\Q}}
  \colon
 D ^\mathrm{b} _\mathrm{coh} (\overset{^\mathrm{*}}{} \smash{\D} ^\dag _{\PP} (\hdag T) _{\Q} )
\times 
D ^\mathrm{b} _\mathrm{coh} (\overset{^\mathrm{g}}{} \smash{\D} ^\dag _{\PP} (\hdag T) _{\Q} )
&
\to D ^\mathrm{b} (\overset{^\mathrm{*}}{} \smash{\D} ^\dag _{\PP} (\hdag T) _{\Q} ).
\end{align}
Par définition, on dispose donc des isomorphismes fonctoriels
\begin{gather}
  \label{coh-otimesO}
    \M
 \smash{\overset{\L}{\otimes}}^{\dag} _{\O _{\PP } ( \hdag T ) _{\Q}}
  \E
\riso
  \underrightarrow{\lim} \,  \M ^{(\bullet)}
 \smash{\overset{\L}{\otimes}}^{\dag} _{\O _{\PP } ( \hdag T ) _{\Q}}
  \E ^{(\bullet)}.
\end{gather}

\end{vide}

\begin{lemm}
\label{hdagTT'}
  Soient $T'$ un second diviseur de $P$,
  $\G ^{(\bullet)}
\in \underrightarrow{LD}  ^\mathrm{b} _{\Q, \mathrm{qc}}
(\overset{^\mathrm{*}}{} \smash{\widehat{\D}} _{\PP} ^{(\bullet)} (T ))$.
Le morphisme canonique 
$\G ^{(\bullet)} (\hdag T') \rightarrow \G ^{(\bullet)} (\hdag T \cup T')$
est alors un isomorphisme
de
$\underrightarrow{LD}  ^\mathrm{b} _{\Q, \mathrm{qc}} (
\smash{\widehat{\D}} _{\PP} ^{(\bullet)} )$.
\end{lemm}

\begin{proof}
D'après 
\cite[3.2.6]{caro-stab-sys-ind-surcoh}, en omettant d'indiquer le foncteur oubli, 
on dispose de l'isomorphisme
$\G ^{(\bullet)} \riso \G ^{(\bullet)} (\hdag T) $.
De plus, par \cite[2.2.14]{caro_surcoherent},
$(\hdag T ') \circ (\hdag T ) \riso (\hdag T \cup T') $.
On en tire l'isomorphisme voulu
$\G ^{(\bullet)} (\hdag T') \riso \G ^{(\bullet)} (\hdag T \cup T')$.
\end{proof}

\begin{prop}
\label{otimesTT'}
  Soient $T'$ un second diviseur de $P$,
  $\G ^{(\bullet)}
\in \underrightarrow{LD}  ^\mathrm{b} _{\Q, \mathrm{qc}}
(\overset{^\mathrm{*}}{} \smash{\widehat{\D}} _{\PP} ^{(\bullet)} (T ))$,
$\E ^{(\bullet)}
\in \underrightarrow{LD}  ^\mathrm{b} _{\Q, \mathrm{qc}}
(\overset{^\mathrm{g}}{} \smash{\widehat{\D}} _{\PP} ^{(\bullet)} (T '))$.
Le morphisme canonique :
\begin{equation}
\label{otimesTT'-iso}
  \G ^{(\bullet)}
\smash{\overset{\L}{\otimes}}   ^{\dag}
_{\O  _{\PP,\Q}}\E ^{(\bullet)}
\to
\G ^{(\bullet)} (\hdag T \cup T')
\smash{\overset{\L}{\otimes}}   ^{\dag} _{\O  _{\PP} (\hdag T \cup T') _\Q}
\E ^{(\bullet)} (\hdag T \cup T')
\end{equation}
est un isomorphisme de $\underrightarrow{LD}  ^\mathrm{b} _{\Q, \mathrm{qc}}
(\overset{^\mathrm{*}}{} \smash{\widehat{\D}} _{\PP} ^{(\bullet)})$. 
\end{prop}

\begin{proof}
D'après \cite[1.1.8]{caro_courbe-nouveau},
le morphisme 
$
\widehat{\B} ^{(\bullet)} _{\PP} ( T \cup T')  \smash{\widehat{\otimes}} ^\L
_{\O _{\PP} } \G ^{(\bullet)}
\to 
\widehat{\D} ^{(\bullet)} _{\PP} ( T \cup T')  \smash{\widehat{\otimes}} ^\L
_{\widehat{\D} ^{(\bullet)} _{\PP}} \G ^{(\bullet)}$
est un isomorphisme. 
De même, pour $\E ^{(\bullet)}$.
Via les propriétés d'associativité du produit tensoriel, 
on en déduit l'isomorphisme canonique:
$(\hdag T \cup T') \left (\G ^{(\bullet)}
\smash{\overset{\L}{\otimes}}   ^{\dag}
_{\O  _{\PP,\Q}}\E ^{(\bullet)} \right )
\riso 
\G ^{(\bullet)} (\hdag T \cup T')
\smash{\overset{\L}{\otimes}}   ^{\dag} _{\O  _{\PP} (\hdag T \cup T') _\Q}
\E ^{(\bullet)} (\hdag T \cup T')$.
Grâce au triangle de localisation par rapport à $T \cup T'$, 
il s'agit donc d'établir 
$\R \underline{\Gamma} ^{\dag} _{T \cup T'}~
(\G ^{(\bullet)}
\smash{\overset{\L}{\otimes}}   ^{\dag}
_{\O  _{\PP,\Q}}\E ^{(\bullet)}) =0$. 
Or, 
$\R  \underline{\Gamma} ^{\dag} _{T}~
(\G ^{(\bullet)}
\smash{\overset{\L}{\otimes}}   ^{\dag}
_{\O  _{\PP,\Q}}\E ^{(\bullet)})
\riso 
\R \underline{\Gamma} ^{\dag} _{T}~
(\G ^{(\bullet)})
\smash{\overset{\L}{\otimes}}   ^{\dag}
_{\O  _{\PP,\Q}}\E ^{(\bullet)}
\riso 0$.
De même, 
$\R \underline{\Gamma} ^{\dag} _{T'}~
(\G ^{(\bullet)}
\smash{\overset{\L}{\otimes}}   ^{\dag}
_{\O  _{\PP,\Q}}\E ^{(\bullet)})
\riso 0$.
Le triangle distingué de Mayer-Vietoris, nous permet de conclure. 
\end{proof}

\begin{vide}
\label{2.1.4-caro-surcoh}
Soient $f \colon \PP' \to \PP$ un morphisme de $\V$-schémas formels lisses, 
$T$ et $T'$ des diviseurs respectifs de $P$ et $P'$ tels que
$f ( P '\setminus T' ) \subset P \setminus T$.
Soient $\E ^{(\bullet)}
\in \underrightarrow{LD}  ^\mathrm{b} _{\Q, \mathrm{coh}}
(\overset{^\mathrm{g}}{} \smash{\widehat{\D}} _{\PP} ^{(\bullet)} (T ))$
et
$\E ^{\prime (\bullet)}
\in \underrightarrow{LD}  ^\mathrm{b} _{\Q, \mathrm{coh}}
(\overset{^\mathrm{g}}{} \smash{\widehat{\D}} _{\PP'} ^{(\bullet)} (T '))$.
D'après \cite[2.1.4]{caro_surcoherent} (et avec \ref{f+!-indTT'}
et \ref{otimesTT'}, les singularités surconvergentes le long de $T$ ou de $T'$ ne sont pas gênantes),
on dispose de l'isomorphisme :
\begin{equation}
\label{2.1.4-caro-surcoh-iso}
f ^{(\bullet)} _{T,T',+} 
(\E ^{\prime (\bullet)})
\smash{\overset{\L}{\otimes}}   ^{\dag} _{\O _{\PP} (\hdag T) _{\Q}}
 \E ^{(\bullet)}
\liso
f ^{(\bullet)} _{T,T',+} \left ( \E ^{\prime (\bullet)}
\smash{\overset{\L}{\otimes}}   ^{\dag} _{\O _{\PP '} (\hdag T' ) _{\Q}}
\L f  ^{(\bullet)*}  _{T',T} (\E ^{(\bullet)}) \right).
\end{equation}
\end{vide}

\begin{lemm}
\label{u!u+=id}
Avec les notations de \ref{2.1.4-caro-surcoh},
on suppose que $f$ est une immersion fermée. 
On dispose alors, pour tous 
$\FF ^{(\bullet)}  , \, \G ^{(\bullet)}
\in \underrightarrow{LD}  ^\mathrm{b} _{\Q, \mathrm{coh}}
(\overset{^\mathrm{g}}{} \smash{\widehat{\D}} _{\PP'} ^{(\bullet)} (T'))$,
de l'isomorphisme canonique
\begin{equation}
\label{u!u+=id-cons}
f ^{(\bullet)} _{T,T',+}  (\FF ^{(\bullet)})  
\smash{\overset{\L}{\otimes}} ^\dag _{\O _{\PP} (\hdag T) _{\Q}} 
f ^{(\bullet)} _{T,T',+}  (\G ^{(\bullet)}) [d _{P'/P}]
\riso
f ^{(\bullet)} _{T,T',+}  (\FF ^{(\bullet)} \smash{\overset{\L}{\otimes}} ^\dag _{\O _{\PP' } (\hdag T') _{\Q}}\G ^{(\bullet)}).
\end{equation}
\end{lemm}

\begin{proof}
On bénéficie de l'isomorphisme canonique:
$$f ^{(\bullet)} _{T,T',+}  (\FF ^{(\bullet)}) \smash{\overset{\L}{\otimes}} ^\dag _{\O _{\PP} (\hdag T) _{\Q}}  f ^{(\bullet)} _{T,T',+} (\G ^{(\bullet)}) [d _{P'/P}]
\underset{\ref{2.1.4-caro-surcoh-iso}}{\riso}
f ^{(\bullet)} _{T,T',+}  (\FF ^{(\bullet)} 
\smash{\overset{\L}{\otimes}} ^\dag _{\O _{\PP' } (\hdag T') _{\Q}}f ^{(\bullet)!} _{T',T}\circ f ^{(\bullet)} _{T,T',+}(\G ^{(\bullet)})).$$
Or, comme $f$ est une immersion fermée, 
on déduit de  \cite[5.3.5.1]{caro-stab-sys-ind-surcoh} 
l'isomorphisme canonique 
$ \G ^{(\bullet)}
\riso f ^{(\bullet)!} _{T',T}\circ f ^{(\bullet)} _{T,T',+}(\G ^{(\bullet)})$.
D'où le résultat.
\end{proof}

\begin{prop}
\label{otimes-ind-TT'}
Soit $D \subset T$ un second diviseur de $P$.
Soient $\E, \FF \in 
D ^\mathrm{b} _\mathrm{coh} (\overset{^\mathrm{g}}{} \smash{\D} ^\dag _{\PP} (\hdag T) _{\Q} ) \cap 
D ^\mathrm{b} _\mathrm{coh} (\overset{^\mathrm{g}}{} \smash{\D} ^\dag _{\PP} (\hdag D) _{\Q} )$.
On dispose de l'isomorphisme canonique 
\begin{equation}
\label{otimesD2otimesTar}
\E
 \smash{\overset{\L}{\otimes}}^{\dag} _{\O _{\PP } ( \hdag D ) _{\Q}}
  \FF
\riso 
  \E
 \smash{\overset{\L}{\otimes}}^{\dag} _{\O _{\PP } ( \hdag T ) _{\Q}}
  \FF.
\end{equation}
\end{prop}

\begin{proof}
D'après \ref{QcohTD},
il existe
$\E ^{(\bullet)}  , \, \FF ^{(\bullet)}
\in \underrightarrow{LD}  ^\mathrm{b} _{\Q, \mathrm{coh}}
(\overset{^\mathrm{g}}{} \smash{\widehat{\D}} _{\PP} ^{(\bullet)} (D ))
\cap 
\underrightarrow{LD}  ^\mathrm{b} _{\Q, \mathrm{coh}}
(\overset{^\mathrm{g}}{} \smash{\widehat{\D}} _{\PP} ^{(\bullet)} (T ))$
tels que
$\E  \riso \underrightarrow{\lim}\, \E ^{(\bullet)} $,
$\FF  \riso \underrightarrow{\lim} \, \FF ^{(\bullet)} $.
Avec \ref{rema-fct-qcoh2coh-coh}, on vérifie alors que l'on dispose des 
isomorphismes canoniques
\begin{equation}
\label{otimesD2otimesTar1}
\E ^{(\bullet)}
 \smash{\overset{\L}{\otimes}}^{\dag} _{\O _{\PP } ( \hdag D ) _{\Q}}
  \FF ^{(\bullet)}
\riso
(\hdag T,D) (\E ^{(\bullet)} )
 \smash{\overset{\L}{\otimes}}^{\dag} _{\O _{\PP } ( \hdag T ) _{\Q}}
(\hdag T,D) (\FF ^{(\bullet)} )
\riso
\E ^{(\bullet)} 
 \smash{\overset{\L}{\otimes}}^{\dag} _{\O _{\PP } ( \hdag T ) _{\Q}}
\FF ^{(\bullet)},
\end{equation}
dont l'image par le foncteur $ \underrightarrow{\lim} \,$ 
donne l'isomorphisme \ref{otimesD2otimesTar} voulu
\end{proof}

Concluons la section par le lemme ci-dessous.
\begin{prop}
\label{fg!tau}
  Soient $f$, $g\colon \PP '\rightarrow \PP$ deux morphismes de $\V$-schémas formels lisses dont les morphismes induits
  sur les fibres spéciales coïncident, soit $a \colon \PP '' \to \PP'$ un morphisme de $\V$-schémas formels lisses, 
  $D$ un diviseur de $P$ tel que $D':=f ^{-1} (D)$ soit un diviseur de $P'$ et $D'':= a ^{-1} (D')$ soit un diviseur de $P''$.
Soient
$\E ^{(\bullet)}, \E ^{\prime (\bullet)} \in \underrightarrow{LD}  ^\mathrm{b} _{\Q, \mathrm{qc}}
( \smash{\widehat{\D}} _{\PP} ^{(\bullet)} (D ))$.
\begin{enumerate}
\item On bénéficie de l'isomorphisme canonique commutant à Frobenius
dans 
$\underrightarrow{LD}  ^\mathrm{b} _{\Q, \mathrm{qc}}
( \smash{\widehat{\D}} _{\PP'} ^{(\bullet)} (D' ))$:
\begin{equation}
\label{fg!prodtens}
f ^{(\bullet)!} (\E^{(\bullet)} \smash{\overset{\L}{\otimes}} ^\dag _{\O _{\PP} (\hdag D) _{\Q}} \E ^{\prime (\bullet)}) [d _{P'/P}]
\riso
f ^{(\bullet)!} (\E^{(\bullet)})  \smash{\overset{\L}{\otimes}} ^\dag _{\O _{\PP'} (\hdag D') _{ \Q}} f ^{(\bullet)!} (\E ^{\prime (\bullet)}).
\end{equation}
\item Les isomorphismes de la forme \ref{fg!prodtens}
sont transitifs, i.e., le diagramme suivant 
\small
\begin{equation}
\label{fg!prodtens-trans}
\xymatrix @R=0,3cm @C=0,6cm{
{a ^{(\bullet)!} \circ  f ^{(\bullet)!} (\E^{(\bullet)} 
\smash{\overset{\L}{\otimes}} ^\dag 
\E ^{\prime (\bullet)}) [d _{P''/P}]
} 
\ar[r] ^-{\ref{fg!prodtens}}
\ar[d] ^-{\sim}
& 
{a ^{(\bullet)!} (f ^{(\bullet)!} (\E^{(\bullet)})  \smash{\overset{\L}{\otimes}} ^\dag 
f ^{(\bullet)!} (\E ^{\prime (\bullet)}))[d _{P'/P}]} 
\ar[r] ^-{\ref{fg!prodtens}}
& 
{a ^{(\bullet)!} \circ   f ^{(\bullet)!} (\E^{(\bullet)})  \smash{\overset{\L}{\otimes}} ^\dag 
a ^{(\bullet)!} \circ  f ^{(\bullet)!} (\E ^{\prime (\bullet)})} 
\ar[d] ^-{\sim}
\\ 
{(f  \circ   a )^{(\bullet)!} (\E^{(\bullet)} 
\smash{\overset{\L}{\otimes}} ^\dag 
\E ^{\prime (\bullet)}) [d _{P''/P}]} 
\ar[rr] ^-{\ref{fg!prodtens}}
& 
{ } 
& 
{(f  \circ   a )^{(\bullet)!} (\E^{(\bullet)})  \smash{\overset{\L}{\otimes}} ^\dag 
(f  \circ   a )^{(\bullet)!} (\E ^{\prime (\bullet)})} 
}
\end{equation}
\normalsize
est commutatif.
\item On dispose de plus du diagramme commutatif :
\begin{equation}
  \label{fg!prodtens-square}
  \xymatrix @R=0,3cm {
  {f ^{(\bullet)!} (\E^{(\bullet)} \smash{\overset{\L}{\otimes}} ^\dag _{\O _{\PP} (\hdag D) _{ \Q}} \E ^{\prime (\bullet)}) [d _{P'/P}]}
  \ar[r] _-\sim ^-{\ref{fg!prodtens}}
  \ar[d] _-\sim ^{\tau _{g,f}}
  &
  {f ^{(\bullet)!} (\E^{(\bullet)})  \smash{\overset{\L}{\otimes}} ^\dag _{\O _{\PP'} (\hdag D') _{ \Q}} f ^{(\bullet)!} (\E ^{\prime (\bullet)}) }
  \ar[d] _-\sim ^{\tau _{g,f} \otimes \tau _{g,f}}
  \\
  {g ^{(\bullet)!} (\E^{(\bullet)} \smash{\overset{\L}{\otimes}} ^\dag _{\O _{\PP} (\hdag D) _{ \Q}} \E ^{\prime (\bullet)}) [d _{P'/P}]}
  \ar[r] _-\sim ^-{\ref{fg!prodtens}}
  &
  {g ^{(\bullet)!} (\E^{(\bullet)})  \smash{\overset{\L}{\otimes}} ^\dag _{\O _{\PP'} (\hdag D') _{ \Q}} g ^{(\bullet)!} (\E ^{\prime (\bullet)}) ,}
  }
\end{equation}
  les isomorphismes de recollement $\tau _{g,f}$ ayant été définis en \cite[2.1.10]{caro-construction}.
\end{enumerate}

\end{prop}

\begin{proof}
Par pleine fidélité du foncteur oubli du diviseur (et avec \ref{otimesTT'}), on se ramène au cas où $D$ est vide. 
L'isomorphisme \ref{fg!prodtens} 
et la vérification de la commutativité des carrés se ramène par construction à l'énoncé analogue au niveau des schémas. 
L'isomorphisme analogue à \ref{fg!prodtens} dans le cadre des morphismes de schémas est trivial. 
De même pour la commutativité du carré \ref{fg!prodtens-trans}.
Pour celle du carré \ref{fg!prodtens-square}, 
cela résulte de la définition de la $m$-PD-stratification associée à un produit tensoriel de $\D ^{(m)}$-modules et
  de la construction de $\tau _{g,f}$ (induit par
  un des isomorphismes de cette $m$-PD-stratification via le diagramme \cite[2.1.1.1]{caro-construction}).
\end{proof}

\subsection{Stabilité par produit tensoriel des isocristaux dans le cas relevable}

\begin{nota}
\label{nota-lambda-mu}
Soient $\lambda, \mu \colon \N \to \N$ deux applications croissantes
telles que $\lambda \geq \mu$, i.e., telles que, pour tout $m\in \N$, on ait $\lambda (m) \geq \mu (m)$.
On note alors $(\lambda \times \mu) ^* \smash{\widehat{\D}} _{\PP} ^{(\bullet)} (T )$ le système inductif 
$(\widehat{\B} ^{(\lambda (m))} _{\PP} ( T)  \smash{\widehat{\otimes}} _{\O _{\PP}} \smash{\widehat{\D}} _{\PP} ^{(\mu (m))}) _{m\in \N}$.
\end{nota}

\begin{prop}
\label{cohD-cohB-lim}
\begin{enumerate}
\item   Soient $\E$ un $\smash{\D} ^\dag _{\PP} (\hdag T) _{\Q}$-module cohérent,
$\O _{\PP} (\hdag T) _{\Q}$-cohérent. 
Avec les notations de \ref{nota-lambda-mu}, 
il existe alors $\lambda \colon \N \to \N$ tel que $\lambda \geq id$, 
$\E ^{(\bullet)}$ soit un $(\lambda \times id) ^* \smash{\widehat{\D}} _{\PP} ^{(\bullet)} (T )$-module localement de présentation finie 
dont le $\lambda ^{*}\widehat{\B} _{\PP } ^{(\bullet)} (T)$-module induit soit  localement de présentation finie et
un isomorphisme $\smash{\D} ^\dag _{\PP} (\hdag T) _{\Q}$-linéaire de la forme
$\underrightarrow{\lim} \,  \E ^{(\bullet)} \riso \E$
\item Soit
$\FF ^{(\bullet)}$ un 
$\smash{\widehat{\D}} _{\PP} ^{(\bullet)} (T )$-module cohérent à lim-ind-isogénie près.
La propriété 
{\og $\underrightarrow{\lim} \,  \FF ^{(\bullet)}$
est $\O _{\PP} (\hdag T) _{\Q}$-cohérent \fg}
équivaut à 
{\og $\FF ^{(\bullet)}$
est un 
$\widehat{\B} _{\PP } ^{(\bullet)} (T)$-module cohérent à lim-ind-isogénie près\fg}.
\end{enumerate}
\end{prop}

\begin{proof}
Comme $\E$ est associé à un isocristal surconvergent sur $(P\setminus T,P)/K$, 
la première assertion résulte de \cite[4.4.7]{Be1}.
Traitons à présent la seconde assertion. 
Il est immédiat que la seconde propriété implique la première. 
Réciproquement, supposons que
$\underrightarrow{\lim} \,  \FF ^{(\bullet)}$
est $\O _{\PP} (\hdag T) _{\Q}$-cohérent .
Comme $\underrightarrow{\lim} \,  \FF ^{(\bullet)}$
est un $\smash{\D} ^\dag _{\PP} (\hdag T) _{\Q}$-module cohérent,
$\O _{\PP} (\hdag T) _{\Q}$-cohérent, d'après ce que l'on vient de voir, 
il existe $\G ^{(\bullet)}$ un $\smash{\widehat{\D}} _{\PP} ^{(\bullet)} (T )$-module cohérent à lim-ind-isogénie près, 
$\widehat{\B} _{\PP } ^{(\bullet)} (T)$-cohérent à lim-ind-isogénie près tel que 
$\underrightarrow{\lim} \,  \G ^{(\bullet)} \riso \underrightarrow{\lim} \,  \FF ^{(\bullet)}$.
Comme le foncteur 
$ \underrightarrow{\lim} $
est pleinement fidèle sur la catégorie 
$\underrightarrow{LM}  _{\Q, \mathrm{coh}}
(\overset{^\mathrm{g}}{} \smash{\widehat{\D}} _{\PP} ^{(\bullet)} (T))$,
on en déduit l'isomorphisme $\G ^{(\bullet)} \riso \FF ^{(\bullet)}$ dans 
$\underrightarrow{LM}  _{\Q, \mathrm{coh}}
(\overset{^\mathrm{g}}{} \smash{\widehat{\D}} _{\PP} ^{(\bullet)} (T))$ et donc
dans 
$\underrightarrow{LM}   _{\Q, \mathrm{coh}}
(\overset{^\mathrm{g}}{} \smash{\widehat{\B}} _{\PP} ^{(\bullet)} (T))$.
\end{proof}

\begin{nota}
\label{nota-LDcapLD}
Notons 
$ \underrightarrow{LD}  ^\mathrm{b} _{\Q, \mathrm{coh}}
(\overset{^\mathrm{g}}{} \smash{\widehat{\D}} _{\PP} ^{(\bullet)} (T))
\cap
\underrightarrow{LD}  ^\mathrm{b} _{\Q, \mathrm{coh}}
(\widehat{\B} _{\PP} ^{(\bullet)} (T ))$,
la sous-catégorie pleine de 
$ \underrightarrow{LD}  ^\mathrm{b} _{\Q, \mathrm{coh}}
(\overset{^\mathrm{g}}{} \smash{\widehat{\D}} _{\PP} ^{(\bullet)} (T))$
des complexes dont l'image par le foncteur oubli
$\underrightarrow{LD}  ^\mathrm{b} _{\Q, \mathrm{qc}}
(\widehat{\D} _{\PP} ^{(\bullet)} (T ))
\to 
\underrightarrow{LD}  ^\mathrm{b} _{\Q, \mathrm{qc}}
(\widehat{\B} _{\PP} ^{(\bullet)} (T ))$
est dans 
$\underrightarrow{LD}  ^\mathrm{b} _{\Q, \mathrm{coh}}
(\widehat{\B} _{\PP} ^{(\bullet)} (T ))$.
En remplaçant {\og $\underrightarrow{LD}  ^\mathrm{b}$\fg} par 
{\og $\underrightarrow{LM} $\fg},
on définit de manière identique
$\underrightarrow{LM} _{\Q, \mathrm{coh}}
(\overset{^\mathrm{g}}{} \smash{\widehat{\D}} _{\PP} ^{(\bullet)} (T))
\cap
\underrightarrow{LM}  _{\Q, \mathrm{coh}}
(\widehat{\B} _{\PP} ^{(\bullet)} (T ))$.

La catégorie $\mathrm{Isoc} ^{\dag \dag}( \PP, T , P /K)$
est celle des $\smash{\D} ^\dag _{\PP} (\hdag T) _{\Q}$-modules cohérents,
$\O _{\PP} (\hdag T) _{\Q}$-cohérents.
Avec les notations \ref{nota-M-eq-isoc-lim}, 
la proposition \ref{cohD-cohB-lim} se traduit par l'égalité
\begin{equation}
\label{isocbullet-carac-PTP}
\mathrm{Isoc} ^{(\bullet)}( \PP, T , P /K)
= 
\underrightarrow{LM} _{\Q, \mathrm{coh}}
(\overset{^\mathrm{g}}{} \smash{\widehat{\D}} _{\PP} ^{(\bullet)} (T))
\cap
\underrightarrow{LM}  _{\Q, \mathrm{coh}}
(\widehat{\B} _{\PP} ^{(\bullet)} (T )).
\end{equation}

\end{nota}

\begin{prop}
\label{stab-isoc-bullet-pas}
Soient
$\E ^{(\bullet)}$
et
$\E ^{\prime (\bullet)}$
deux objets de 
$\mathrm{Isoc} ^{(\bullet)}( \PP, T , P /K)$.

\begin{enumerate}
\item On bénéficie des isomorphismes de 
$ \underrightarrow{LD}  ^\mathrm{b} _{\Q, \mathrm{coh}}
(\overset{^\mathrm{g}}{} \smash{\widehat{\D}} _{\PP} ^{(\bullet)} (T))
\cap
\underrightarrow{LD}  ^\mathrm{b} _{\Q, \mathrm{coh}}
(\widehat{\B} _{\PP} ^{(\bullet)} (T ))$ :
\begin{equation}
\label{form-proj-alpha}
\E  ^{(\bullet)} \otimes _{\widehat{\B} _{\PP } ^{(\bullet)} (T )}\E  ^{\prime(\bullet)}
\liso
\E  ^{(\bullet)} \otimes ^{\L} _{\widehat{\B} _{\PP } ^{(\bullet)} (T )}\E  ^{\prime(\bullet)}
\riso 
\E  ^{(\bullet)} 
\smash{\overset{\L}{\otimes}} ^\dag _{\O _{\PP} (\hdag T ) _{\Q}} 
\E  ^{\prime(\bullet)}.
\end{equation}
De plus, 
$\E  ^{(\bullet)} \otimes _{\widehat{\B} _{\PP } ^{(\bullet)} (T )}\E  ^{\prime(\bullet)}
\in \mathrm{Isoc} ^{(\bullet)}( \PP, T , P /K)$.

\item On dispose des isomorphismes :
\begin{equation}
  \label{rema15}
\E \otimes _{\O _{\PP} (\hdag T) _{\Q}} \E' 
\liso
\E \otimes ^\L _{\O _{\PP} (\hdag T) _{\Q}} \E '
\riso
\E \smash{\overset{\L}{\otimes}}   ^{\dag} _{\O _{\PP} (\hdag T) _{\Q}} \E '.
\end{equation}

\end{enumerate}

\end{prop}

\begin{proof}
D'après \ref{cohD-cohB-lim}, on peut supposer qu'il existe
$\lambda \geq id\colon \N\to \N$ telle que 
$\E ^{(\bullet)}$ et 
$\E ^{\prime (\bullet)}$
soient deux 
$(\lambda \times id) ^* \smash{\widehat{\D}} _{\PP} ^{(\bullet)} (T )$-module localement de présentation finie 
dont le $\lambda ^{*}\widehat{\B} _{\PP } ^{(\bullet)} (T)$-module induit soit  localement de présentation finie.
Comme le morphisme canonique 
$\E  ^{(\bullet)} \otimes _{\widehat{\B} _{\PP } ^{(\bullet)} (T )}\E  ^{\prime(\bullet)}
\to
\E  ^{(\bullet)} \otimes _{\lambda ^* \widehat{\B} _{\PP } ^{(\bullet)} (T )}\E  ^{\prime(\bullet)}$
est un isomorphisme de 
$\underrightarrow{LD}  ^\mathrm{b} _{\Q, \mathrm{coh}}
(\overset{^\mathrm{g}}{} \smash{\widehat{\D}} _{\PP} ^{(\bullet)} (T))$, 
on obtient alors que
$\E  ^{(\bullet)} \otimes _{\widehat{\B} _{\PP } ^{(\bullet)} (T )}\E  ^{\prime(\bullet)}$
est un $(\lambda \times id) ^* \smash{\widehat{\D}} _{\PP} ^{(\bullet)} (T )$-module localement de présentation finie 
dont le $\lambda ^{*}\widehat{\B} _{\PP } ^{(\bullet)} (T)$-module induit est localement de présentation finie.
En particulier,
$\E  ^{(\bullet)} \otimes _{\widehat{\B} _{\PP } ^{(\bullet)} (T )}\E  ^{\prime(\bullet)}\in \mathrm{Isoc} ^{(\bullet)}( \PP, T , P /K)$.
De plus, on vérifie aussi facilement que le second morphisme canonique de \ref{form-proj-alpha} est un isomorphisme.
En lui appliquant le foncteur $\underrightarrow{\lim}$,
on en déduit la même chose pour \ref{rema15}. 
Comme $\E$ est (associé à) un isocristal surconvergent sur $(P\setminus T,P)/K$,
alors $\E$ est  un $\O _{\PP} (\hdag T) _{\Q}$-module localement projectif de type fini et par conséquent on obtient 
le premier isomorphisme de \ref{rema15}. 
Comme le foncteur $\underrightarrow{\lim}$ est pleinement fidèle sur 
$\underrightarrow{LD}  ^\mathrm{b} _{\Q, \mathrm{coh}}
(\widehat{\D} _{\PP} ^{(\bullet)} (T ))$, on en déduit que le premier 
morphisme de \ref{form-proj-alpha} est aussi un isomorphisme.

\end{proof}

\subsection{Produits tensoriels externes: définition et propriétés}

On pose $\S := \Spf \V$. Soient $\X$, $\Y$ deux $\V$-schémas formels lisses, $\ZZ := \X \times _\S \Y$, 
et $f\colon \ZZ \to \X$, 
$g \colon \ZZ \to \Y$ les deux projections.
Soient 
$T _1$ un diviseur de $X$, $T _2$ un diviseur de $Y$
et $T := f ^{-1}(T _1) \cup g ^{-1} (T _2)$. 

\begin{vide}
Soient
$\E ^{(\bullet)}
\in \underrightarrow{LD}  ^\mathrm{b} _{\Q, \mathrm{qc}}
(\overset{^\mathrm{g}}{} \smash{\widehat{\D}} _{\X} ^{(\bullet)} (T _1))$
et 
$\FF ^{(\bullet)}
\in \underrightarrow{LD}  ^\mathrm{b} _{\Q, \mathrm{qc}}
(\overset{^\mathrm{g}}{} \smash{\widehat{\D}} _{\Y} ^{(\bullet)}(T _2))$.
De manière analogue à \cite[2.3.1, 3.4.7, 4.3.5]{Beintro2} et en utilisant
\ref{otimesTT'}, on définit le bifoncteur produit tensoriel externe
\begin{equation}
\label{prodtens-ext-qc}
\smash{\overset{\L}{\boxtimes}}   ^{\dag} _{\O _\S} 
\colon 
\underrightarrow{LD}  ^\mathrm{b} _{\Q, \mathrm{qc}}
(\overset{^\mathrm{g}}{} \smash{\widehat{\D}} _{\X} ^{(\bullet)} (T _1) )
\times 
\underrightarrow{LD}  ^\mathrm{b} _{\Q, \mathrm{qc}}
(\overset{^\mathrm{g}}{} \smash{\widehat{\D}} _{\Y} ^{(\bullet)}(T _2) )
\to
\underrightarrow{LD}  ^\mathrm{b} _{\Q, \mathrm{qc}}
(\overset{^\mathrm{g}}{} \smash{\widehat{\D}} _{\ZZ} ^{(\bullet)}(T)),
\end{equation}
en posant 
$\E ^{(\bullet)}
\smash{\overset{\L}{\boxtimes}}   ^{\dag} _{\O _\S}
\FF ^{(\bullet)}:=
f _{T _1}^{(\bullet)*} (\E ^{(\bullet)})
\smash{\overset{\L}{\otimes}}   ^{\dag} _{\O _{\ZZ,\Q}}
g _{T _2} ^{(\bullet)*} (\FF ^{(\bullet)})$
(voir les notations de \ref{nota-f*TT'} et \ref{f+!-indTT'}).
Modulo les foncteurs oublis des diviseurs 
$oub _{T _1}$, $oub _{T _2}$, $oub _{T}$
qui sont pleinement fidèles,
ces bifoncteurs sont identiques et ne dépendent pas des diviseurs $T _1$ ou $T _2$.
Il est donc anodin de ne pas indiquer les diviseurs $T _1$ et $T _2$ dans la notation
$\smash{\overset{\L}{\boxtimes}}   ^{\dag} _{\O _\S} $. 
\end{vide}

\begin{vide}
\label{calcul-DT1T2-lim}
On note 
$\smash{\widehat{\B}} _{\ZZ} ^{(m)} (T _1,T _2):=
f ^{(m)*} (\smash{\widehat{\B}} _{\X} ^{(m)}  ( T _1))
\smash{\widehat{\otimes}}  _{\O _{\ZZ}}
g ^{(m)*} ( \smash{\widehat{\B}} _{\Y} ^{(m)}  ( T _2))$, 
qui est une $\O _{\ZZ}$-algèbre munie d'une structure 
 compatible canonique de $\smash{\widehat{\D}} _{\ZZ} ^{(m)}$-module à gauche.
Via \cite[5.2.1]{caro-stab-sys-ind-surcoh}, on dispose de l'isomorphisme canonique de la forme
$\smash{\widehat{\B}} _{\ZZ} ^{(m)} (T _1,T _2)
\riso
\smash{\widehat{\B}} _{\ZZ} ^{(m)}  ( f ^{-1}(T _1))
\smash{\widehat{\otimes}}  _{\O _{\ZZ}}
 \smash{\widehat{\B}} _{\ZZ} ^{(m)}  ( g ^{-1} (T _2)) $.
Par \cite[3.2.7.2]{caro-stab-sys-ind-surcoh}, 
on en déduit que le morphisme canonique
$f ^{(m)*} (\smash{\widehat{\B}} _{\X} ^{(m)}  ( T _1))
\smash{\widehat{\otimes}} ^{\L} _{\O _{\ZZ}}
g ^{(m)*} ( \smash{\widehat{\B}} _{\Y} ^{(m)}  ( T _2))
\to 
\smash{\widehat{\B}} _{\ZZ} ^{(m)} (T _1,T _2)$
est un isomorphisme.

On obtient un $\smash{\widehat{\D}} _{\ZZ} ^{(m)}$-module à gauche en posant
$\smash{\widehat{\D}} _{\ZZ} ^{(m)} (T _1,T _2):=
f ^{(m)*} (\smash{\widehat{\D}} _{\X} ^{(m)}  ( T _1))
\smash{\widehat{\otimes}}  _{\O _{\ZZ}}
g ^{(m)*} ( \smash{\widehat{\D}} _{\Y} ^{(m)}  ( T _2))$. 
On dispose de l'isomorphisme canonique de $\smash{\widehat{\D}} _{\ZZ} ^{(m)}$-modules à gauche
$\smash{\widehat{\B}} _{\ZZ} ^{(m)} (T _1,T _2) 
\smash{\widehat{\otimes}}  _{\O _{\ZZ}}
\smash{\widehat{\D}} _{\ZZ} ^{(m)} 
\riso \smash{\widehat{\D}} _{\ZZ} ^{(m)} (T _1,T _2)$.
Il en résulte que le morphisme naturel
$f ^{(m)*} (\smash{\widehat{\D}} _{\X} ^{(m)}  ( T _1))
\smash{\widehat{\otimes}}  ^{\L} _{\O _{\ZZ}}
g ^{(m)*} ( \smash{\widehat{\D}} _{\Y} ^{(m)}  ( T _2))
\to
\smash{\widehat{\D}} _{\ZZ} ^{(m)} (T _1,T _2)$
est un isomorphisme.

On pose 
$\smash{\widehat{\B}} _{\ZZ} ^{(\bullet)} (T _1,T _2):=
f ^{(\bullet)*} (\smash{\widehat{\B}} _{\X} ^{(\bullet)}  ( T _1))
\smash{\overset{\L}{\otimes}}   ^{\dag} _{\O _{\ZZ,\Q}}
g ^{(\bullet)*} ( \smash{\widehat{\B}} _{\Y} ^{(\bullet)}  ( T _2))$
et
$\smash{\widehat{\D}} _{\ZZ} ^{(\bullet)} (T _1,T _2):=
f ^{(\bullet)*} (\smash{\widehat{\D}} _{\X} ^{(\bullet)}  ( T _1))
\smash{\overset{\L}{\otimes}}   ^{\dag} _{\O _{\ZZ,\Q}}
g ^{(\bullet)*} ( \smash{\widehat{\D}} _{\Y} ^{(\bullet)}  ( T _2))$.
Avec \cite[3.2.7.4]{caro-stab-sys-ind-surcoh},
le morphisme canonique
$\smash{\widehat{\B}} _{\ZZ} ^{(\bullet)} (T _1,T _2)
\to
\smash{\widehat{\B}} _{\ZZ} ^{(\bullet)} (T )$ 
est un monomorphisme de $D (\overset{^\mathrm{g}}{} \smash{\widehat{\D}} _{\ZZ} ^{(\bullet)})$
et un isomorphisme de
$\underrightarrow{LD}  ^\mathrm{b} _{\Q, \mathrm{qc}}
(\overset{^\mathrm{g}}{} \smash{\widehat{\D}} _{\ZZ} ^{(\bullet)})$.

Avec \cite[2.3.1]{Beintro2}, on dispose de l'isomorphisme canonique $\smash{\widehat{\D}} _{\ZZ} ^{(\bullet)}$-linéaire 
$
f ^{(\bullet)*} (\smash{\widehat{\D}} _{\X} ^{(\bullet)})
\smash{\overset{\L}{\otimes}}   ^{\dag} _{\O _{\ZZ,\Q}}
g ^{(\bullet)*} ( \smash{\widehat{\D}} _{\Y} ^{(\bullet)} )
\riso
\smash{\widehat{\D}} _{\ZZ} ^{(\bullet)}$.
On obtient alors le morphisme canonique
$\smash{\widehat{\D}} _{\ZZ} ^{(\bullet)} (T _1,T _2)
\to
\smash{\widehat{\D}} _{\ZZ} ^{(\bullet)} (T )$
qui est un monomorphisme de $D (\overset{^\mathrm{g}}{} \smash{\widehat{\D}} _{\ZZ} ^{(\bullet)})$
et un isomorphisme de
$\underrightarrow{LD}  ^\mathrm{b} _{\Q, \mathrm{qc}}
(\overset{^\mathrm{g}}{} \smash{\widehat{\D}} _{\ZZ} ^{(\bullet)})$. 
\end{vide}

\begin{vide}
[Produit tensoriel externe de complexes cohérents]
\label{nota-prodext-coh}

Soient $\E ^{(\bullet)} 
\in \underrightarrow{LD}  ^\mathrm{b} _{\Q, \mathrm{coh}}
(\overset{^\mathrm{g}}{} \smash{\widehat{\D}} _{\X} ^{(\bullet)} (T _1))$,
$\FF ^{(\bullet)} 
\in \underrightarrow{LD}  ^\mathrm{b} _{\Q, \mathrm{coh}}
(\overset{^\mathrm{g}}{} \smash{\widehat{\D}} _{\Y} ^{(\bullet)} (T _2))$.
Posons
$\E:  = \underrightarrow{\lim}\, \E ^{(\bullet)} \in D ^\mathrm{b} _\mathrm{coh} ( \smash{\D} ^\dag _{\X} (\hdag T _1) _{\Q} ) $,
$\FF:  = \underrightarrow{\lim} \, \FF ^{(\bullet)} \in 
D ^\mathrm{b} _\mathrm{coh} ( \smash{\D} ^\dag _{\Y} (\hdag T _2) _{\Q} )$.

De manière analogue à 
 \cite[4.3.5]{Beintro2} 
 (on utilise l'isomorphisme $\smash{\widehat{\D}} _{\ZZ} ^{(\bullet)} (T _1,T _2)
\riso
\smash{\widehat{\D}} _{\ZZ} ^{(\bullet)} (T )$ de \ref{calcul-DT1T2-lim}), on vérifie que le bifoncteur produit tensoriel externe de 
 \ref{prodtens-ext-qc} induit le bifoncteur
 \begin{equation}
\label{prodtens-ext-coh}
\smash{\overset{\L}{\boxtimes}}   ^{\dag} _{\O _\S} 
\colon 
\underrightarrow{LD}  ^\mathrm{b} _{\Q, \mathrm{coh}}
(\overset{^\mathrm{g}}{} \smash{\widehat{\D}} _{\X} ^{(\bullet)} (T _1) )
\times 
\underrightarrow{LD}  ^\mathrm{b} _{\Q, \mathrm{coh}}
(\overset{^\mathrm{g}}{} \smash{\widehat{\D}} _{\Y} ^{(\bullet)}(T _2) )
\to
\underrightarrow{LD}  ^\mathrm{b} _{\Q, \mathrm{coh}}
(\overset{^\mathrm{g}}{} \smash{\widehat{\D}} _{\ZZ} ^{(\bullet)}(T)).
\end{equation}
Avec les notations de \ref{fct-qcoh2coh}, 
on obtient le bifoncteur 
$\smash{\overset{\L}{\boxtimes}}   ^{\dag} _{\O _\S, T _1, T _2}
:=
\mathrm{Coh} _{T _1, T _2} 
(\smash{\overset{\L}{\boxtimes}}   ^{\dag} _{\O _\S})$ de la forme:
\begin{equation}
\label{prod-tens-extlimcoh}
\smash{\overset{\L}{\boxtimes}}   ^{\dag} _{\O _\S, T _1, T _2} \colon
D ^\mathrm{b} _\mathrm{coh} ( \smash{\D} ^\dag _{\X} (\hdag T _1) _{\Q} ) 
\times
D ^\mathrm{b} _\mathrm{coh} ( \smash{\D} ^\dag _{\Y} (\hdag T _2) _{\Q} )
\to 
D ^\mathrm{b} _\mathrm{coh} ( \smash{\D} ^\dag _{\ZZ} (\hdag T ) _{\Q} ).
\end{equation}
Par définition, on a donc l'isomorphisme
$\E
\smash{\overset{\L}{\boxtimes}}   ^{\dag} _{\O _\S, T _1, T _2}
\FF
\riso
 \underrightarrow{\lim} ( \E ^{(\bullet)}
\smash{\overset{\L}{\boxtimes}}   ^{\dag} _{\O _\S}
\FF ^{(\bullet)})$.
Par exemple
$
\smash{\D} ^\dag _{\X} (\hdag T _1) _{\Q} 
\smash{\overset{\L}{\boxtimes}}   ^{\dag} _{\O _\S}
\smash{\D} ^\dag _{\Y} (\hdag T _2) _{\Q} 
=\smash{\D} ^\dag _{\ZZ} (\hdag T ) _{\Q} .$
\end{vide}

\begin{rema}
Soient
$D _1\subset T _1$ un second diviseur de $X$,
$D _2\subset T _2$ un second diviseur de $Y$. 
Soient 
$\E \in 
D ^\mathrm{b} _\mathrm{coh} ( \smash{\D} ^\dag _{\X} (\hdag D _1) _{\Q} )
\cap 
D ^\mathrm{b} _\mathrm{coh} ( \smash{\D} ^\dag _{\X} (\hdag T _1) _{\Q} )$ 
et 
$\FF\in 
D ^\mathrm{b} _\mathrm{coh} ( \smash{\D} ^\dag _{\Y} (\hdag D _2) _{\Q} )
\cap 
D ^\mathrm{b} _\mathrm{coh} ( \smash{\D} ^\dag _{\Y} (\hdag T _2) _{\Q} )$.
D'après \ref{fct-qcoh2coh}, 
le morphisme canonique 
\begin{equation}
\label{boxtimesDT-fleche-surcoh}
\E
\smash{\overset{\L}{\boxtimes}}   ^{\dag} _{\O _\S, D _1, D _2}\,
\FF
\to 
\E 
\smash{\overset{\L}{\boxtimes}}   ^{\dag} _{\O _\S, T _1, T _2}\,
\FF
\end{equation}
est un isomorphisme.
\end{rema}

\begin{lemm}
\label{lemm-boxtimes-TouT'}
\begin{enumerate}
\item Pour tous 
$\E ^{(\bullet)} 
\in \underrightarrow{LD}  ^\mathrm{b} _{\Q, \mathrm{qc}}
(\overset{^\mathrm{g}}{} \smash{\widehat{\D}} _{\X} ^{(\bullet)} (T _1))$,
$\FF ^{(\bullet)} 
\in \underrightarrow{LD}  ^\mathrm{b} _{\Q, \mathrm{qc}}
(\overset{^\mathrm{g}}{} \smash{\widehat{\D}} _{\Y} ^{(\bullet)} (T _2))$,
on bénéficie de l'isomorphisme canonique dans $\underrightarrow{LD}  ^\mathrm{b} _{\Q, \mathrm{coh}}
(\overset{^\mathrm{g}}{} \smash{\widehat{\D}} _{\ZZ} ^{(\bullet)}(T))$:
\begin{equation}
\label{lemm-boxtimes-TouT'1pre}
\E ^{(\bullet)}
\smash{\overset{\L}{\boxtimes}}   ^{\dag} _{\O _\S}
\FF ^{(\bullet)}
\riso 
(\hdag T ) \circ f ^{(\bullet)*} (\E ^{(\bullet)})
\smash{\overset{\L}{\otimes}}   ^{\dag} _{\O _{\ZZ} (\hdag T ) _{\Q}}
(\hdag T )\circ  g ^{(\bullet)*} (\FF^{(\bullet)}).
\end{equation}

\item On dispose, 
pour tous
$\E \in D ^\mathrm{b} _\mathrm{coh} ( \smash{\D} ^\dag _{\X} (\hdag T _1) _{\Q} ) $,
$\FF\in 
D ^\mathrm{b} _\mathrm{coh} ( \smash{\D} ^\dag _{\Y} (\hdag T _2) _{\Q} )$,
de l'isomorphisme canonique dans $D ^\mathrm{b} _\mathrm{coh} ( \smash{\D} ^\dag _{\ZZ} (\hdag T ) _{\Q} )$:
\begin{equation}
\label{lemm-boxtimes-TouT'1}
\E 
\smash{\overset{\L}{\boxtimes}}   ^{\dag} _{\O _\S, T _1, T _2}
\FF
\riso 
(\hdag T )\circ  f ^{*} _{T _1}(\E)
\smash{\overset{\L}{\otimes}}   ^{\dag} _{\O _{\ZZ} (\hdag T ) _{\Q}}
(\hdag T ) \circ  g ^* _{T _2} (\FF).
\end{equation}

\end{enumerate}

\end{lemm}

\begin{proof}
L'isomorphisme \ref{lemm-boxtimes-TouT'1pre} découle de \ref{otimesTT'}. 
De plus, pour tous complexes
$\E ^{(\bullet)} 
\in \underrightarrow{LD}  ^\mathrm{b} _{\Q, \mathrm{coh}}
(\overset{^\mathrm{g}}{} \smash{\widehat{\D}} _{\X} ^{(\bullet)} (T _1))$,
$\FF ^{(\bullet)} 
\in \underrightarrow{LD}  ^\mathrm{b} _{\Q, \mathrm{coh}}
(\overset{^\mathrm{g}}{} \smash{\widehat{\D}} _{\Y} ^{(\bullet)} (T _2))$
tels que 
$\E  = \underrightarrow{\lim}\, \E ^{(\bullet)}$
et
$\FF  = \underrightarrow{\lim}\, \FF ^{(\bullet)}$,
comme 
$(\hdag T ) \circ f ^{(\bullet)*} (\E ^{(\bullet)}) $
et
$(\hdag T ) \circ  g ^{(\bullet)*} (\FF ^{(\bullet)})$
sont des objets de
$\underrightarrow{LD}  ^\mathrm{b} _{\Q, \mathrm{coh}}
(\overset{^\mathrm{g}}{} \smash{\widehat{\D}} _{\ZZ} ^{(\bullet)}(T))$, 
en appliquant  le foncteur $\underrightarrow{\lim} $
à 
\ref{lemm-boxtimes-TouT'1pre}, 
on obtient \ref{lemm-boxtimes-TouT'1}.
\end{proof}

\begin{lemm}
Soient 
$T '_1 \supset T _1$ un second diviseur de $X$, $T' _2 \supset T _2$ un second diviseur de $Y$
et $T' := f ^{-1}(T' _1) \cup g ^{-1} (T' _2)$. 
\begin{enumerate}
\item 
Pour tous 
$\E ^{(\bullet)} 
\in \underrightarrow{LD}  ^\mathrm{b} _{\Q, \mathrm{qc}}
(\overset{^\mathrm{g}}{} \smash{\widehat{\D}} _{\X} ^{(\bullet)} (T _1))$,
$\FF ^{(\bullet)} 
\in \underrightarrow{LD}  ^\mathrm{b} _{\Q, \mathrm{qc}}
(\overset{^\mathrm{g}}{} \smash{\widehat{\D}} _{\Y} ^{(\bullet)} (T _2))$,
on bénéficie de l'isomorphisme canonique dans 
$\underrightarrow{LD}  ^\mathrm{b} _{\Q, \mathrm{qc}}
(\overset{^\mathrm{g}}{} \smash{\widehat{\D}} _{\ZZ} ^{(\bullet)}(T'))$:
\begin{equation}
\label{lemm-boxtimes-TouT'2pre}
(\hdag T') \left (\E ^{(\bullet)}
\smash{\overset{\L}{\boxtimes}}   ^{\dag} _{\O _\S}
\FF ^{(\bullet)} \right ) 
\riso 
(\hdag T _1 ')(\E ^{(\bullet)})
\smash{\overset{\L}{\boxtimes}}   ^{\dag} _{\O _\S}
(\hdag T _2 ') (\FF ^{(\bullet)}).
\end{equation}
\item 
On dispose, 
pour tous
$\E \in D ^\mathrm{b} _\mathrm{coh} ( \smash{\D} ^\dag _{\X} (\hdag T _1) _{\Q} ) $,
$\FF\in 
D ^\mathrm{b} _\mathrm{coh} ( \smash{\D} ^\dag _{\Y} (\hdag T _2) _{\Q} )$, 
de l'isomorphisme canonique dans 
$D ^\mathrm{b} _\mathrm{coh} ( \smash{\D} ^\dag _{\ZZ} (\hdag T ') _{\Q} )$:
\begin{equation}
\label{lemm-boxtimes-TouT'2}
(\hdag T') \left (\E 
\smash{\overset{\L}{\boxtimes}}   ^{\dag} _{\O _\S, T _1, T _2}
\FF \right ) 
\riso 
(\hdag T _1 ')(\E )
\smash{\overset{\L}{\boxtimes}}   ^{\dag} _{\O _\S, T '_1, T '_2}
(\hdag T _2 ') (\FF ).
\end{equation}

\end{enumerate}
\end{lemm}

\begin{proof}

Comme les foncteurs images inverses extraordinaires commutent aux foncteurs
de localisation en dehors d'un diviseur (voir \cite[1.1.9--10]{caro_courbe-nouveau}), 
on obtient le dernier des isomorphismes:
\begin{gather}
\notag
(\hdag T') \left ( 
f ^{(\bullet)*} (\E ^{(\bullet)})
\smash{\overset{\L}{\otimes}}   ^{\dag} _{\O _{\ZZ,\Q}}
g ^{(\bullet)*} (\FF ^{(\bullet)})
\right )
\riso 
(\hdag T') \circ  f ^{(\bullet)*} (\E ^{(\bullet)})
\smash{\overset{\L}{\otimes}}   ^{\dag} _{\O _{\ZZ} (\hdag T ') _{\Q}}
(\hdag T') \circ  g ^{(\bullet)*} (\FF  ^{(\bullet)})
\\
\notag
\underset{\ref{otimesTT'}}{\riso}
(\hdag f ^{-1} (T '_1) ) \circ  f ^{(\bullet)*} (\E ^{(\bullet)})
\smash{\overset{\L}{\otimes}}   ^{\dag} _{\O _{\ZZ,\Q}}
(\hdag g ^{-1} (T '_2) ) \circ  g ^{(\bullet)*} (\FF  ^{(\bullet)})
\\
\label{lemm-boxtimes-TouT'2-iso}
\riso
 f ^{(\bullet)*} \circ (\hdag T ' _1) (\E ^{(\bullet)})
\smash{\overset{\L}{\otimes}}   ^{\dag} _{\O _{\ZZ,\Q}}
g ^{(\bullet)*} \circ (\hdag T ' _2)(\FF  ^{(\bullet)}),
\end{gather}
dont  la composition est l'isomorphisme \ref{lemm-boxtimes-TouT'2pre}.
Si
$\E ^{(\bullet)} 
\in \underrightarrow{LD}  ^\mathrm{b} _{\Q, \mathrm{coh}}
(\overset{^\mathrm{g}}{} \smash{\widehat{\D}} _{\X} ^{(\bullet)} (T _1))$,
$\FF ^{(\bullet)} 
\in \underrightarrow{LD}  ^\mathrm{b} _{\Q, \mathrm{coh}}
(\overset{^\mathrm{g}}{} \smash{\widehat{\D}} _{\Y} ^{(\bullet)} (T _2))$,
comme alors
$(\hdag T ' _1) (\E ^{(\bullet)})
\in 
\underrightarrow{LD}  ^\mathrm{b} _{\Q, \mathrm{coh}}
(\overset{^\mathrm{g}}{} \smash{\widehat{\D}} _{\X} ^{(\bullet)} (T ' _1) )$
et
$(\hdag T ' _2)(\FF  ^{(\bullet)})
\in 
\underrightarrow{LD}  ^\mathrm{b} _{\Q, \mathrm{coh}}
(\overset{^\mathrm{g}}{} \smash{\widehat{\D}} _{\Y} ^{(\bullet)}(T ' _2) )
$,
en appliquant  le foncteur $\underrightarrow{\lim} $
à la composition des isomorphismes de \ref{lemm-boxtimes-TouT'2-iso}, on obtient alors \ref{lemm-boxtimes-TouT'1}.
\end{proof}

\begin{prop}
Soient $u\colon \X' \to \X$ et $v\colon \Y' \to \Y$ deux morphismes de $\V$-schémas formels lisses tels que
$T' _1 := u ^{-1}(T _1 )$ soit un diviseur de $X'$
et
$T' _2 := v ^{-1}(T _2 )$ soit un diviseur de $Y'$.
On note 
$\ZZ ' := \X '\times \Y'$,
$w:= (u, v) \colon \ZZ ' \to \ZZ$ le morphisme induit
et $T' := w ^{-1} (T)$.

\begin{enumerate}
\item Pour tous
$\E ^{\prime(\bullet)}
\in \underrightarrow{LD}  ^\mathrm{b} _{\Q, \mathrm{qc}}
(\overset{^\mathrm{g}}{} \smash{\widehat{\D}} _{\X'} ^{(\bullet)} ( T '_1))$
et 
$\FF ^{\prime (\bullet)}
\in \underrightarrow{LD}  ^\mathrm{b} _{\Q, \mathrm{qc}}
(\overset{^\mathrm{g}}{} \smash{\widehat{\D}} _{\Y'} ^{(\bullet)} (T '_2))$, 
on dispose dans $\underrightarrow{LD}  ^\mathrm{b} _{\Q, \mathrm{qc}}
(\overset{^\mathrm{g}}{} \smash{\widehat{\D}} _{\ZZ} ^{(\bullet)}(T))$
de l'isomorphisme:
\begin{equation}
\label{boxtimes-v+}
w ^{(\bullet)} _+ (\E ^{\prime(\bullet)}
\smash{\overset{\L}{\boxtimes}}   ^{\dag} _{\O _\S}
\FF ^{\prime(\bullet)})
\riso
u ^{(\bullet)}_+ (\E ^{\prime(\bullet)})
\smash{\overset{\L}{\boxtimes}}   ^{\dag} _{\O _\S}
v ^{(\bullet)}_+ (\FF ^{\prime(\bullet)}).
\end{equation}

\item Supposons $u$ et $v$ propres. 
Pour tous
$\E' \in D ^\mathrm{b} _\mathrm{coh} ( \smash{\D} ^\dag _{\X'} (\hdag T '_1) _{\Q} ) $, 
$\FF'\in 
D ^\mathrm{b} _\mathrm{coh} ( \smash{\D} ^\dag _{\Y'} (\hdag T '_2) _{\Q} )$,
on dispose alors de l'isomorphisme dans $D ^\mathrm{b} _\mathrm{coh} ( \smash{\D} ^\dag _{\ZZ} (\hdag T ) _{\Q} ) $:
\begin{equation}
\label{boxtimes-v+coh}
w _{T,+} (\E '
\smash{\overset{\L}{\boxtimes}}   ^{\dag} _{\O _\S, T '_1, T '_2}\,
\FF ')
\riso
u _{T_1,+} (\E ')
\smash{\overset{\L}{\boxtimes}}   ^{\dag} _{\O _\S, T _1, T _2}\,
v _{T_2,+} (\FF ').
\end{equation}
\end{enumerate}

\end{prop}

\begin{proof}
Établissons d'abord l'isomorphisme \ref{boxtimes-v+}.
Par symétrie et composition, il suffit de le prouver
lorsque $v $ est l'identité. On note alors $\FF ^{(\bullet)}$ à la place de $\FF ^{\prime(\bullet)}$.
On pose $f '\colon \ZZ' \to \X'$ et $g ' \colon \ZZ' \to \Y$ les projections canoniques. 
Il résulte de 
\ref{otimesTT'-iso}
 l'isomorphisme:
\begin{equation}
\notag
u ^{(\bullet)} _{+} (\E ^{\prime (\bullet)})
\smash{\overset{\L}{\boxtimes}}   ^{\dag} _{\O _\S}\,
\FF ^{(\bullet)}
\riso
(\hdag T  ) f ^{(\bullet)*} (u  ^{(\bullet)} _{+} (\E ^{\prime (\bullet)}))
\smash{\overset{\L}{\otimes}}   ^{\dag} _{\O _{\ZZ} (\hdag T) _{\Q}}
(\hdag T )  g ^{(\bullet)*} (\FF ^{(\bullet)}).
\end{equation}
De plus, $(\hdag T)\circ  f ^{(\bullet)*} \circ  u ^{(\bullet)} _{+} (\E ^{\prime (\bullet)})
\underset{\cite[5.7]{Abe-Frob-Poincare-dual}}{\riso}   
(\hdag T) \circ w ^{(\bullet)} _{+} \circ f ^{\prime *}  (\E ^{\prime (\bullet)})
\underset{\cite[2.2.18.2]{caro_surcoherent}}{\riso}   
w ^{(\bullet)} _{+}\circ (\hdag T' ) \circ  f ^{\prime *} (\E ^{\prime (\bullet)})$.
D'où, 
\begin{equation}
\notag
u ^{(\bullet)} _{+} (\E ^{\prime (\bullet)})
\smash{\overset{\L}{\boxtimes}}   ^{\dag} _{\O _\S}\,
\FF ^{(\bullet)}
\riso
w ^{(\bullet)} _{+}\circ (\hdag T' ) \circ  f ^{\prime *}  (\E ^{\prime (\bullet)})
\smash{\overset{\L}{\otimes}}   ^{\dag} _{\O _{\ZZ} (\hdag T) _{\Q}}
(\hdag T)  g ^{(\bullet)*} (\FF ^{(\bullet)}).
\end{equation}
D'après \ref{2.1.4-caro-surcoh},
on dispose de l'isomorphisme :
$$w ^{(\bullet)} _{+}\circ (\hdag T') \circ f ^{\prime *}  (\E ^{\prime (\bullet)})
\smash{\overset{\L}{\otimes}}   ^{\dag} _{\O _{\ZZ} (\hdag T) _{\Q}}
(\hdag T)  g ^{(\bullet)*} (\FF ^{(\bullet)})
\liso
w ^{(\bullet)} _{+}\left ( (\hdag T' ) \circ  f ^{\prime *} (\E ^{\prime (\bullet)})
\smash{\overset{\L}{\otimes}}   ^{\dag} _{\O _{\ZZ '} (\hdag T' ) _{\Q}}
\L w  ^{(\bullet)*}  \circ (\hdag T)  \circ g ^{(\bullet)*} (\FF ^{(\bullet)}) \right)
.$$
Or, $\L w  ^{(\bullet)*}  \circ (\hdag T)  \circ g ^{(\bullet)*} (\FF ^{(\bullet)})
\underset{\cite[2.2.18.1]{caro_surcoherent}}{\riso}      
(\hdag T') \circ \L w  ^{(\bullet)*}  \circ g ^{(\bullet)*} (\FF ^{(\bullet)})
\riso
(\hdag T') \circ g ^{\prime *} (\FF ^{(\bullet)})$.
Cela donne
\begin{gather}
\notag
u ^{(\bullet)} _{+} (\E ^{\prime (\bullet)})
\smash{\overset{\L}{\boxtimes}}   ^{\dag} _{\O _\S}\,
\FF ^{(\bullet)}
\riso
w ^{(\bullet)} _{+}\left ( (\hdag T' ) \circ   f ^{\prime *} (\E ^{\prime (\bullet)})
\smash{\overset{\L}{\otimes}}   ^{\dag} _{\O _{\ZZ '} (\hdag T') _{\Q}}
(\hdag T' ) \circ g ^{\prime *} (\FF ^{(\bullet)}) \right)
\\
\underset{\ref{otimesTT'-iso}}{\riso}
w ^{(\bullet)} _{+} \left (\E ^{\prime (\bullet)}
\smash{\overset{\L}{\boxtimes}}   ^{\dag} _{\O _\S}\,
\FF ^{(\bullet)} \right ).
\end{gather}
Par stabilité de la cohérence par l'image directe d'un morphisme propre,
on obtient \ref{boxtimes-v+coh} en appliquant le foncteur $\underrightarrow{\lim} $
à \ref{boxtimes-v+}.
\end{proof}

\begin{lemm}
\label{prodextD}
Pour tous
$\E \in D ^\mathrm{b} _\mathrm{coh} ( \smash{\D} ^\dag _{\X} (\hdag T _1) _{\Q} ) $,
$\FF\in 
D ^\mathrm{b} _\mathrm{coh} ( \smash{\D} ^\dag _{\Y} (\hdag T _2) _{\Q} )$, 
on dispose de l'isomorphisme canonique 
dans $D ^\mathrm{b} _\mathrm{coh} ( \smash{\D} ^\dag _{\ZZ} (\hdag T ) _{\Q} )$:
\begin{equation}
\label{prodextDisopre}
\left ( 
\smash{\D} ^\dag _{\ZZ} (\hdag T ) _{\Q} 
\otimes _{g ^{-1}\smash{\D} ^\dag _{\Y } (\hdag T _2) _{\Q} }
g ^{-1} \FF
\right )
\otimes _{f ^{-1} \smash{\D} ^\dag _{\X} (\hdag T _1) _{\Q} } ^{\L}
f ^{-1} \E
\riso
\E
\smash{\overset{\L}{\boxtimes}}   ^{\dag} _{\O _\S, T _1, T _2}
\FF.
\end{equation}
De plus, 
en prenant  $\E = \smash{\D} ^\dag _{\X} (\hdag T _1 ) _{\Q}$, 
l'isomorphisme induit 
\begin{equation}
\label{prodextDiso}
\smash{\D} ^\dag _{\ZZ} (\hdag T ) _{\Q} 
\otimes _{g ^{-1}\smash{\D} ^\dag _{\Y } (\hdag T _2) _{\Q} }
g ^{-1} \FF
\riso
\smash{\D} ^\dag _{\X} (\hdag T _1) _{\Q} 
\smash{\overset{\L}{\boxtimes}}   
^{\dag} _{\O _\S, T _1, T _2}
\FF
\end{equation}
est $(\smash{\D} ^\dag _{\ZZ} (\hdag T ) _{\Q} , f ^{-1} \smash{\D} ^\dag _{\X} (\hdag T _1) _{\Q} )$-bilinéaire.
\end{lemm}

\begin{proof}
On vérifie par fonctorialité que le faisceau d'anneaux $\smash{\widehat{\D}} _{\ZZ} ^{(m)} (T _1,T _2)$
est muni  d'une structure canonique de
$(\overset{^\mathrm{g}}{} \smash{\widehat{\D}} _{\ZZ} ^{(m)} (T _1,T _2), 
f ^{-1} (\smash{\widehat{\D}} _{\X} ^{(m)}  ( T _1)) \overset{^\mathrm{d}}{},
g ^{-1} ( \smash{\widehat{\D}} _{\Y} ^{(m)}  ( T _2))\overset{^\mathrm{d}}{})$-trimodule
et 
$\smash{\D} ^\dag _{\ZZ} (\hdag T ) _{\Q} $ de
$(\overset{^\mathrm{g}}{}\smash{\D} ^\dag _{\ZZ} (\hdag T ) _{\Q} , f ^{-1} \smash{\D} ^\dag _{\X} (\hdag T _1) _{\Q} \overset{^\mathrm{d}}{}, 
g ^{-1}\smash{\D} ^\dag _{\Y } (\hdag T _2) _{\Q} \overset{^\mathrm{d}}{})$-trimodule.
De plus, l'égalité et le morphisme canonique 
$f ^{(\bullet)*} (\smash{\widehat{\D}} _{\X} ^{(\bullet)}(T _1))
\smash{\overset{\L}{\otimes}}   ^{\dag} _{\O _{\ZZ,\Q}}
g ^{(\bullet)*} ( \smash{\widehat{\D}} _{\Y} ^{(\bullet)} (T _2))
=
\smash{\widehat{\D}} _{\ZZ} ^{(\bullet)} (T _1,T _2)
\to
\smash{\widehat{\D}} _{\ZZ} ^{(\bullet)} (T )$
sont 
$(\overset{^\mathrm{g}}{} \smash{\widehat{\D}} _{\ZZ} ^{(\bullet)} (T _1,T _2), 
f ^{-1} (\smash{\widehat{\D}} _{\X} ^{(\bullet)}  ( T _1)) \overset{^\mathrm{d}}{},
g ^{-1} ( \smash{\widehat{\D}} _{\Y} ^{(\bullet)}  ( T _2))\overset{^\mathrm{d}}{})$-trilinéaire pour les structures canoniques.
Comme ce dernier est en fait un isomorphisme (voir \ref{calcul-DT1T2-lim}), on en déduit l'avant-dernier isomorphisme:
\small
\begin{gather}
\notag
\left (\smash{\widehat{\D}} _{\ZZ} ^{(\bullet)} (T )
\otimes  ^{\L }_{g ^{-1} ( \smash{\widehat{\D}} _{\Y} ^{(\bullet)}   (T _2))}
g ^{-1} (\FF ^{(\bullet)} ) \right)
\otimes  ^{\L }_{f ^{-1} ( \smash{\widehat{\D}} _{\X} ^{(\bullet)}   (T _1))}
f ^{-1} (\E ^{(\bullet)} )
\riso
\left(\smash{\widehat{\D}} _{\ZZ} ^{(\bullet)} (T )
\smash{\overset{\L}{\otimes}}   ^{\dag}  _{g ^{-1} ( \smash{\widehat{\D}} _{\Y} ^{(\bullet)}   (T _2))}
g ^{-1} (\FF ^{(\bullet)} )\right)
\smash{\overset{\L}{\otimes}}   ^{\dag}  _{f ^{-1} ( \smash{\widehat{\D}} _{\X} ^{(\bullet)}   (T _1))}
f ^{-1} (\E ^{(\bullet)} )
\\
\notag
\liso 
\left(  (f ^{(\bullet)*} (\smash{\widehat{\D}} _{\X} ^{(\bullet)}(T _1))
\smash{\overset{\L}{\otimes}}   ^{\dag} _{\O _{\ZZ,\Q}}
g ^{(\bullet)*} ( \smash{\widehat{\D}} _{\Y} ^{(\bullet)} (T _2)) )
\smash{\overset{\L}{\otimes}}   ^{\dag}  _{g ^{-1} ( \smash{\widehat{\D}} _{\Y} ^{(\bullet)}   (T _2))}
g ^{-1} (\FF ^{(\bullet)} )\right)
\smash{\overset{\L}{\otimes}}   ^{\dag}  _{f ^{-1} ( \smash{\widehat{\D}} _{\X} ^{(\bullet)}   (T _1))}
f ^{-1} (\E ^{(\bullet)} )
\riso
\\
\label{prodextDisopre1}
\riso
f ^{(\bullet)*} (\E ^{(\bullet)} )
\smash{\overset{\L}{\otimes}}   ^{\dag} _{\O _{\ZZ,\Q}}
g ^{(\bullet)*} (\FF ^{(\bullet)}),
\end{gather}
\normalsize
le dernier se déduisant des isomorphismes 
$f ^{(\bullet)*} (\E ^{(\bullet)})
\riso
f ^{(\bullet)*} ( \smash{\widehat{\D}} _{\X} ^{(\bullet)}(T _1) )
\smash{\overset{\L}{\otimes}}   ^{\dag}  _{f ^{-1} ( \smash{\widehat{\D}} _{\X} ^{(\bullet)}  (T _1))}
f ^{-1} (\E ^{(\bullet)} )$
et
$g ^{(\bullet)*} (\FF ^{(\bullet)})
\riso
g ^{(\bullet)*} ( \smash{\widehat{\D}} _{\Y} ^{(\bullet)}(T _2) )
\smash{\overset{\L}{\otimes}}   ^{\dag}  _{g ^{-1} ( \smash{\widehat{\D}} _{\Y} ^{(\bullet)}  (T _2))}
g ^{-1} (\FF ^{(\bullet)} )$,
En appliquant le foncteur $\underrightarrow{\lim} $ à 
\ref{prodextDisopre1}, 
on obtient l'isomorphisme \ref{prodextDisopre}.
Lorsque $\E ^{(\bullet)}=\smash{\widehat{\D}} _{\X} ^{(\bullet)}  ( T _1)$, on vérifie que tous les isomorphismes
sont $(\overset{^\mathrm{g}}{} \smash{\widehat{\D}} _{\ZZ} ^{(\bullet)} (T _1,T _2), 
f ^{-1} (\smash{\widehat{\D}} _{\X} ^{(\bullet)}  ( T _1)) \overset{^\mathrm{d}}{})$-bilinéaire, ce qui donne la bilinéarité de 
\ref{prodextDiso}. 

\end{proof}

\begin{rema}
\label{rema-otimes-boxtimes}
Supposons $\X = \Y$, $T _1 = T _2$.
En notant $\delta \colon \X \hookrightarrow \X \times \X$ l'immersion diagonale, 
pour tous 
$\E ^{(\bullet)} ,~
\FF ^{(\bullet)} 
\in \underrightarrow{LD}  ^\mathrm{b} _{\Q, \mathrm{qc}}
(\overset{^\mathrm{g}}{} \smash{\widehat{\D}} _{\X} ^{(\bullet)} (T _1))$,
on dispose alors de l'isomorphisme canonique 
\begin{equation}
\label{otimes-boxtimes}
\E ^{(\bullet)}
\smash{\overset{\L}{\otimes}}   ^{\dag}
_{\O _{\X} (\hdag T _1) _\Q}\FF ^{(\bullet)}
\riso 
\delta 
^! (\E ^{(\bullet)}
\smash{\overset{\L}{\boxtimes}}   ^{\dag} _{\O _\S}
\FF^{(\bullet)}).
\end{equation}
Il en résulte,
pour tous
$\E ,\FF\in D ^\mathrm{b} _\mathrm{coh} ( \smash{\D} ^\dag _{\X} (\hdag T _1) _{\Q} ) $,
l'isomorphisme canonique
\begin{equation}
\label{lim-otimes-boxtimes}
\E
\smash{\overset{\L}{\otimes}}   ^{\dag}
_{\O  _{\X} (\hdag T  _1 ) _\Q}\FF
\riso 
\delta ^!
(\E
\smash{\overset{\L}{\boxtimes}}   ^{\dag} _{\O _\S, T _1, T _1}
\FF).
\end{equation}
\end{rema}

\section{Commutation au produit tensoriel des foncteurs de la forme $\sp _+$}

\subsection{Cas du cadre lisse}
\label{nota-cohPXT}
Nous garderons dans toute cette section les notations suivantes: 
soit $(\PP, T,X,Y)$ un cadre lisse
(voir les conventions de \ref{defi-cadre}).
Supposons de plus que $T \cap X$ soit un diviseur de $X$ (hypothèse non restrictive).
Nous avions appelé {\og cas de la compactification lisse\fg} l'ensemble de ces hypothèses.
Fixons $(\PP _{\alpha}) _{\alpha \in \Lambda}$ un recouvrement ouvert de $\PP$.
On note $\PP _{\alpha \beta}:= \PP _\alpha \cap \PP _\beta$,
$\PP _{\alpha \beta \gamma}:= \PP _\alpha \cap \PP _\beta \cap \PP _\gamma$,
$X _\alpha := X \cap P _\alpha$,
$X_{\alpha \beta } := X _\alpha \cap X _\beta$ et
$X_{\alpha \beta \gamma } := X _\alpha \cap X _\beta \cap X _\gamma $.
On suppose de plus que pour tout $\alpha\in \Lambda$, $X _\alpha$ est affine.
Comme $P$ est  séparé, pour tous $\alpha,\beta ,\gamma \in \Lambda$,
$X_{\alpha \beta }$ et $X_{\alpha \beta \gamma }$ sont donc affines.

Pour tout triplet $(\alpha, \, \beta,\, \gamma)\in \Lambda ^3$, choisissons
$\X _\alpha$ (resp. $\X _{\alpha \beta}$, $\X _{\alpha \beta \gamma}$)
des $\V$-schémas formels lisses relevant $X _\alpha$
(resp. $X _{\alpha \beta}$, $X _{\alpha \beta \gamma}$),
$u _{\alpha}\colon \X _{\alpha} \to \PP _{\alpha}$,
(resp. $p _1 ^{\alpha \beta}\colon \X  _{\alpha \beta} \rightarrow \X _{\alpha}$,
$p _2 ^{\alpha \beta}\colon \X  _{\alpha \beta} \rightarrow \X _{\beta}$)
des relèvements de
$X _{\alpha} \to P _{\alpha}$
(resp. $X  _{\alpha \beta} \rightarrow X _{\alpha}$,
$X  _{\alpha \beta} \rightarrow X _{\beta}$).
Rappelons que grâce à Elkik (\cite{elkik} de tels relèvements existent bien.

De même, pour tout triplet $(\alpha,\,\beta,\,\gamma )\in \Lambda ^3$, on choisit des relèvements
$p _{12} ^{\alpha \beta \gamma}\colon \X  _{\alpha \beta \gamma} \rightarrow \X  _{\alpha \beta} $,
$p _{23} ^{\alpha \beta \gamma}\colon \X  _{\alpha \beta \gamma} \rightarrow \X  _{\beta \gamma} $,
$p _{13} ^{\alpha \beta \gamma}\colon \X  _{\alpha \beta \gamma} \rightarrow \X  _{\alpha \gamma} $,
$p _1 ^{\alpha \beta \gamma}\colon \X  _{\alpha \beta \gamma} \rightarrow \X  _{\alpha} $,
$p _2 ^{\alpha \beta \gamma}\colon \X  _{\alpha \beta \gamma} \rightarrow \X  _{\beta} $,
$p _3 ^{\alpha \beta \gamma}\colon \X  _{\alpha \beta \gamma} \rightarrow \X  _{\gamma} $
induisant les morphismes canoniques au niveau des fibres spéciales.

Les deux catégories suivantes ont été construites dans \cite{caro-construction} (mais on modifie quelque peu les notations afin de mettre en exergue l'indépendance par rapport à $\PP$):
\begin{defi}
\label{donnee-recol}
Pour tout $\alpha \in \Lambda$, donnons-nous $\E _\alpha$,
un $(F\text{-})\smash{\D} ^{\dag} _{\X _{\alpha} } (\hdag T  \cap X _{\alpha}) _{\Q}$-module cohérent.
Une \textit{donnée de recollement} sur la famille $(\E _{\alpha})_{\alpha \in \Lambda}$ est 
la donnée pour tous $\alpha,\,\beta \in \Lambda$ d'un isomorphisme
$(F\text{-})\smash{\D} ^{\dag} _{\X _{\alpha \beta} }(\hdag T  \cap X _{\alpha \beta}) _{ \Q}$-linéaire de la forme
$ \theta _{  \alpha \beta} \colon   p _2  ^{\alpha \beta !} (\E _{\beta}) \riso p  _1 ^{\alpha \beta !} (\E _{\alpha}),$
ceux-ci vérifiant la condition de cocycle :
$\theta _{13} ^{\alpha \beta \gamma }=
\theta _{12} ^{\alpha \beta \gamma }
\circ
\theta _{23} ^{\alpha \beta \gamma }$,
où $\theta _{12} ^{\alpha \beta \gamma }$, $\theta _{23} ^{\alpha \beta \gamma }$
et $\theta _{13} ^{\alpha \beta \gamma }$ sont définis par les diagrammes commutatifs
\begin{equation}
  \label{diag1-defindonnederecol}
\xymatrix  @R=0,3cm {
{  p _{12} ^{\alpha \beta \gamma !} p  _2 ^{\alpha \beta !}  (\E _\beta )}
\ar[r] ^-{\tau} _{\sim}
\ar[d] ^-{p _{12} ^{\alpha \beta \gamma !} (\theta _{\alpha \beta})} _{\sim}
&
{p _2 ^{\alpha \beta \gamma!}  (\E _\beta )}
\ar@{.>}[d] ^-{\theta _{12} ^{\alpha \beta \gamma }}
\\
% 2-ième ligne
{ p _{12} ^{\alpha \beta \gamma !}  p  _1 ^{\alpha \beta !}  (\E _\alpha)}
\ar[r]^{\tau} _{\sim}
&
{p _1 ^{\alpha \beta \gamma!}(\E _\alpha),}
}
%
%
%2-ième diagramma :
\xymatrix  @R=0,3cm {
{  p _{23} ^{\alpha \beta \gamma !} p  _2 ^{\beta \gamma!}  (\E _\gamma )}
\ar[r] ^-{\tau} _{\sim}
\ar[d] ^-{p _{23} ^{\alpha \beta \gamma !} (\theta _{ \beta \gamma})} _{\sim}
&
{p _3 ^{\alpha \beta \gamma!}  (\E _\gamma )}
\ar@{.>}[d] ^-{\theta _{23} ^{\alpha \beta \gamma }}
\\
% 2-ième ligne
{ p _{23} ^{\alpha \beta \gamma !}  p  _1 ^{ \beta \gamma !}  (\E _\beta)}
\ar[r]^{\tau} _{\sim}
&
{p _2 ^{\alpha \beta \gamma!}(\E _\beta),}
}
%
%
%3-ième diagramme :
\xymatrix  @R=0,3cm {
{  p _{13} ^{\alpha \beta \gamma !} p  _2 ^{\alpha \gamma !}  (\E _\gamma )}
\ar[r] ^-{\tau} _{\sim}
\ar[d] ^-{p _{13} ^{\alpha \beta \gamma !} (\theta _{\alpha \gamma})} _{\sim}
&
{p _3 ^{\alpha \beta \gamma!}  (\E _\gamma )}
\ar@{.>}[d]^{\theta _{13} ^{\alpha \beta \gamma }}
\\
% 2-ième ligne
{ p _{13} ^{\alpha \beta \gamma !}  p  _1 ^{\alpha \gamma !}  (\E _\alpha)}
\ar[r]^{\tau} _{\sim}
&
{p _1 ^{\alpha \beta \gamma!}(\E _\alpha),}
}
\end{equation}
où les isomorphismes de la forme $\tau $ désignent les isomorphismes canoniques de recollement  (voir \cite[2.1.10]{caro-construction}).
\end{defi}

\begin{defi}
\label{251caroconstruction}
On construit la catégorie $(F\text{-})\mathrm{Isoc} ^{\dag \dag} (X,\, (\X _\alpha) _{\alpha \in \Lambda},\, T\cap X/K)$ de la manière
suivante :

\begin{itemize}
\item Un objet est une famille $(\E _\alpha) _{\alpha \in \Lambda}$
de $(F\text{-})\smash{\D} ^{\dag} _{\X _{\alpha} } (\hdag T  \cap X _{\alpha}) _{\Q}$-modules cohérents, $\O _{\X _{\alpha} } (\hdag T  \cap X _{\alpha}) _{\Q}$-cohérents, 
munie d'une donnée de recollement $ (\theta _{\alpha\beta}) _{\alpha ,\beta \in \Lambda}$.

\item Un morphisme
$((\E _{\alpha})_{\alpha \in \Lambda},\, (\theta _{\alpha\beta}) _{\alpha ,\beta \in \Lambda})
\rightarrow
((\E ' _{\alpha})_{\alpha \in \Lambda},\, (\theta '_{\alpha\beta}) _{\alpha ,\beta \in \Lambda})$
est une famille de morphismes $(F\text{-})\smash{\D} ^{\dag} _{\X _{\alpha} } (\hdag T  \cap X _{\alpha}) _{\Q}$-linéaire $f _\alpha\colon \E _\alpha \rightarrow \E '_\alpha$
commutant aux données de recollement,
i.e., telle que le diagramme suivant soit commutatif :
\begin{equation}
  \label{diag2-defindonnederecol}
\xymatrix  @R=0,3cm {
{ p _2  ^{\alpha \beta !} (\E _{\beta}) }
\ar[d] _-{p _2  ^{\alpha \beta !} (f _{\beta}) }
\ar[r] ^-{\theta _{\alpha\beta}} _{\sim}
&
{  p  _1 ^{\alpha \beta !} (\E _{\alpha}) }
\ar[d] ^-{p  _1 ^{\alpha \beta !} (f _{\alpha})}
\\
{p _2  ^{\alpha \beta !} (\E '_{\beta})  }
\ar[r]^{\theta '_{\alpha\beta}} _{\sim}
&
{ p  _1 ^{\alpha \beta !} (\E '_{\alpha})  .}
}
\end{equation}

\end{itemize}
\end{defi}

\begin{vide}
\label{const-sp+lisse}
On dispose de manière analogue à \ref{251caroconstruction} de la catégorie 
$(F\text{-})\mathrm{Isoc} ^{\dag} (X,\, (\X _\alpha) _{\alpha \in \Lambda},\, T\cap X/K)$ 
des familles d'objets de 
$\mathrm{Isoc} ^{\dag} (\X _\alpha , T \cap X _\alpha, X _\alpha /K)$ (voir la notation de \ref{nota-IsocDag})
munis d'une donnée de recollement.
Cette catégorie est construite dans \cite[2.5.6]{caro-construction} et y est notée 
$\mathrm{Isoc} ^\dag (Y,\,X,\, (\PP _\alpha) _{\alpha \in \Lambda}/K)$.
Nous préférons cette nouvelle notation pour mettre en évidence l'indépendance par rapport à $\PP$.
Par respectivement \cite[2.5.4, preuve de 2.5.7 et 2.5.9]{caro-construction}, 
les foncteurs canoniques de la forme
\begin{align}
\label{const-sp+lisse1}
\mathcal{L}oc&\colon \mathrm{Isoc} ^{\dag \dag} (\PP ,T, X /K)
\rightarrow
\mathrm{Isoc} ^{\dag \dag} (X,\, (\X _\alpha) _{\alpha \in \Lambda},\, T\cap X/K),
\\
\label{const-sp+lisse2}
\mathcal{L}oc&\colon 
\mathrm{Isoc} ^{\dag} (\PP ,T, X /K)
\to 
\mathrm{Isoc} ^{\dag} (X,\, (\X _\alpha) _{\alpha \in \Lambda},\, T\cap X/K),
\\
\label{const-sp+lisse3}
\sp _* &\colon 
\mathrm{Isoc} ^{\dag} (X,\, (\X _\alpha) _{\alpha \in \Lambda},\, T\cap X/K)
\to 
\mathrm{Isoc} ^{\dag \dag} (X,\, (\X _\alpha) _{\alpha \in \Lambda},\, T\cap X/K),
\end{align}
induits respectivement les foncteurs images inverses extraordinaires, images inverses et image directe par morphisme de spécialisation,
sont des équivalences de catégories. 
On dispose d'ailleurs 
d'un foncteur quasi-inverse canonique au foncteur \ref{const-sp+lisse1} dit de recollement
$\mathcal{R}ecol\colon 
\mathrm{Isoc} ^{\dag\dag } (X,\, (\X _\alpha) _{\alpha \in \Lambda},\, T\cap X/K)
\rightarrow
\mathrm{Isoc} ^{\dag \dag} (\PP ,T, X /K)$ (cela correspond à une {\og version recollée\fg} du théorème de Berthelot-Kashiwara).  
L'équivalence de catégories
$
\sp _{X \hookrightarrow \PP, T,+} \colon 
\mathrm{Isoc} ^{\dag}( \PP, T , X /K)
\cong 
\mathrm{Isoc} ^{\dag \dag}( \PP, T , X /K)$
de  \ref{eq-cat-spxPT+}
se construit alors, dans ce cas de la compactification lisse,
en posant 
\begin{equation}
\label{sp-cas-lisse}
\sp _{X \hookrightarrow \PP, T,+}
:=
\mathcal{R}ecol \circ \sp _* \circ \mathcal{L}oc.
\end{equation}

\end{vide}

\begin{prop}
\label{defotimescohalpha}
Avec les notations 
de \ref{251caroconstruction},
on définit le bifoncteur produit tensoriel
\begin{equation}
  \label{defotimescohalpha2}
-\otimes - \ : \
\mathrm{Isoc} ^{\dag \dag} (X,\, (\X _\alpha) _{\alpha \in \Lambda},\, T\cap X/K) \times
\mathrm{Isoc} ^{\dag \dag} (X,\, (\X _\alpha) _{\alpha \in \Lambda},\, T\cap X/K)
\rightarrow
\mathrm{Isoc} ^{\dag \dag} (X,\, (\X _\alpha) _{\alpha \in \Lambda},\, T\cap X/K),
\end{equation}
en posant,
pour tous $((\E _{\alpha}) _{\alpha \in \Lambda}, (\theta _{\alpha \beta }) _{\alpha, \beta \in \Lambda})$,
$((\E '_{\alpha}) _{\alpha \in \Lambda}, (\theta '_{\alpha \beta }) _{\alpha, \beta \in \Lambda})
\in \mathrm{Isoc} ^{\dag \dag} (X,\, (\X _\alpha) _{\alpha \in \Lambda},\, T\cap X/K)$,
$$((\E _{\alpha}) _{\alpha \in \Lambda}, (\theta _{\alpha \beta }) _{\alpha, \beta \in \Lambda})
\otimes
((\E '_{\alpha}) _{\alpha \in \Lambda}, (\theta '_{\alpha \beta }) _{\alpha, \beta \in \Lambda})
:=
(
(
\E _{\alpha}
\otimes _{\O _{\X _{\alpha}} (\hdag T  \cap X _{\alpha}) _{\Q}}
\E '_{\alpha}
) _{\alpha \in \Lambda}, (\theta ''_{\alpha \beta }  ) _{\alpha, \beta \in \Lambda}
),$$
où $\theta ''_{\alpha \beta }$ est l'unique morphisme induisant le diagramme commutatif :
\begin{equation}
  \label{deftheta''}
  \xymatrix  @R=0,4cm {
{ p  _2 ^{\alpha \beta !} (\E _{\beta}
\otimes _{\O _{\X _{\beta}} (\hdag T  \cap X _{\beta}) _{\Q}}
\E '_{\beta}) }
\ar[r] _-\sim
\ar@{.>}[d] _-\sim ^-{\theta ''_{\alpha \beta }}
&
{ p  _2 ^{\alpha \beta !} (\E _{\beta} )
\otimes _{\O _{\X _{\alpha \beta }} (\hdag T  \cap X _{\alpha \beta}) _{\Q}}
p  _2 ^{\alpha \beta !} (\E '_{\beta}) }
\ar[d] _-\sim ^-{\theta _{\alpha \beta } \otimes \theta '_{\alpha \beta }}
\\
{ p  _1 ^{\alpha \beta !} (\E _{\alpha}
\otimes _{\O _{\X _{\alpha}} (\hdag T  \cap X _{\alpha}) _{\Q}}
\E '_{\alpha}) }
\ar[r] _-\sim
&
{ p  _1 ^{\alpha \beta !} (\E _{\alpha} )
\otimes _{\O _{\X _{\alpha \beta }} (\hdag T  \cap X _{\alpha \beta}) _{\Q}}
p  _1 ^{\alpha \beta !} (\E '_{\alpha}) ,}
}
\end{equation}
dont les isomorphismes horizontaux découlent de la commutation des produits tensoriels aux images inverses extraordinaires (voir \ref{fg!prodtens})
et des isomorphismes 
\ref{rema15}.
\end{prop}

\begin{proof}
Pour vérifier que ce bifoncteur produit tensoriel a bien un sens,
il s'agit d'établir que les isomorphismes $\theta ''_{\alpha \beta }$ satisfont à
la condition de cocycle (voir \ref{donnee-recol}).
Considérons le diagramme commutatif :
\begin{equation}
\label{deftheta''cocycle}
  \xymatrix @R =0,4cm{
{p _2 ^{\alpha \beta \gamma!} ( \E _\beta \otimes \E' _\beta)}
\ar[r] _-\sim ^-{\tau}
\ar@{.>}[d] _-\sim ^-{\theta _{12} ^{\prime \prime \alpha \beta \gamma }}
&
{p _{12 } ^{\alpha \beta !} p  _2 ^{\alpha \beta !} (\E _{\beta}\otimes \E '_{\beta})}
\ar[r] _-\sim
\ar[d] _-\sim ^-{p _{12 } ^{\alpha \beta !}(\theta ''_{\alpha \beta })}
&
{p _{12 } ^{\alpha \beta !} p  _2 ^{\alpha \beta !} (\E _{\beta}) \otimes p _{12 } ^{\alpha \beta !}
p  _2 ^{\alpha \beta !} (\E '_{\beta})}
\ar[d] _-\sim ^-{p _{12 } ^{\alpha \beta !}(\theta _{\alpha \beta }) \otimes p _{12 } ^{\alpha \beta !}(\theta '_{\alpha \beta })}
&
{p _2 ^{\alpha \beta \gamma!} ( \E _\beta )\otimes p _2 ^{\alpha \beta \gamma!} ( \E' _\beta)}
\ar[d] _-\sim ^-{\theta _{12} ^{\alpha \beta \gamma }\otimes \theta _{12} ^{\prime \alpha \beta \gamma }}
\ar[l] _-\sim ^-{\tau\otimes \tau}
\\
{p _1 ^{\alpha \beta \gamma!} ( \E _\alpha \otimes \E' _\alpha)}
\ar[r] _-\sim ^-{\tau}
&
{p _{12 } ^{\alpha \beta !} p  _1 ^{\alpha \beta !} (\E _{\alpha}\otimes \E '_{\alpha})}
\ar[r] _-\sim
&
{p _{12 } ^{\alpha \beta !} p  _1 ^{\alpha \beta !} (\E _{\alpha})
\otimes p _{12 } ^{\alpha \beta !} p  _1 ^{\alpha \beta !} ( \E '_{\alpha})}
&
{p _1 ^{\alpha \beta \gamma!} ( \E _\alpha )\otimes p _1 ^{\alpha \beta \gamma} ( \E' _\alpha),}
\ar[l] _-\sim ^-{\tau\otimes \tau}
}
\end{equation}
où, d'après 
\ref{diag1-defindonnederecol}
et avec ses notations, les carrés de droite et de gauche sont commutatifs par définition.
En appliquant le foncteur $p _{12 } ^{\alpha \beta !}$ à 
\ref{deftheta''} puis par fonctorialité de la commutation du produit tensoriel aux images inverses extraordinaire,
on obtient le carré du milieu, qui est donc commutatif.  
Or, il découle de \ref{fg!prodtens-square} que les
isomorphismes composés horizontaux de \ref{deftheta''cocycle}
sont les isomorphismes canoniques de commutation des images inverses extraordinaires aux produits tensoriels. 
Avec les deux diagrammes analogues à \ref{deftheta''cocycle}, 
comme les familles d'isomorphismes $\theta _{\alpha \beta }$ et $\theta '_{\alpha \beta }$ satisfont aux conditions de cocycle,
il en est de même de
$\theta ''_{\alpha \beta }$. 
\end{proof}

\begin{lemm}
\label{loc-otimes-coh}
Avec les notations de 
\ref{nota-M-eq-isoc-lim},
soient $\E ^{(\bullet)},\,\E ^{\prime (\bullet)}\in \mathrm{Isoc} ^{(\bullet)}(\PP ,T, X /K)$
et $\E:=
\underrightarrow{\lim} \, 
\E ^{(\bullet)}$ et $\E':=
\underrightarrow{\lim} \, 
\E ^{\prime (\bullet)}$.
\begin{enumerate}
\item \label{loc-otimes-cohi)} 
Pour tout  $j \not =0$, 
$\mathcal{H} ^{j} 
(\E ^{(\bullet)} 
\smash{\overset{\L}{\otimes}} ^\dag _{\O _{\PP} (\hdag T) _{\Q}} 
\E ^{\prime (\bullet)} [d _{Y/P}])\riso 0$ dans 
$\underrightarrow{LD}  ^\mathrm{b} _{\Q, \mathrm{coh}}
(\smash{\widehat{\D}} _{\PP} ^{(\bullet)} (T))$
et 
$\mathcal{H} ^{0} (\E ^{(\bullet)} 
\smash{\overset{\L}{\otimes}} ^\dag _{\O _{\PP} (\hdag T) _{\Q}} 
\E ^{\prime (\bullet)} [d _{Y/P}])
\in
\mathrm{Isoc} ^{(\bullet)}( \PP, T, X/K)$.

\item \label{loc-otimes-cohii)}  
En identifiant $\E  
\smash{\overset{\L}{\otimes}} ^\dag _{\O _{\PP} (\hdag T) _{\Q}} \E '[d _{Y/P}]$
et
$\mathcal{H} ^{0} (\E  
\smash{\overset{\L}{\otimes}} ^\dag _{\O _{\PP} (\hdag T) _{\Q}} \E '[d _{Y/P}])$,
on dispose de l'isomorphisme canonique commutant à Frobenius :
\begin{equation}
\notag
\mathcal{L}oc (  \E  
\smash{\overset{\L}{\otimes}} ^\dag _{\O _{\PP} (\hdag T) _{\Q}} \E '[d _{Y/P}])
\riso \mathcal{L}oc (  \E ) \otimes \mathcal{L}oc (  \E ').
\end{equation}

\end{enumerate}

\end{lemm}

\begin{proof}
Posons $\E ^{\prime \prime (\bullet)} :=\E   ^{(\bullet)}
\smash{\overset{\L}{\otimes}} ^\dag _{\O _{\PP} (\hdag T) _{\Q}} 
 \E ^{\prime (\bullet)} [d _{Y/P}]$
 et
 $\E ^{\prime \prime } := 
\underrightarrow{\lim}  \E ^{\prime \prime (\bullet)} 
 \riso
 \E   
\smash{\overset{\L}{\otimes}} ^\dag _{\O _{\PP} (\hdag T) _{\Q}} 
 \E '[d _{Y/P}]  $.
 Traitons d'abord la partie \ref{loc-otimes-cohi)} du lemme.
Grâce à \cite[5.3.5.2]{caro-stab-sys-ind-surcoh}, les complexes 
$\E ^{(\bullet)} _\alpha:= \mathcal{H} ^{0} u _{\alpha } ^{(\bullet)!} (\E^{(\bullet)}   |_{\PP _\alpha} )$ et 
$\E ^{\prime (\bullet)}  _\alpha:=\mathcal{H} ^{0}  u _{\alpha } ^{(\bullet)!} (\E ^{\prime (\bullet)} |_{\PP _\alpha} )$ sont 
des objets de
$\underrightarrow{LM}   _{\Q, \mathrm{coh}}
(\overset{^\mathrm{g}}{} \smash{\widehat{\D}} _{\X _\alpha} ^{(\bullet)} (T   \cap X _\alpha))$.
Comme $\underrightarrow{\lim} \,  \E ^{(\bullet)} _\alpha$ 
et
$\underrightarrow{\lim} \,  \E ^{\prime (\bullet)} _\alpha$ sont
$\O _{\X _{\alpha}} (\hdag T \cap X _\alpha) _{\Q}$-cohérents
(voir la caractérisation \cite[2.5.10]{caro-construction} de 
$\mathrm{Isoc} ^{\dag \dag} (\PP ,T, X /K)$),
on a alors
$\E ^{(\bullet)} _\alpha, \,  \E ^{\prime (\bullet)} _\alpha
\in 
\mathrm{Isoc} ^{(\bullet)}(\PP  _\alpha ,T \cap X _\alpha, X  _\alpha /K)$.
On dispose alors des isomorphismes dans 
$\underrightarrow{LD} ^{\mathrm{b}}   _{\Q, \mathrm{coh}}
(\overset{^\mathrm{g}}{} \smash{\widehat{\D}} _{\PP _\alpha} ^{(\bullet)} (T \cap P_\alpha))$:
\begin{equation}
\label{loc-otimes-cohi)pre-iso1}
\E ^{\prime \prime (\bullet)} |_{\PP _\alpha} 
\underset{\cite[5.3.5.2]{caro-stab-sys-ind-surcoh}}{\riso} 
u _{\alpha +} ^{(\bullet)} (\E _\alpha ^{(\bullet)})  \smash{\overset{\L}{\otimes}} ^\dag 
_{\O _{\PP _\alpha}} 
u _{\alpha +} ^{(\bullet)} (\E _\alpha ^{\prime(\bullet)}) [d _{Y/P}]
\underset{\ref{u!u+=id-cons}}{\riso} 
u _{\alpha +} ^{(\bullet)} (\E _\alpha ^{(\bullet)} 
\smash{\overset{\L}{\otimes}} ^\dag 
_{\O _{\X _\alpha}} 
\E _\alpha ^{\prime(\bullet)})
\underset{\ref{form-proj-alpha}}{\riso}
u _{\alpha +} ^{(\bullet)}
(\E _\alpha ^{(\bullet)} \otimes _{\widehat{\B} _{\X _\alpha} ^{(\bullet)} (T   \cap X_\alpha)}\E _\alpha ^{\prime(\bullet)}).
\end{equation}
Comme, 
d'après \ref{form-proj-alpha},
$\E _\alpha ^{(\bullet)} \otimes _{\widehat{\B} _{\X _\alpha} ^{(\bullet)} (T   \cap X_\alpha)} \E _\alpha ^{\prime(\bullet)}
\in \mathrm{Isoc} ^{(\bullet)}(\PP  _\alpha ,T \cap X _\alpha, X  _\alpha /K)$, 
il en résulte que 
$\E ^{\prime \prime (\bullet)} \in \smash{\underrightarrow{LD}} ^{\mathrm{b}} _{\Q ,\mathrm{coh}}
(\smash{\widehat{\D}} _{\PP} ^{(\bullet)}(T))$.
En appliquant le foncteur $\underrightarrow{\lim}$
à l'isomorphisme \ref{loc-otimes-cohi)pre-iso1},
on en déduit que 
$\E '' \riso \mathcal{H} ^{0}(\E '') $
et
$\mathcal{H} ^{0}(\E '') |_{\PP _\alpha}\in
\mathrm{Isoc} ^{\dag \dag} (\PP _\alpha ,T\cap X _\alpha , X _\alpha /K)$.
Il en résulte que
$\mathcal{H} ^{0}(\E '')\in \mathrm{Isoc} ^{\dag \dag} (\PP ,T, X /K)$.
La pleine fidélité du foncteur $\underrightarrow{\lim}$ sur 
$ \smash{\underrightarrow{LD}} ^{\mathrm{b}} _{\Q ,\mathrm{coh}}
(\smash{\widehat{\D}} _{\PP} ^{(\bullet)}(T))$ nous permet de conclure la partie \ref{loc-otimes-cohi)} du lemme.

Vérifions à présent la partie \ref{loc-otimes-cohii)} du lemme. 
On dispose des isomorphismes:
\begin{equation}
  \label{pre-loc-otimes-coh-iso}
 u ^{(\bullet)!} _{\alpha } ( \E ^{\prime \prime (\bullet)}  |_{\PP _\alpha} )
 \underset{\ref{fg!prodtens}}{\riso} 
 u ^{(\bullet)!} _{\alpha }  (\E ^{(\bullet)} |_{\PP _\alpha} )
\smash{\overset{\L}{\otimes}} ^\dag _{\O _{\X _\alpha} (\hdag T\cap X_\alpha ) _{\Q}} u ^{(\bullet) !} _{\alpha } 
 (\E ^{\prime (\bullet)} |_{\PP _\alpha}) 
 \underset{\ref{form-proj-alpha}}{\riso}
 u ^{(\bullet)!} _{\alpha }  (\E ^{(\bullet)} |_{\PP _\alpha} )
\otimes _{\widehat{\B} _{\X _\alpha} ^{(\bullet)} (T   \cap X_\alpha)}
u ^{(\bullet) !} _{\alpha } 
 (\E ^{\prime (\bullet)} |_{\PP _\alpha}).
\end{equation}
En lui appliquant le foncteur $\underrightarrow{\lim} \, $, 
on obtient: 
\begin{equation}
  \label{loc-otimes-coh-iso}
 u ^{!} _{\alpha } ( \E ^{\prime \prime }  |_{\PP _\alpha} )
\riso 
 u ^{!} _{\alpha }  (\E ^{} |_{\PP _\alpha} )
\otimes _{\O _{\X _\alpha}(\hdag T \cap X _\alpha ) _{ \Q} } 
u ^{ !} _{\alpha } 
 (\E ^{\prime } |_{\PP _\alpha}).
\end{equation}
 Il reste à présent à vérifier que les isomorphismes \ref{loc-otimes-coh-iso}
 commutent aux isomorphismes de recollement.
 Considérons le diagramme suivant
 \small
\begin{equation}
\label{loc-otimes-coh-diagpre}
\xymatrix @R=0,5cm @C=0,4cm{
{p  _2 ^{\alpha \beta (\bullet)!}  u _{\beta } ^{(\bullet)!} ( \E ^{\prime \prime (\bullet)}  |_{\PP _\beta})}
\ar[d] _-\sim ^-{\tau}
\ar[r] _-\sim ^-{\ref{fg!prodtens}}
&
{ p  _2 ^{\alpha \beta (\bullet)!}  ( u ^{(\bullet)!} _{\beta }  (\E ^{(\bullet)} |_{\PP _\beta} )
\smash{\overset{\L}{\otimes}} ^\dag 
_{\O _{\X _{\beta}} }
 (\E ^{\prime (\bullet)} |_{\PP _\beta}) )}
\ar[r] _-\sim ^-{\ref{fg!prodtens}}
&
{p  _2 ^{\alpha \beta (\bullet)!}  ( u ^{(\bullet)!} _{\beta }  (\E ^{(\bullet)} |_{\PP _\beta} ))
\smash{\overset{\L}{\otimes}} ^\dag 
_{\O _{\X _{\alpha \beta}} } 
 p  _2 ^{\alpha \beta (\bullet)!}  ( u ^{(\bullet) !} _{\beta } (\E ^{\prime (\bullet)} |_{\PP _\beta}) )}
\ar[d] _-\sim ^-{\tau \otimes \tau }
\\
{ p  _1 ^{\alpha \beta (\bullet) !} u ^{(\bullet)!} _{\alpha } ( \E ^{\prime \prime (\bullet)}  |_{\PP _\alpha} )}
\ar[r] _-\sim ^-{\ref{fg!prodtens}}
&
{ p  _1 ^{\alpha \beta (\bullet)!}  ( u ^{(\bullet)!} _{\alpha }  (\E ^{(\bullet)} |_{\PP _\alpha} )
\smash{\overset{\L}{\otimes}} ^\dag 
_{\O _{\X _{\alpha}} }
 (\E ^{\prime (\bullet)} |_{\PP _\alpha}) )}
\ar[r] _-\sim ^-{\ref{fg!prodtens}}
&
{ p  _1 ^{\alpha \beta (\bullet)!}  ( u ^{(\bullet)!} _{\alpha }  (\E ^{(\bullet)} |_{\PP _\alpha} ))
\smash{\overset{\L}{\otimes}} ^\dag 
_{\O _{\X _{\alpha \beta}} }
p  _1 ^{\alpha \beta (\bullet)!}  ( u ^{(\bullet) !} _{\alpha } 
 (\E ^{\prime (\bullet)} |_{\PP _\alpha}) )}.
}
\end{equation}
 Modulo les identifications $p  _2 ^{\alpha \beta (\bullet)!}  \circ u _{\beta } ^{(\bullet)!} \riso (u _{\beta }\circ  p  _2 ^{\alpha \beta}  )^{(\bullet)!}$
 et
 $ p  _1 ^{\alpha \beta (\bullet)!}  \circ u _{\alpha } ^{(\bullet)!} \riso  (u _{\alpha }  \circ p  _1 ^{\alpha \beta} )^{(\bullet)!}$
 et par transitivité des isomorphismes \ref{fg!prodtens} (voir \ref{fg!prodtens-trans}),
 le rectangle \ref{loc-otimes-coh-diagpre} est de la forme \ref{fg!prodtens-square}. D'où sa commutativité.
En appliquant le foncteur $\underrightarrow{\lim} \, $
au rectangle \ref{loc-otimes-coh-diagpre},
modulo les isomorphismes canoniques de la forme \ref{rema15},
on obtient le contour du diagramme 
\begin{equation}
\label{loc-otimes-coh-diag}
\xymatrix {
{p  _2 ^{\alpha \beta !}  u _{\beta } ^{!} (\E  '' |_{\PP _\beta})}
\ar[d] _-\sim ^-{\tau}
\ar[r] _-\sim ^-{\ref{loc-otimes-coh-iso}}
&
\ar[r] _-\sim
{ p  _2 ^{\alpha \beta !}  ( u _{\beta } ^! (\E |_{\PP _\beta} )
\otimes _{\O _{\X _\beta}(\hdag T \cap X _\beta ) _{ \Q} } 
u _{\beta } ^! (\E '|_{\PP _\beta}))}
\ar[d] _-\sim
\ar[r] _-\sim
&
{p  _2 ^{\alpha \beta !}   u _{\beta } ^! (\E |_{\PP _\beta} )
\otimes _{\O _{\X _{\alpha \beta }} (\hdag T  \cap X _{\alpha \beta}) _{\Q}}
p  _2 ^{\alpha \beta !}   u _{\beta } ^! (\E |_{\PP _\beta} ) }
\ar[d] _-\sim ^-{\tau \otimes \tau }
\\
{ p  _1 ^{\alpha \beta !}  u _{\alpha } ^!
(\E '' |_{\PP _\alpha}) }
\ar[r] _-\sim ^-{\ref{loc-otimes-coh-iso}}
&
{ p  _1 ^{\alpha \beta !}  ( u _{\alpha } ^! (\E |_{\PP _\alpha} )
\otimes _{\O _{\X _\alpha}(\hdag T \cap X _\alpha ) _{ \Q} } 
u _{\alpha } ^! (\E '|_{\PP _\alpha}))}
\ar[r] _-\sim
&
{p  _1 ^{\alpha \beta !}   u _{\alpha } ^! (\E |_{\PP _\alpha} )
\otimes _{\O _{\X _{\alpha \beta }} (\hdag T  \cap X _{\alpha \beta}) _{\Q}}
p  _1 ^{\alpha \beta !}   u _{\alpha } ^! (\E |_{\PP _\alpha} ) }
}
\end{equation}
dont l'isomorphisme vertical du milieu est 
 défini de telle manière que le carré de droite soit commutatif.
 Par définition, cet isomorphisme définit la structure canonique de recollement
 de la famille $\left (u _{\alpha } ^! (\E |_{\PP _\alpha} )
\otimes _{\O _{\X _\alpha}(T \cap X _\alpha ) _{ \Q} } 
u _{\alpha } ^! (\E '|_{\PP _\alpha}) \right ) _\alpha$ (voir \ref{deftheta''}).  
Il s'agit ainsi de vérifier que le carré de gauche du diagramme \ref{loc-otimes-coh-diag} est commutatif, ce qui découle de celle de son contour et 
du carré de droite. 
\end{proof}

\begin{vide}
\label{nota-loc-otimes-isoc}
  Avec les notations de \ref{const-sp+lisse},
  de façon similaire à \ref{defotimescohalpha2}, on définit le bifoncteur produit tensoriel
  \begin{equation}
    \notag
    -\otimes -\ : \
    \mathrm{Isoc} ^{\dag} (X,\, (\X _\alpha) _{\alpha \in \Lambda},\, T\cap X/K) 
    \times
  \mathrm{Isoc} ^{\dag} (X,\, (\X _\alpha) _{\alpha \in \Lambda},\, T\cap X/K)
  \rightarrow
\mathrm{Isoc} ^{\dag} (X,\, (\X _\alpha) _{\alpha \in \Lambda},\, T\cap X/K).
  \end{equation}
  Comme pour \ref{loc-otimes-coh}, on construit alors, pour tous
  $E, E ' \in \mathrm{Isoc} ^{\dag} (\PP ,T, X /K)$, l'isomorphisme canonique commutant à Frobenius :
\begin{equation}
  \label{loc-otimes-isoc}
  \mathcal{L}oc (  E  \otimes _{j ^\dag \O _{] X [ _{\PP}} } E ' )
\riso \mathcal{L}oc (E) \otimes \mathcal{L}oc (E').
\end{equation}
\end{vide}

\begin{lemm}
\label{sp*otimes}
  Soient $((E _{\alpha})_{\alpha \in \Lambda},\, (\eta _{\alpha\beta}) _{\alpha ,\beta \in \Lambda}),
  ((E '_{\alpha})_{\alpha \in \Lambda},\, (\eta ' _{\alpha\beta}) _{\alpha ,\beta \in \Lambda})\in
\mathrm{Isoc} ^{\dag} (X,\, (\X _\alpha) _{\alpha \in \Lambda},\, T\cap X/K)$. 
  Avec les notations de \ref{const-sp+lisse} et \ref{nota-loc-otimes-isoc},
  on bénéficie de l'isomorphisme canonique commutant à Frobenius :
  $$\sp  _*  ( (E _{\alpha},\, \eta _{\alpha\beta})
  \otimes
(E '_{\alpha},\, \eta '_{\alpha\beta}) )
  \riso
\sp  _*  ( E _{\alpha},\, \eta _{\alpha\beta})
  \otimes
\sp _* ( E '_{\alpha},\, \eta '_{\alpha\beta}) .$$
\end{lemm}

\begin{proof}
D'après \cite[2.5.9]{caro-construction}, en notant
$(\E _{\alpha},\, \theta _{\alpha\beta}) :=\sp  _*  ( E _{\alpha},\, \eta _{\alpha\beta})$
et de même avec des primes, il revient au même d'établir l'isomorphisme
$$\sp  ^*  (( \E _{\alpha},\, \theta _{\alpha\beta})
  \otimes
(\E '_{\alpha},\, \theta '_{\alpha\beta}) )
  \riso
\sp  ^*  ( \E _{\alpha},\, \theta _{\alpha\beta})
  \otimes
\sp ^* ( \E '_{\alpha},\, \theta '_{\alpha\beta} ).$$
Pour cela, on vérifie que les isomorphismes canoniques
$\sp  ^*  (\E _{\alpha}   \otimes \E '_{\alpha})
  \riso
\sp  ^*  ( \E _{\alpha})
  \otimes
\sp ^* ( \E '_{\alpha})$
commutent aux isomorphismes de recollement respectifs.
\end{proof}

\begin{prop}
\label{lissestableotimes}
Pour tous $E ,E ' \in \mathrm{Isoc} ^{\dag} (\PP , T , X /K)$,
on bénéficie de l'isomorphisme canonique dans $ \mathrm{Isoc} ^{\dag \dag} (\PP , T , X /K)$:
\begin{equation}
  \label{lissestableotimes-iso}
  \sp _{X \hookrightarrow \PP, T ,+} (E  \otimes _{j ^\dag \O _{] X [ _{\PP}} } E ') \riso
  \sp _{X \hookrightarrow \PP, T,+} (E  )
  \smash{\overset{\L}{\otimes}} ^\dag _{\O _{\PP} (\hdag T) _{\Q}} 
  \sp _{X \hookrightarrow \PP, T,+} (E ' )[d _{Y/P}]
\end{equation}
qui commute à Frobenius.
\end{prop}

\begin{proof}
Grâce à la première assertion de \ref{loc-otimes-coh}, on vérifie que le terme de droite de \ref{lissestableotimes-iso} est aussi un élément
de $\mathrm{Isoc} ^{\dag \dag} (\PP , T , X /K)$.
Pour construire \ref{lissestableotimes-iso}, 
comme on dispose du foncteur 
pleinement fidèle $\mathcal{L}oc $, 
il est alors équivalent de construire un isomorphisme de la forme
\begin{equation}
  \label{28iso}
  \mathcal{L}oc \circ \sp _{X \hookrightarrow \PP, T ,+} (E  \otimes _{j ^\dag \O _{] X [ _{\PP}} } E ') \riso
\mathcal{L}oc (  \sp _{X \hookrightarrow \PP, T,+} (E  )
  \smash{\overset{\L}{\otimes}} ^\dag _{\O _{\PP} (\hdag T) _{\Q}} 
  \sp _{X \hookrightarrow \PP, T,+} (E ' )[d _{Y/P}] ).
\end{equation}
Or, on bénéficie de l'isomorphisme canonique
  $\mathcal{L}oc \circ \sp _{X \hookrightarrow \PP, T ,+}
  \riso
\sp _* \circ \mathcal{L}oc $ (voir la construction de $\sp _{X \hookrightarrow \PP, T ,+} $ 
rappelée dans \ref{const-sp+lisse}),
celui-ci commutant aux actions de Frobenius.
Par \ref{loc-otimes-coh}.\ref{loc-otimes-cohii)}, \ref{loc-otimes-isoc} et \ref{sp*otimes},
on obtient alors l'isomorphisme \ref{28iso}.

\end{proof}

\begin{nota}
\label{nota-P-P'lisse}
Soient $(\PP, T,X,Y)$ et $(\PP', T',X',Y')$ deux cadres tels que $X$ et $X'$ soient lisses, 
$T \cap X$ soit un diviseur de $X$ et $T '\cap X'$ soit un diviseur de $X'$.
On pose
$\PP'':=\PP \times \PP '$, $X '' := X \times X'$, $Y '' := Y \times Y'$,
$j\colon Y \subset X$, $j'\colon Y '\subset X'$ et $j''\colon Y ''\subset X''$
les inclusions canoniques. 
On note 
$\theta = (p, a, b)  \colon (\PP'', T'',X'',Y'') \to (\PP, T,X,Y)$
et
$\theta' = (p', a', b')  \colon (\PP'', T'',X'',Y'') \to
(\PP', T',X',Y')$
les morphismes de cadres 
induits par les projections canoniques, où
$T''= p ^{-1} (T) \cup p ^{\prime -1} (T')$. 
Soient $E  \in \mathrm{Isoc} ^{\dag}( \PP, T, X/K)$ et
$E' \in \mathrm{Isoc} ^{\dag }( \PP', T', X'/K)$.
Avec les notations de \ref{corr-424525}, 
on définit le bifoncteur 
$- \boxtimes -\colon 
\mathrm{Isoc} ^{\dag}( \PP, T, X/K) \times 
\mathrm{Isoc} ^{\dag }( \PP', T', X'/K) 
\rightarrow
\mathrm{Isoc} ^{\dag }( \PP'', T'', X''/K)$
en posant 
$$E \boxtimes E ' := \theta ^* (E ) \otimes _{j ^{\prime \prime \dag} \O _{] X'' [ _{\PP''} }} \theta ^{\prime *} (E ').$$

\end{nota}

\begin{prop}
\label{preboxplusproplisse}
Avec les notations \ref{nota-P-P'lisse}, on dispose de l'isomorphisme canonique commutant à Frobenius dans
$\mathrm{Isoc} ^{\dag \dag}( \PP'', T'', X''/K)$:
\begin{equation}
\label{preboxplusproplisse-iso}
\sp _{X'' \hookrightarrow \PP '' , T'', +} (E \boxtimes E ')
\riso
\sp _{X\hookrightarrow \PP , T , +} (E)
\smash{\overset{\L}{\boxtimes}}   ^{\dag} _{\O _\S, T, T '}\,
\sp _{X' \hookrightarrow \PP ', T' , +} (E').
\end{equation}
\end{prop}

\begin{proof}
Construisons d'abord le morphisme \ref{preboxplusproplisse-iso}.
Notons $\E :=\sp _{X\hookrightarrow \PP , T , +} (E)$ et
$\E ' := \sp _{X' \hookrightarrow \PP ', T' , +} (E')$.
D'après \ref{corr-424525},
il vient
\begin{gather}
\notag
\sp _{X'' \hookrightarrow \PP '' , T'', +} ( \theta ^* (E))
\riso
 (\hdag T '') \R \underline{\Gamma} ^\dag _{X ''} p ^! (\E  [-d _{Y '}] )
=: \theta ^{*} (\E) ,
\\
\notag
\sp _{X'' \hookrightarrow \PP '' , T'', +} ( \theta ^{\prime *} (E'))
\riso
(\hdag T '')  \R \underline{\Gamma} ^\dag _{X''} p ^{\prime !} (\E' [-d _{Y }])
=:
\theta ^{\prime *} (\E').
\end{gather}
Grâce à \ref{lissestableotimes}, 
on en déduit l'isomorphisme dans $D ^\mathrm{b} _\mathrm{coh} (\D ^\dag _{\PP ''} (\hdag T'') _\Q)$
\begin{equation}
\label{boxplusproplisse-iso1pre}
\sp _{X'' \hookrightarrow \PP '' , T'', +} (E \boxtimes E ')
\riso
\R \underline{\Gamma} ^\dag _{X ''} (\hdag T '')   p ^! (\E )
\smash{\overset{\L}{\otimes}}   ^{\dag} _{\O _{\PP''} (\hdag T '' ) _{\Q}}
\R \underline{\Gamma} ^\dag _{X ''} (\hdag T '')  p ^{\prime !} (\E') [-d _P - d _{P'}].
\end{equation}
Or, d'après  \ref{lemm-boxtimes-TouT'1},
on dispose de l'isomorphisme canonique 
\begin{equation}
\label{preboxplusproplisse-iso1}
\E
\smash{\overset{\L}{\boxtimes}}   ^{\dag} _{\O _\S, T, T '}\,
\E '
\riso 
(\hdag T '' ) p ^{!} (\E)
\smash{\overset{\L}{\otimes}}   ^{\dag} _{\O _{\PP''} (\hdag T '' ) _{\Q}}
(\hdag T '') p ^{\prime !} (\E) [-d _P - d _{P'}].
\end{equation}
Comme le terme de droite de \ref{preboxplusproplisse-iso1} est à support dans $X''$, 
il est donc isomorphe au terme de droite de \ref{boxplusproplisse-iso1pre}.
D'où le résultat. 
 \end{proof}

\subsection{Cas du cadre lisse en dehors du diviseur}

Soient $(\PP, T,X,Y)$ et $(\PP', T',X',Y')$ deux cadres lisses en dehors de leur diviseur (voir la définition \ref{defi-cadre}).
On pose
$\PP'':=\PP \times \PP '$, $X '' := X \times X'$, $Y '' := Y \times Y'$,
$j\colon Y \subset X$, $j'\colon Y '\subset X'$ et $j''\colon Y ''\subset X''$
les inclusions canoniques. 
On note 
$\theta = (p, a, b)  \colon (\PP'', T'',X'',Y'') \to (\PP, T,X,Y)$
et
$\theta' = (p', a', b')  \colon (\PP'', T'',X'',Y'') \to
(\PP', T',X',Y')$
les morphismes de cadres induits par les projections canoniques, où
$T''= p ^{-1} (T) \cup p ^{\prime -1} (T')$.

\begin{lemm}
\label{boxtimesDT}
Soient 
$\E ^{(\bullet)}
\in 
\mathrm{Isoc} ^{(\bullet)}(\PP ,T, X /K)$, 
$\E:=
\underrightarrow{\lim} \, 
\E ^{(\bullet)}$,
$\E ^{\prime (\bullet)}\in \mathrm{Isoc} ^{(\bullet)}(\PP ',T', X '/K)$,
$\E':=
\underrightarrow{\lim}~
\E ^{\prime (\bullet)}$.

\begin{enumerate}
\item Pour tout entier $j \not =0$, 
on a $\mathcal{H} ^{j} 
(\E ^{(\bullet)} 
\smash{\overset{\L}{\boxtimes}}   ^{\dag} _{\O _\S}\,
\E ^{\prime (\bullet)} )\riso 0$ dans 
$\underrightarrow{LD}  ^\mathrm{b} _{\Q, \mathrm{coh}}
(\smash{\widehat{\D}} _{\PP''} ^{(\bullet)} (T''))$.
De plus,
$\mathcal{H} ^{0} (\E ^{(\bullet)} 
\smash{\overset{\L}{\boxtimes}}   ^{\dag} _{\O _\S}\,
\E ^{\prime (\bullet)} )
\in
\mathrm{Isoc} ^{(\bullet)}( \PP'', T'', X''/K)$.

\item Pour tout entier $j \not =0$, on a l'isomorphisme $\mathcal{H} ^{j} (\E 
\smash{\overset{\L}{\boxtimes}}   ^{\dag}  _{\O _\S, T, T '}\,
\E ' )\riso 0$ 
et
$\mathcal{H} ^{0} (\E 
\smash{\overset{\L}{\boxtimes}}   ^{\dag} _{\O _\S, T, T '}\,
\E ')
\in
\mathrm{Isoc} ^{\dag \dag}( \PP'', T'', X''/K)$.
\end{enumerate}
\end{lemm}

\begin{proof}
Comme $\E ^{(\bullet)} 
\smash{\overset{\L}{\boxtimes}}   ^{\dag} _{\O _\S}\,
\E ^{\prime (\bullet)} )\in \underrightarrow{LD}  ^\mathrm{b} _{\Q, \mathrm{coh}}
(\smash{\widehat{\D}} _{\PP''} ^{(\bullet)} (T''))$, 
comme le foncteur 
$\underrightarrow{\lim}$
est pleinement fidèle sur cette dernière catégorie, 
on se ramène alors à vérifier la seconde partie du lemme. 
Notons $\U:= \PP \setminus T$, $\U':= \PP '\setminus T'$, $\U'':= \PP ''\setminus T''$.
Comme $\E 
\smash{\overset{\L}{\boxtimes}}   ^{\dag} _{\O _\S, T, T '}\,
\E' | \U''
= \E |\U
\smash{\overset{\L}{\boxtimes}}   ^{\dag} _{\O _\S}\,
\E' |\U'$, 
d'après le cas des cadres lisses déjà traités (voir \ref{preboxplusproplisse}), 
on a donc $\E 
\smash{\overset{\L}{\boxtimes}}   ^{\dag} _{\O _\S, T, T '}\,
\E' | \U''
\in \mathrm{Isoc} ^{\dag \dag}( \U'', Y''/K)$.
Il suffit alors de prouver que 
$\E
\smash{\overset{\L}{\boxtimes}}   ^{\dag} _{\O _\S, T, T '}\,
\E'$
est $\smash{\D} ^\dag _{\PP''} (\hdag T ''  ) _\Q$-surcohérent.
Comme $\E
\smash{\overset{\L}{\boxtimes}}   ^{\dag} _{\O _\S, T, T '}\,
\E'$
est $\smash{\D} ^\dag _{\PP''} (\hdag T ''  ) _\Q$-cohérent, 
l'idée est comme d'habitude de procéder par descente:
avec le lemme \cite[3.1.9]{caro-pleine-fidelite} 
(et la remarque \cite[3.1.10]{caro-pleine-fidelite}), 
on se ramène au cas où
$Y$ est dense dans $X$ avec $X$ intègre et 
où $Y'$ est dense dans $X'$ avec $X'$ intègre.
D'après \cite[5.3.1]{caro-pleine-fidelite} (on y utilise le cas particulier où le morphisme de départ est l'identité),
il résulte du théorème de désingularisation de de Jong qu'il existe
un morphisme de cadres 
de la forme 
$\alpha = (f, g, h)  \colon (\widetilde{\PP} , \widetilde{T} ,\widetilde{X} ,\widetilde{Y}) \to (\PP, T,X,Y)$
où
$\widetilde{X} $ est lisse, 
$\widetilde{T} = f ^{-1} (T)$ et $\widetilde{T} \cap \widetilde{X} $ est un diviseur à croisements normaux de $\widetilde{X}$,
$f$ est un morphisme propre et lisse de $\V$-schémas formels séparés et lisses,
$g$ est un morphisme propre, surjectif, génériquement fini et étale de $k$-variétés.
On pose
$\widetilde{\PP} '':= \widetilde{\PP}\times \PP '$, $\widetilde{X} '' := \widetilde{X} \times  X'$, $\widetilde{Y} '' := \widetilde{Y} \times Y '$,
$\alpha ''= (f'', g'', h'')  \colon (\widetilde{\PP} '', \widetilde{T} '',\widetilde{X} '',\widetilde{Y} '') \to (\PP'', T'',X'',Y'')$
le morphisme de cadres 
induit par $\alpha$,
où
$\widetilde{T} ''= f ^{\prime \prime -1} (T'') $. 
On note 
$\widetilde{\theta}  = (\widetilde{p}, \widetilde{a}, \widetilde{b})  \colon (\widetilde{\PP} '', \widetilde{T} '',\widetilde{X} '',\widetilde{Y} '')
 \to (\widetilde{\PP} , \widetilde{T},\widetilde{X} ,\widetilde{Y} )$
et
$\widetilde{\theta} ' = (\widetilde{p}', \widetilde{a}', \widetilde{b}') \colon (\widetilde{\PP} '', \widetilde{T} '',\widetilde{X} '',\widetilde{Y} '') \to
(\PP', T',X',Y')$
les morphismes de cadres 
induits par les projections canoniques.
Posons $\widetilde{\E}:= \R \underline{\Gamma} ^\dag _{\widetilde{X}} f _T ^! (\E )$. 
Par stabilité de la surcohérence, 
$\widetilde{\E} \in  D ^\mathrm{b} _\mathrm{surcoh} ( \smash{\D} ^\dag _{\widetilde{\PP}} (\hdag \widetilde{T}) _{\Q} )$.

D'après \cite[5.3.1]{caro-pleine-fidelite},
$\E $ est un facteur direct de
$ f _{T+}  (\widetilde{\E})$ dans 
$ D ^\mathrm{b} _\mathrm{coh} ( \smash{\D} ^\dag _{\PP} (\hdag T) _{\Q} )$.
Cela implique que 
$\E
\smash{\overset{\L}{\boxtimes}}   ^{\dag} _{\O _\S, T, T '}\,
\E'
$
est un facteur direct de 
$f _{T+} (\widetilde{\E})
\smash{\overset{\L}{\boxtimes}}   ^{\dag} _{\O _\S, T, T '}\,
\E'$ dans 
$ D ^\mathrm{b} _\mathrm{coh} ( \smash{\D} ^\dag _{\PP''} (\hdag T'') _{\Q} )$.
Or, d'après \ref{boxtimes-v+coh}, on dispose de l'isomorphisme
$f _{T,+} (\widetilde{\E})
\smash{\overset{\L}{\boxtimes}}   ^{\dag} _{\O _\S, T, T '}\,
\E'
\riso
f ''_{T'',+} \left (\widetilde{\E}
\smash{\overset{\L}{\boxtimes}}   ^{\dag} _{\O _\S, \widetilde{T}, T '}\,
\E' \right )$.
Par stabilité de la surcohérence par l'image directe d'un morphisme propre, 
on se ramène ainsi au cas où $X$ est lisse et $T \cap X$ est un diviseur à croisements normaux de $X$.
De même,
on se ramène ainsi au cas où $X'$ est lisse et $T' \cap X'$ est un diviseur à croisements normaux de $X'$. 
Dans ce cas, d'après \ref{preboxplusproplisse}, on obtient
$\E 
\smash{\overset{\L}{\boxtimes}}   ^{\dag} _{\O _\S, T, T '}\,
\E' \in \mathrm{Isoc} ^{\dag \dag}( \PP'', T'', X''/K)$.
\end{proof}

\begin{lemm}
\label{lemm-dagdagstblotimes}
Soient $\E ^{(\bullet)},\,\E ^{\prime (\bullet)}\in \mathrm{Isoc} ^{(\bullet)}(\PP ,T, X /K)$
et $\E:=
\underrightarrow{\lim} \, 
\E ^{(\bullet)}$ et $\E':=
\underrightarrow{\lim} \, 
\E ^{\prime (\bullet)}$.

\begin{enumerate}
\item Pour tout  $j \not =0$, 
$\mathcal{H} ^{j} 
(\E ^{(\bullet)} 
\smash{\overset{\L}{\otimes}} ^\dag _{\O _{\PP} (\hdag T) _{\Q}} 
\E ^{\prime (\bullet)} [d _{Y/P}])\riso 0$ dans 
$\underrightarrow{LD}  ^\mathrm{b} _{\Q, \mathrm{coh}}
(\smash{\widehat{\D}} _{\PP} ^{(\bullet)} (T))$
et 
$\mathcal{H} ^{0} (\E ^{(\bullet)} 
\smash{\overset{\L}{\otimes}} ^\dag _{\O _{\PP} (\hdag T) _{\Q}} 
\E ^{\prime (\bullet)} [d _{Y/P}])
\in
\mathrm{Isoc} ^{(\bullet)}( \PP, T, X/K)$.

\item On obtient alors 
$\E   
\smash{\overset{\L}{\otimes}} ^\dag _{\O _{\PP} (\hdag T) _{\Q} } 
\E ' [d _{Y/P}]
\riso
\mathcal{H} ^{0}(\E   
\smash{\overset{\L}{\otimes}} ^\dag _{\O _{\PP} (\hdag T) _{\Q} } 
\E ' [d _{Y/P}])
\in\mathrm{Isoc} ^{\dag \dag}( \PP, T , X /K)$.

\end{enumerate}

\end{lemm}

\begin{proof}
Notons 
$\delta \colon \PP \hookrightarrow \PP \times \PP$ l'immersion diagonale.
D'après \ref{otimes-boxtimes}
$\E ^{(\bullet)}
\smash{\overset{\L}{\otimes}}   ^{\dag}
_{\O _{\PP} (\hdag T ) _\Q}
\E ^{\prime (\bullet)}
\riso 
\delta 
^! (\E ^{(\bullet)}
\smash{\overset{\L}{\boxtimes}}   ^{\dag} _{\O _\S}
\E ^{\prime (\bullet)})$.
Il résulte du lemme \ref{boxtimesDT} que 
$\E ^{(\bullet)}
\smash{\overset{\L}{\boxtimes}}   ^{\dag} _{\O _\S}
\E ^{\prime (\bullet)}
\in 
\smash{\underrightarrow{LD}} ^{\mathrm{b}} _{\Q ,\mathrm{surcoh}}
(\smash{\widehat{\D}} _{\PP''} ^{(\bullet)}(T''))$.
Par stabilité de la surcohérence (voir \cite{caro-stab-sys-ind-surcoh}), cela entraîne que
$\E ^{(\bullet)}
\smash{\overset{\L}{\otimes}}   ^{\dag}
_{\O _{\PP} (\hdag T ) _\Q}
\E ^{\prime (\bullet)}
\in \smash{\underrightarrow{LD}} ^{\mathrm{b}} _{\Q ,\mathrm{surcoh}}
(\smash{\widehat{\D}} _{\PP} ^{(\bullet)}(T))$.
Il suffit par conséquent de vérifier le lemme en dehors de $T$ 
(e.g. on utilise \cite[4.3.12]{Be1}, la caractérisation 
des catégories de la forme $\mathrm{Isoc} ^{\dag \dag}$ et le fait que le foncteur $\underrightarrow{\lim}$ est pleinement fidèle
sur $\smash{\underrightarrow{LD}} ^{\mathrm{b}} _{\Q ,\mathrm{surcoh}}
(\smash{\widehat{\D}} _{\PP} ^{(\bullet)}(T))$
et commute aux foncteurs $\mathcal{H} ^{j}$).
Notons $\U:= \PP \setminus T$.
Il découle alors du cas traité dans
\ref{loc-otimes-coh} que
$\E   \smash{\overset{\L}{\otimes}} ^\dag _{\O _{\PP} (\hdag T) _{\Q} } \E '[d _{Y/P}] |_{\U}
= \E   | _{\U} \smash{\overset{\L}{\otimes}} ^\dag _{\O _{\U,\Q}} \E ' |_{\U} [d _{Y/U}]
\riso 
\mathcal{H} ^{0} 
(\E   | _{\U} \smash{\overset{\L}{\otimes}} ^\dag _{\O _{\U,\Q}} \E ' |_{\U} [d _{Y/U}])
\in \mathrm{Isoc} ^{\dag \dag}( \U, Y /K)$.
\end{proof}

\begin{lemm}
\label{sp+-comm-dual}
Soient $E ,E'\in \mathrm{Isoc} ^{\dag}( \PP, T, X/K)$. On dispose de l'isomorphisme canonique dans
$\mathrm{Isoc} ^{\dag \dag}( \PP, T, X/K)$:
$$  \sp _{X \hookrightarrow \PP, T ,+} (E  \otimes _{j ^\dag \O _{] X [ _{\PP}} } E ') 
\riso
 \sp _{X \hookrightarrow \PP, T,+} (E  )
\smash{\overset{\L}{\otimes}} ^\dag _{\O _{\PP} (\hdag T) _{\Q} } 
 \sp _{X \hookrightarrow \PP, T,+} (E ' )[d _{Y/P}]. $$
\end{lemm}

\begin{proof}
D'après \cite[4.2.2]{caro-pleine-fidelite},
on dispose de l'équivalence de catégories
$ \sp _{X \hookrightarrow \PP, T,+} \colon 
\mathrm{Isoc} ^{\dag}( \PP, T , X /K)
\cong 
\mathrm{Isoc} ^{\dag \dag}( \PP, T , X /K)$.
D'après le lemme \ref{lemm-dagdagstblotimes}, il en résulte que les deux termes 
sont bien dans $\mathrm{Isoc} ^{\dag \dag}( \PP, T, X/K)$.
Avec le lemme \cite[3.1.9]{caro-pleine-fidelite} 
(et la remarque \cite[3.1.10]{caro-pleine-fidelite}), 
on se ramène au cas où $X$ est intègre
et $Y$ est dense dans $X$.
Pour se ramener au cas où $X$ est lisse (cas déjà traité dans \ref{lissestableotimes}), on utilise le théorème de pleine fidélité 
\cite[3.4.2]{caro-pleine-fidelite} de la manière suivante :
d'après le théorème de désingularisation de de Jong 
(et avec l'aide du lemme \cite[3.5.9]{caro-pleine-fidelite}),
il existe un diviseur $\widetilde{T}$ contenant $T$
et un diagramme de la forme 
\begin{equation}
\label{formel-II-diag}
\xymatrix@R=0,3cm {
{\widetilde{Y}^{(0)}} \ar@{^{(}->}[r]  ^-{l ^{(0)}} \ar[d] ^-{c}
& 
{Y^{(0)}} \ar@{^{(}->}[r]  ^-{j ^{(0)}} \ar[d] ^-{b}
& 
{X^{(0)} } \ar@{^{(}->}[r] ^-{u^{(0)}} \ar[d] ^-{a}
& 
{\PP ^{(0)} }  \ar[d] ^-{f}
\\ 
{\widetilde{Y}} \ar@{^{(}->}[r] ^-{l}
&
{Y} \ar@{^{(}->}[r] ^-{j }
&
 {X } \ar@{^{(}->}[r] ^-{u} & {\PP  ,}
}
\end{equation}
où les deux carrés de gauche sont cartésiens, 
$\widetilde{Y} = X \setminus  \widetilde{T}$,
$X ^{(0)}$ est lisse, 
$f$ est un morphisme propre et lisse de $\V$-schémas formels séparés et lisses,
$a$ est un morphisme propre, surjectif de $k$-variétés, 
$b$ est un morphisme de $k$-variétés lisses, 
$c$ est un morphisme fini et étale, 
$l$, $l^{(0)}$, $j$ et $j^{(0)}$ sont des immersions ouvertes, $u$ et $u^{(0)}$ sont des immersions fermées,
$\widetilde{Y}$ est dense dans $Y$
et 
$\widetilde{Y} ^{(0)}$ est dense dans $Y ^{(0)}$.
Comme le foncteur $(\hdag \widetilde{T})$ est pleinement fidèle sur $\mathrm{Isoc} ^{\dag \dag}( \PP, T, X/K)$ (voir \cite[3.4.2]{caro-pleine-fidelite}), on se ramène au cas 
où $T = \widetilde{T}$. Notons $\theta:=(f,a,b)$ le morphisme de cadres.
On dispose des isomorphismes:
\begin{align}
\label{sp+-comm-dual-iso1}
\theta ^{*} \circ 
\sp _{X \hookrightarrow \PP, T ,+} (E  \otimes _{j ^\dag \O _{] X [ _{\PP}} } E ') 
\underset{\ref{iso-theta*com-sp}}{\riso} 
&
\sp _{X ^{(0)}\hookrightarrow \PP ^{(0)}, T ^{(0)},+} \circ \theta ^{*}  (  E \otimes _{j ^\dag \O _{] X [ _{\PP}} } E ' )
\\
\riso &
\sp _{X ^{(0)}\hookrightarrow \PP^{(0)}, T ^{(0)},+} \left (\theta ^{*}  (E ) \otimes _{j ^{(0)\dag} \O _{] X ^{(0)} [ _{\PP ^{(0)}}} } \theta ^{*}(E ') 
\right )
\\
\underset{\ref{lissestableotimes-iso}}{\riso}
&
\sp _{X ^{(0)}\hookrightarrow \PP^{(0)}, T ^{(0)},+}(\theta ^{*}  (E )) 
\smash{\overset{\L}{\otimes}} ^\dag _{\O _{\PP^{(0)}} (\hdag T ^{(0)}) _{\Q} } 
\sp _{X ^{(0)}\hookrightarrow \PP^{(0)}, T ^{(0)},+}  
( \theta ^{*}(E ') )[d _{Y ^{(0)}/P ^{(0)}}]
\\
\underset{\ref{iso-theta*com-sp}}{\riso} 
 &
\theta ^{*} \circ 
\sp _{X \hookrightarrow \PP, T,+} (E  )
 \smash{\overset{\L}{\otimes}} ^\dag _{\O _{\PP^{(0)}} (\hdag T ^{(0)}) _{\Q} } 
\theta ^{*} \circ  \sp _{X \hookrightarrow \PP, T,+} (E ' )[d _{Y ^{(0)}/P ^{(0)}}]
\\
 \riso
 &
\theta ^{*} \left ( 
\sp _{X \hookrightarrow \PP, T,+} (E  )
\smash{\overset{\L}{\otimes}} ^\dag _{\O _{\PP} (\hdag T) _{\Q} } 
 \sp _{X \hookrightarrow \PP, T,+} (E ' )[d _{Y/P}]
\right ).
\end{align}
En dehors du diviseur $T ^{(0)}$, cet isomorphisme 
est \ref{lissestableotimes-iso}.
Comme le foncteur $(\theta ^{*}, |\U)$ est pleinement fidèle (voir \cite[3.4.2]{caro-pleine-fidelite}), 
on en déduit la proposition.

\end{proof}

\begin{nota}
\label{nota-P-P'}
Soient $E  \in \mathrm{Isoc} ^{\dag}( \PP, T, X/K)$ et
$E' \in \mathrm{Isoc} ^{\dag }( \PP', T', X'/K)$.
Avec les notations de \ref{corr-424525},
on définit le bifoncteur 
$- \boxtimes -\colon 
\mathrm{Isoc} ^{\dag}( \PP, T, X/K) \times 
\mathrm{Isoc} ^{\dag }( \PP', T', X'/K) 
\rightarrow
\mathrm{Isoc} ^{\dag }( \PP'', T'', X''/K)$
en posant 
$$E \boxtimes E ' := \theta ^* (E ) \otimes _{j ^{\prime \prime \dag} \O _{] X'' [ _{\PP''} }} \theta ^{\prime *} (E ').$$

\end{nota}

\begin{prop}
\label{boxplusproplisse}
Avec les notations \ref{nota-P-P'}, on dispose de l'isomorphisme canonique 
dans $\mathrm{Isoc} ^{\dag \dag}( \PP'', T'', X''/K)$:
\begin{equation}
\label{boxplusproplisse-iso}
\sp _{X'' \hookrightarrow \PP '' , T'', +} (E \boxtimes E ')
\riso
\sp _{X\hookrightarrow \PP , T , +} (E)
\smash{\overset{\L}{\boxtimes}}   ^{\dag} _{\O _\S, T, T '}\,
\sp _{X' \hookrightarrow \PP ', T' , +} (E').
\end{equation}
\end{prop}

\begin{proof}
On construit le morphisme \ref{boxplusproplisse-iso} de manière analogue au début de la preuve de \ref{preboxplusproplisse} (on y remplace l'utilisation de la proposition \ref{lissestableotimes} par le lemme \ref{sp+-comm-dual}).
Grâce à \ref{boxtimesDT}, le terme de droite de \ref{boxplusproplisse-iso} est $\D ^\dag _{\PP ''} (\hdag T'') _\Q$-surcohérent.
Comme il en va de même de celui de gauche, pour vérifier 
que la flèche \ref{boxplusproplisse-iso} est un isomorphisme, 
il suffit
de l'établir en dehors de $T''$, ce qui nous ramène au cas où $T$ et $T'$ sont vides, et donc au cas où $X$ et $X'$ sont lisses.
On conclut en invoquant la proposition \ref{preboxplusproplisse} qui a résolu ce cas.
 \end{proof}

\begin{theo}
\label{stbintnlisse}
On suppose $\PP ' = \PP$ et que $(\PP, T  \cup T ', X \cap X', Y \cap Y')$ soit 
un cadre lisse en dehors du diviseur. 
Notons $i\colon (\PP, T  \cup T ', X \cap X', Y \cap Y') \to (\PP, T,X,Y)$, 
$i'\colon  (\PP, T  \cup T ', X \cap X', Y \cap Y') \to (\PP, T',X',Y')$ les morphismes canoniques de cadres
et $\widetilde{j}\colon Y \cap Y' \subset X \cap X'$ l'inclusion canonique. 

\begin{enumerate}
\item Pour tous 
$\E ^{(\bullet)}\in \mathrm{Isoc} ^{(\bullet)}(\PP ,T, X /K)$,
$\E ^{\prime (\bullet)}\in \mathrm{Isoc} ^{(\bullet)}(\PP,T', X '/K)$,
on bénéficie de l'isomorphisme canonique
dans 
$\mathrm{Isoc} ^{(\bullet)}( \PP,T  \cup T ', X \cap X'/K)$ de la forme:
\begin{equation}
\label{stbintnlisse-bullet}
\E ^{(\bullet)} 
\smash{\overset{\L}{\otimes}} ^\dag _{\O _{\PP}} 
\E ^{\prime (\bullet)} 
[d _{Y} + d _{Y'} - d _{Y\cap Y'} -d _{P}]
\riso
i ^{*}(\E ^{(\bullet)} )
\smash{\overset{\L}{\otimes}} ^\dag _{\O _{\PP}} 
i ^{\prime *}(\E ^{\prime (\bullet)} )
[d _{Y\cap Y'/P} ].
\end{equation}

\item Pour tous 
$E  \in \mathrm{Isoc} ^{\dag}( \PP, T, X/K)$ et
$E' \in \mathrm{Isoc} ^{\dag }( \PP', T', X'/K)$,
on dispose de l'isomorphisme canonique dans 
$\mathrm{Isoc} ^{\dag \dag}( \PP,T  \cup T ', X \cap X'/K)$:
\begin{gather}
\notag
\sp _{X \cap X'  \hookrightarrow \PP , T  \cup T ',+}   (i ^* (E) \otimes _{\widetilde{j} ^\dag\O _{ ] X \cap X ' [ _{\PP} } } i ^{\prime *} (E'))
\\
\label{stbintnisocorlisse}
\riso
(\hdag T \cup T') \circ \sp _{X \hookrightarrow \PP, T  ,+} (E)
\smash{\overset{\L}{\otimes}} ^\dag _{\O _{\PP} (\hdag T \cup T') _{\Q} }
(\hdag T \cup T') \circ \sp _{X ' \hookrightarrow \PP, T ' ,+} (E')
[d _{Y} + d _{Y'} - d _{Y\cap Y'} -d _{P}].
\end{gather}

\end{enumerate}
\end{theo}

\begin{proof}
Traitons d'abord \ref{stbintnlisse-bullet}.
Considérons les isomorphismes:
\begin{gather}
\notag
\E ^{(\bullet)} 
\smash{\overset{\L}{\otimes}} ^\dag _{\O _{\PP}} 
\E ^{\prime (\bullet)} 
\underset{\ref{otimesTT'}}{\riso}
(\hdag T \cup T') (\E ^{(\bullet)} )
\smash{\overset{\L}{\otimes}} ^\dag _{\O _{\PP} (\hdag T \cup T') _{\Q} }
(\hdag T \cup T')(\E ^{\prime (\bullet)} )
\riso
\\
\notag
\riso
\R\underline{\Gamma} ^\dag _{X \cap X'}
\circ (\hdag T \cup T') (\E ^{(\bullet)} )
\smash{\overset{\L}{\otimes}} ^\dag _{\O _{\PP} (\hdag T \cup T') _{\Q} }
\R\underline{\Gamma} ^\dag _{X \cap X'} \circ (\hdag T \cup T')(\E ^{\prime (\bullet)} )
=
i ^{*}(\E ^{(\bullet)} )
\smash{\overset{\L}{\otimes}} ^\dag _{\O _{\PP}} 
i ^{\prime *}(\E ^{\prime (\bullet)} )
[2 d _{Y\cap Y'}- (d _{Y} + d _{Y'})],
\end{gather}
le dernier isomorphisme résultant des isomorphismes
$\R\underline{\Gamma} ^\dag _{X} \E ^{(\bullet)} \riso \E ^{(\bullet)} $,
$\R\underline{\Gamma} ^\dag _{X'} \E ^{\prime (\bullet)} \riso \E ^{\prime(\bullet)} $
et 
$\R\underline{\Gamma} ^\dag _{X} \circ \R\underline{\Gamma} ^\dag _{X'}
\riso 
\R\underline{\Gamma} ^\dag _{X \cap X'}$.
Il découle de \ref{stabIsoc*inv-i-bullet} que
$i ^{*}(\E ^{(\bullet)} ),~
i ^{\prime *}(\E ^{\prime (\bullet)} )\in 
\mathrm{Isoc} ^{(\bullet)}( \PP,T  \cup T ', X \cap X'/K)$.
Le lemme \ref{lemm-dagdagstblotimes} 
nous permet de conclure la validation de \ref{stbintnlisse-bullet}. 
Enfin, l'isomorphisme
\ref{stbintnisocorlisse}
se construit en composant les isomorphismes ci-dessous:
\begin{gather}
\notag
\sp _{X \cap X'  \hookrightarrow \PP , T  \cup T ',+}   (i ^* (E) \otimes _{\widetilde{j} ^\dag \O _{ ] X \cap X ' [ _{\PP} } } i ^{\prime *} (E'))
\riso 
  \sp _{X \cap X'  \hookrightarrow \PP , T  \cup T ',+}   (i ^* (E) )
\smash{\overset{\L}{\otimes}} ^\dag _{\O _{\PP} (\hdag T \cup T') _{\Q} }
  \sp _{X \cap X'  \hookrightarrow \PP , T  \cup T ',+}   (i ^{\prime *} (E'))
[  d _{Y\cap Y'/P}]
  \\
  \notag
\underset{\ref{sp+-comm-dual}}{\riso}  
  i ^* \circ  \sp _{X \hookrightarrow \PP, T  ,+} (E)
\smash{\overset{\L}{\otimes}} ^\dag _{\O _{\PP} (\hdag T \cup T') _{\Q} }
i ^{\prime *} \circ \sp _{X ' \hookrightarrow \PP, T ' ,+} (E')
[  d _{Y\cap Y'/P}]
\\
\notag
=
\R\underline{\Gamma} ^\dag _{X \cap X'}
\circ (\hdag T \cup T') \circ \sp _{X \hookrightarrow \PP, T  ,+} (E)
\smash{\overset{\L}{\otimes}} ^\dag _{\O _{\PP} (\hdag T \cup T') _{\Q} }
\R\underline{\Gamma} ^\dag _{X \cap X'} \circ (\hdag T \cup T')\circ \sp _{X ' \hookrightarrow \PP, T ' ,+} (E')
[d _{Y} + d _{Y'} - d _{Y\cap Y'} -d _{P}]
\\
\notag
\riso
(\hdag T \cup T') \circ \sp _{X \hookrightarrow \PP, T  ,+} (E)
\smash{\overset{\L}{\otimes}} ^\dag _{\O _{\PP} (\hdag T \cup T') _{\Q} }
(\hdag T \cup T') \circ \sp _{X ' \hookrightarrow \PP, T ' ,+} (E')
[d _{Y} + d _{Y'} - d _{Y\cap Y'} -d _{P}],
\end{gather}
les raisons du dernier isomorphisme ayant déjà été données en début de preuve.
\end{proof}

\subsection{Cas des schémas faiblement formels}

\begin{vide}
\label{nota-spY->UdagT}
Soient $P ^\dag$ un $\V$-schéma formel faible séparé et lisse de fibre spéciale $P$, $T $ un diviseur de $P $, $U ^\dag$ l'ouvert de
$P ^\dag $ complémentaire de $T $, $U$ la fibre spéciale de $U ^\dag$, 
 $j$ : $U ^\dag \hookrightarrow P ^\dag$ l'immersion ouverte, 
  $v$ : $Y \hookrightarrow U $ une immersion fermée de $k$-schémas lisses.
On note $X $ l'adhérence schématique de $Y $ dans $P $. 
D'après \cite[6.1.8.1]{caro-pleine-fidelite},
on dispose du foncteur 
\begin{equation}
\label{spY->UdagT}
\sp _{Y \hookrightarrow U ^{\dag}, T,+}\colon
\mathrm{Isoc} ^\dag (Y/K) \to \mathrm{Isoc} ^{\dag\dag} (\PP, T, X /K),
\end{equation}
construit par recollement.
\end{vide}
Le théorème \cite[6.1.10]{caro-pleine-fidelite} reste valable 
sans l'hypothèse de lissité de $X$, i.e., on dispose du théorème suivant \ref{sp+essent}: 
\begin{theo}
 \label{sp+essent}
Avec les notations de \ref{nota-spY->UdagT},
on dispose du diagramme 
\begin{equation}
\xymatrix{
{\mathrm{Isoc} ^\dag (Y/K)}
\ar[rr] ^-{}
\ar[rd] _-{\sp _{Y  \hookrightarrow U ^{\dag},T  , + }} 
&&
{\mathrm{Isoc} ^{\dag} (\PP, T, X /K)} 
\ar[ld] ^-{ \sp _{X  \hookrightarrow \PP,T  , + }} _-{\cong}
\\ 
& 
{\mathrm{Isoc} ^{\dag\dag} (\PP, T, X /K)} 
& 
}
\end{equation}
commutatif à isomorphisme canonique près. 
\end{theo}

\begin{proof}
Soit $E \in \mathrm{Isoc} ^{\dag} (Y/K)$
et $E |(Y,X)$ l'objet de $\mathrm{Isoc} ^{\dag} (Y,X /K)=\mathrm{Isoc} ^{\dag} (\PP, T, X /K)$ induit. 
Avec le lemme \cite[3.1.9]{caro-pleine-fidelite} 
(et la remarque \cite[3.1.10]{caro-pleine-fidelite}), 
on se ramène au cas où 
$X$ est intègre et 
$Y$ est dense dans $X$.
Lorsque $X$ est lisse, cette proposition a été établie dans \cite[6.1.10]{caro-pleine-fidelite}.
Pour se ramener au cas où $X$ est lisse, comme pour la preuve de \ref{sp+-comm-dual}, on utilise le théorème de pleine fidélité 
\cite[3.4.2]{caro-pleine-fidelite} de la manière suivante:
d'après le théorème de désingularisation de de Jong 
(et avec l'aide du lemme \cite[3.5.9]{caro-pleine-fidelite}),
il existe un diviseur $\widetilde{T}$ contenant $T$
et un diagramme de la forme 
\begin{equation}
\label{formel-II-diag-dag}
\xymatrix@R=0,3cm {
{\widetilde{Y}^{(0)}} \ar@{^{(}->}[r]  ^-{l ^{(0)}} \ar[d] ^-{c}
& 
{Y^{(0)}} \ar@{^{(}->}[r]  ^-{j ^{(0)}} \ar[d] ^-{b}
& 
{X^{(0)} } \ar@{^{(}->}[r] ^-{u^{(0)}} \ar[d] ^-{a}
& 
{P ^{(0)\dag} }  \ar[d] ^-{f}
\\ 
{\widetilde{Y}} \ar@{^{(}->}[r] ^-{l}
&
{Y} \ar@{^{(}->}[r] ^-{j }
&
 {X } \ar@{^{(}->}[r] ^-{u} & {P ^\dag  ,}
}
\end{equation}
où les deux carrés de gauche sont cartésiens, 
$f$ est un morphisme propre et lisse de $\V$-schémas formels faibles séparés et lisses,
$a$ est un morphisme propre, surjectif de $k$-variétés avec $X ^{(0)}$ lisse, 
$c$ est un morphisme fini et étale, 
$l$, $l^{(0)}$, $j$ et $j^{(0)}$ sont des immersions ouvertes, $u$ et $u^{(0)}$ sont des immersions fermées,
$\widetilde{Y}$ est dense dans $X$ et
$\widetilde{Y}= Y \setminus \widetilde{T}$.
Comme le foncteur $(\hdag \widetilde{T})$ est pleinement fidèle sur la catégorie des isocristaux, on se ramène alors au cas 
où $T = \widetilde{T}$. 
Or, d'après respectivement \cite[6.1.9]{caro-pleine-fidelite} et 
\ref{corr-424525},
on bénéficie des isomorphismes canoniques:
\begin{gather}
\notag
\sp _{Y ^{(0) }  \hookrightarrow U ^{(0) \dag},T ^{(0) } , + }  \circ  b ^{*} (E) 
\riso
(\hdag T^{(0)} ) \circ \R \underline{\Gamma} ^{\dag} _{X ^{(0)} } \circ  f ^{!} \circ  \sp _{Y  \hookrightarrow U ^{\dag},T  , + } (E);
\\
\label{formel-II-diag-dag-iso}
\sp _{X ^{(0) } \hookrightarrow \PP^{(0) },T ^{(0) } , + } (a ^* (E |(Y,X)))
\riso
(\hdag T^{(0)} ) \circ \R \underline{\Gamma} ^{\dag} _{X ^{(0)} } \circ  f ^{!}
\sp _{X  \hookrightarrow \PP,T  , + } (E |(Y,X)).
\end{gather}
Comme $X ^{(0)}$ est lisse, les termes de gauche de \ref{formel-II-diag-dag-iso} sont canoniquement isomorphes.
Il en est donc de même des termes de droite. 
De plus, toujours d'après le cas lisse traité dans
\cite[6.1.10]{caro-pleine-fidelite}, on obtient l'isomorphisme canonique du milieu
\begin{equation}
\notag
 \sp _{Y  \hookrightarrow U ^{\dag},T  , + } (E) |\U 
 \riso 
  \sp _{Y  \hookrightarrow U ^{\dag},\emptyset  , + } (E) 
  \riso
 \sp _{X  \hookrightarrow \U,\emptyset  , + } (E |(Y,Y))  
  \riso
 \sp _{X  \hookrightarrow \PP,T  , + } (E |(Y,X))  |\U.
\end{equation}
Grâce au théorème de pleine fidélité \cite[3.4.2]{caro-pleine-fidelite}, on en déduit le résultat. 
\end{proof}

\begin{vide}\label{boxplus}
    Soient $P ^\dag$, $P ^{\prime \dag}$ deux $\V$-schémas formels faibles lisses et séparés,
    $T $ (resp. $T '$) un diviseur de $P $ (resp. $P '$), $U ^\dag$ (resp. $U ^{\prime \dag }$)
    l'ouvert de $P ^\dag $ (resp. $P ^{\prime \dag}$)
    complémentaire de $T $ (resp. $T '$), $j\colon U ^\dag \hookrightarrow P ^\dag$
    (resp. $j'\colon U ^{\prime \dag} \hookrightarrow P ^{\prime \dag}$)
    l'immersion ouverte et $v\colon Y \hookrightarrow U $
    (resp. $v '\colon Y ^{\prime } \hookrightarrow U ^{\prime}$)
    une immersion fermée de $k$-schémas lisses.
On note $P ^{\prime \prime \dag} := P ^\dag \times P ^{\prime \dag}$,
$U ^{\prime \prime \dag} := U ^\dag \times U ^{\prime \dag}$,
$T ''$ le diviseur réduit de $P '' $ d'espace topologique $P''  \setminus U ''$,
$Y ^{\prime \prime} := Y  \times Y ^{\prime}$,
$b \colon Y ^{\prime \prime} \rightarrow  Y $
et
$b '\colon Y ^{\prime \prime} \rightarrow  Y ^{\prime }$
les projections canoniques.

Soient $E  \in\mathrm{Isoc} ^\dag ( Y /K)$ et
$E' \in \mathrm{Isoc} ^\dag ( Y '/K)$.
On dispose des foncteurs canoniques
$b  ^* \colon \mathrm{Isoc} ^\dag ( Y /K) \rightarrow
\mathrm{Isoc} ^\dag ( Y '' /K)$
et
$b  ^{\prime *} \colon \mathrm{Isoc} ^\dag ( Y '/K) \rightarrow
\mathrm{Isoc} ^\dag ( Y '' /K)$ (voir \cite[2.3.6]{Berig} et \cite[1.4.1]{caro-2006-surcoh-surcv})).
Le produit tensoriel externe de $E$ et $E'$ est défini en posant
$E \boxtimes E ' := b  ^* (E ) \otimes  b  ^{\prime *} (E ')$, ce qui donne le bifoncteur 
$$-\boxtimes - \colon 
\mathrm{Isoc} ^\dag ( Y /K) 
\times 
\mathrm{Isoc} ^\dag ( Y '/K) 
\to 
\mathrm{Isoc} ^\dag ( Y''/K) .$$
\end{vide}

\begin{prop}\label{boxplusprop}
  Avec les notations \ref{boxplus}, on dispose d'un isomorphisme canonique
\begin{equation}
  \label{boxplusprop-iso}
\sp _{Y ^{\prime \prime } \hookrightarrow U ^{\prime \prime \dag} , T''  +}
(E \boxtimes E ') 
\riso
\sp _{Y  \hookrightarrow U ^{\dag} , T  +} (E)
\smash{\overset{\L}{\boxtimes}}   ^{\dag} _{\O _\S, T, T '}\,
\sp _{Y ^{\prime} \hookrightarrow U ^{\prime \dag} , T ' +} (E').
\end{equation}
\end{prop}

\begin{proof}
Notons 
$X $ (resp. $X'$, resp. $X''$) l'adhérence schématique de $Y $ dans $P $ (resp. de $Y'$ dans $P'$, resp. de $Y''$ dans $P''$). 
Grâce à \ref{boxplusproplisse}, le carré de droite du diagramme canonique
\small
\begin{equation}
\xymatrix{
{\mathrm{Isoc} ^\dag ( Y /K) 
\times 
\mathrm{Isoc} ^\dag ( Y '/K) }
\ar[d] ^-{-\boxtimes - }
\ar[r] ^-{}
& 
{\mathrm{Isoc} ^{\dag} (\PP, T, X /K) \times \mathrm{Isoc} ^{\dag} (\PP', T', X' /K)} 
\ar[d] ^-{-\boxtimes - }
\ar[r] ^-{\sp _{ + } \times \sp _{+ }}
&
{\mathrm{Isoc} ^{\dag \dag} (\PP, T, X /K) \times \mathrm{Isoc} ^{\dag} (\PP', T', X' /K)} 
\ar[d] ^-{-\boxtimes - }
\\ 
{\mathrm{Isoc} ^\dag ( Y''/K)}  
\ar[r] ^-{} 
& 
{\mathrm{Isoc} ^{\dag} (\PP'', T'', X'' /K)} 
\ar[r] ^-{\sp _{X '' \hookrightarrow \PP'',T '' , + }}
& 
{\mathrm{Isoc} ^{\dag \dag} (\PP'', T'', X'' /K).} 
}
\end{equation}
\normalsize
est commutatif, à isomorphisme canonique près. Comme le carré de gauche l'est aussi, 
le théorème \ref{sp+essent} nous permet de conclure. 
\end{proof}

\begin{prop}\label{stbintn}
  Avec les notations de \ref{boxplus}, on suppose de plus $P^\dag   = P^{\prime \dag} $ et
$Y  \cap Y '$ lisse (e.g., $Y  = Y' $). 
En notant
$i  ^* \colon \mathrm{Isoc} ^\dag ( Y  /K)\rightarrow
\mathrm{Isoc} ^\dag ( Y \cap Y ' /K) $
et
$i  ^{\prime *} \colon \mathrm{Isoc} ^\dag ( Y ' /K) \rightarrow
\mathrm{Isoc} ^\dag (  Y \cap Y '  /K) $ 
les foncteurs canoniques,
on dispose alors de l'isomorphisme canonique :
\begin{gather}
  \sp _{Y  \cap Y ^{\prime }  \hookrightarrow U ^{\dag} \cap U ^{\prime \dag} , T  \cup T ' ,+}
  (i ^* (E) \otimes i ^{\prime *} (E'))
  \\
  \label{stbintnisocor}
    \riso
(\hdag T \cup T') \circ \sp _{Y  \hookrightarrow U ^{\dag} , T  +}   (E )
\smash{\overset{\L}{\otimes}} ^\dag _{\O _{\PP} (\hdag T \cup T') _{\Q} }
(\hdag T \cup T') \circ \sp _{Y ^{\prime} \hookrightarrow U ^{\prime \dag} , T '  +} (E')
[d _{Y} + d _{Y'} - d _{Y\cap Y '} -d _{P}].
\end{gather}
\end{prop}

\begin{proof}
Cela résulte de \ref{stbintnlisse} et de \ref{sp+essent}.
\end{proof}

\section{Stabilité par produits tensoriels}

\subsection{Dévissabilité en isocristaux surconvergents pour les systèmes inductifs}

Soit $\PP$ un $\V$-schéma formel séparé et lisse. 
Lorsque $\PP$ est propre, nous avions défini 
dans \cite[7.3.2]{caro-2006-surcoh-surcv}
une sous-catégorie pleine de
$\smash{\underrightarrow{LD}} ^{\mathrm{b}} _{\Q ,\mathrm{qc}}
(\smash{\widehat{\D}} _{\PP} ^{(\bullet)})$
des complexes dévissables en isocristaux surconvergents.
Grâce à \cite{caro-pleine-fidelite}, 
on peut étendre naturellement cette définition en général, i.e. 
sans supposer $\PP$ propre (de manière analogue au paragraphe \ref{rappel-dev-coh}).
Nous en profitons pour introduire cette notion de dévissabilité de manière plus agréable
mais équivalente (voir la remarque \ref{rema-dev=dev-plong})
ainsi que quelques améliorations.

\begin{defi}
Soit $Y$ une sous-variété de $P$.
On dit que $Y $ est {\og $d$-plongeable dans $P$\fg} s'il existe un diviseur $T$ de $P$ tel que $Y$ soit fermé dans $P\setminus T$.
On remarque que cela équivaut à supposer qu'il existe un cadre 
de la forme $(\PP, T,X,Y)$.
De plus, 
si $Y', Y$ sont deux sous variétés $d$-plongeables dans $P$, 
alors $Y \cap Y'$ est $d$-plongeable dans $P$. 
\end{defi}

\begin{defi}
Soit $Y$ une sous-variété de $P$.

\begin{enumerate}
\item Une {\og stratification de $Y$\fg} 
est la donnée
(pour un certain entier $r\geq 1$) 
 de $r$ sous-variétés  
$(Y _1, Y _2,\dots, Y _{r})$ 
telle que, en posant $Y _{0}:=\emptyset$, 
pour tout entier $i$ vérifiant $1\leq i \leq r-1$, la variété $Y _{i}$ soit un ouvert de $Y  \setminus  (\cup  _{0\leq j \leq i-1}Y _{j})$ et telle que 
$Y _{r}= Y\setminus  (\cup  _{0\leq j \leq r-1}Y _{j})$.
Autrement dit, on dispose de la somme directe  
$Y =\sqcup _{i=1,\dots , r} Y _i$
telle que,
pour tout $1\leq i \leq r-1$,
la variété $Y _i$ est un ouvert de  
$\sqcup _{j=i,\dots , r} Y _j$. 
En prenant garde à l'ordre, 
on dira aussi qu'une telle décomposition
$Y =\sqcup _{i=1,\dots , r} Y _i$
est une stratification. 

\item Soit $Y =\sqcup _{i=1,\dots , r} Y _i$ une stratification.
On dira que $Y =\sqcup _{i=1,\dots , r} Y _i$ ou $(Y _1, Y _2,\dots, Y _{r})$ 
est une
{\og $d$-stratification de $Y$ dans $P$\fg} 
(resp. une stratification lisse, resp. une 
{\og $d$-stratification lisse de $Y$ dans $P$\fg}) 
si, pour tout $1\leq i \leq r$, la variété $Y _i$ est $d$-plongeable dans $P$
(resp. est lisse, resp. est lisse et $d$-plongeable dans $P$).

\end{enumerate}
\end{defi}

\begin{rema}
\label{rema-dev-tridist}
Soient $Y$ une sous-variété de $P$ 
et $Y =\sqcup _{i=1,\dots , r} Y _i$ une stratification.
Pour tout,
$\E ^{(\bullet)} \in \smash{\underrightarrow{LD}} ^{\mathrm{b}} _{\Q ,\mathrm{qc}}
(\smash{\widehat{\D}} _{\PP} ^{(\bullet)})$,
pour tout $1\leq i \leq r$, 
on dispose du triangle distingué de localisation
\begin{equation}
\notag
\R\underline{\Gamma} ^\dag _{\sqcup _{j=i+1,\dots , r} Y _j} (\E ^{(\bullet)}) \rightarrow \R\underline{\Gamma} ^\dag _{\sqcup _{j=i,\dots , r} Y _j} (\E ^{(\bullet)})
\rightarrow \R\underline{\Gamma} ^\dag _{Y _i} (\E ^{(\bullet)}) \rightarrow +1
\end{equation}

\end{rema}

\begin{nota}
Soit $Y$ une sous-variété lisse $d$-plongeable dans $P$. 
Choisissons $X$ un sous-schéma fermé de $P$, $T$ un diviseur de $P$ tels que $Y = X \setminus T$.
On note 
$\smash{\underrightarrow{LD}} ^{\mathrm{b}} _{\Q ,\mathrm{isoc},X}
(\smash{\widehat{\D}} _{\PP} ^{(\bullet)} (T))$
la sous-catégorie strictement pleine triangulée de 
$\smash{\underrightarrow{LD}} ^{\mathrm{b}} _{\Q ,\mathrm{coh}}
(\smash{\widehat{\D}} _{\PP} ^{(\bullet)} (T))$ des 
complexes 
 $\E ^{(\bullet)}$ à support dans $X$
tels que, pour tout $j \in \Z$, on ait
$$\mathcal{H} ^j (\E ^{(\bullet)} )
\in 
\mathrm{Isoc} ^{(\bullet)} ( \PP, T ,X /K).$$
Cette catégorie ne dépend  ni du choix du diviseur $T$, ni de celui du sous-schéma fermé $X$ tels que
$Y = X \setminus T$.
L'indépendance par rapport à $X$  résulte du fait qu'un complexe
$\E ^{(\bullet)}\in \smash{\underrightarrow{LD}} ^{\mathrm{b}} _{\Q ,\mathrm{coh}}
(\smash{\widehat{\D}} _{\PP} ^{(\bullet)} (T))$ 
est à support dans $X$ si et seulement si
$\R\underline{\Gamma} ^\dag _{Y}  (\E ^{(\bullet)}) \riso \E ^{(\bullet)}$. 
Soit $T'$ un diviseur tel que $X \setminus T'= Y$.
Quitte à considérer $T \cup T'$, on peut supposer $T \subset T'$.
Or, les foncteurs exacts $(\hdag T',T)$ et $oub _{T,T'}$ induisent des équivalences quasi-inverses entre
$\mathrm{Isoc} ^{\dag \dag} ( \PP, T ,X /K)$ et 
$\mathrm{Isoc} ^{\dag \dag} ( \PP, T ',X /K)$.
On déduit alors de \cite[3.5.2]{caro-stab-sys-ind-surcoh} 
que l'on dispose de la factorisation
$oub _{T,T'}
\colon 
\mathrm{Isoc} ^{(\bullet)} ( \PP, T',X /K)\to 
\mathrm{Isoc} ^{(\bullet)} ( \PP, T ,X /K)$.
Comme la factorisation induite par $(\hdag T',T)$ est triviale, 
il en résulte que 
les foncteurs $(\hdag T',T)$ et $oub _{T,T'}$ induisent des équivalences quasi-inverses entre
$\mathrm{Isoc} ^{(\bullet)} ( \PP, T,X /K)$
et
$\mathrm{Isoc} ^{(\bullet)} ( \PP, T',X /K)$.
Comme ces foncteurs $(\hdag T',T)$ et $oub _{T,T'}$ sont exacts, 
il en est de même pour 
$\smash{\underrightarrow{LD}} ^{\mathrm{b}} _{\Q ,\mathrm{isoc},X}
(\smash{\widehat{\D}} _{\PP} ^{(\bullet)} (T))$
et
$\smash{\underrightarrow{LD}} ^{\mathrm{b}} _{\Q ,\mathrm{isoc},X}
(\smash{\widehat{\D}} _{\PP} ^{(\bullet)}(T'))$.
On pourra donc la noter simplement 
$\smash{\underrightarrow{LD}} ^{\mathrm{b}} _{\Q ,\mathrm{isoc},Y}
(\smash{\widehat{\D}} _{\PP} ^{(\bullet)})$. 

\end{nota}

\begin{rema}
Soit $Y$ une sous-variété lisse $d$-plongeable dans $P$. 
Les objets de $\smash{\underrightarrow{LD}} ^{\mathrm{b}} _{\Q ,\mathrm{isoc},Y}
(\smash{\widehat{\D}} _{\PP} ^{(\bullet)})$ sont appelés, 
conformément à la terminologie de \cite[3.2.2]{caro-2006-surcoh-surcv} 
(voir les précisions de la remarque \ref{rema-precision-in-isoc} concernant cette terminologie)
ceux dont les espaces de cohomologie sont des isocristaux surconvergents sur $(Y,X)/K$.
On bénéficie de plus de l'équivalence de catégories:
\begin{equation}
  \label{dev=qwcdev+cohlim-pre}
\underrightarrow{\lim} \colon
\smash{\underrightarrow{LD}} ^{\mathrm{b}} _{\Q ,\mathrm{isoc},Y}
(\smash{\widehat{\D}} _{\PP} ^{(\bullet)})
\cong
D ^\mathrm{b} _\mathrm{isoc} (\PP, T, X/K)
\end{equation}
où la catégorie à droite a été définie dans \cite[1.2.5]{caro-image-directe}
et se définit de manière analogue.
Cette équivalence 
caractérise d'ailleurs la sous-catégorie 
$D ^\mathrm{b} _\mathrm{isoc} (\PP, T, X/K)$ strictement pleine de 
$D ^\mathrm{b} _\mathrm{coh} (\D ^{\dag}_{\PP} (\hdag T) _\Q)$.

\end{rema}

\begin{defi}
\label{coro-caradev}
Soit $Y$ une sous-variété de $P$.
Soit $\E ^{(\bullet)} \in \underrightarrow{LD}  ^\mathrm{b} _{\Q, \mathrm{qc}}
(\overset{^\mathrm{g}}{} \smash{\widehat{\D}} _{\PP} ^{(\bullet)})$.

\begin{enumerate}
\item Le complexe 
$\E ^{(\bullet)} $
 {\og se dévisse au dessus de $Y$ en isocristaux surconvergents\fg}
s'il existe une $d$-stratification lisse de $Y$ dans $P$ de la forme
$Y =\sqcup _{i=1,\dots , r} Y _i$
telle que, 
pour tout $i=1,\dots, r$,
on ait 
$\R\underline{\Gamma} ^\dag _{Y _i}  (\E ^{(\bullet)}) \in 
\smash{\underrightarrow{LD}} ^{\mathrm{b}} _{\Q ,\mathrm{isoc},Y _i}
(\smash{\widehat{\D}} _{\PP} ^{(\bullet)})$.
On dira aussi que le complexe {\og $\E ^{(\bullet)}$ se dévisse au-dessus de
la $d$-stratification lisse 
$Y =\sqcup _{i= 0,\dots , r-1} Y _i$ dans $P$ en isocristaux surconvergents\fg}.

\item Lorsque $Y=P$, on dit simplement que $\E ^{(\bullet)}$ se dévisse en isocristaux surconvergents. 
Pour tout diviseur $T$ de $P$,
on notera $\underrightarrow{LD}  ^\mathrm{b} _{\Q, \textrm{dév}}
(\overset{^\mathrm{g}}{} \smash{\widehat{\D}} _{\PP} ^{(\bullet)} (T ))$ 
la sous-catégorie pleine de 
$\underrightarrow{LD}  ^\mathrm{b} _{\Q, \mathrm{qc}}
(\overset{^\mathrm{g}}{} \smash{\widehat{\D}} _{\PP} ^{(\bullet)} (T ))$
des complexes dévissables en isocristaux surconvergents.
\end{enumerate}

\end{defi}

\begin{rema}
\label{rema-dev=dev-plong}
Avec les notations de \ref{coro-caradev}, supposons $P$ propre. 

1) Lorsque $Y$ est une sous-variété $d$-plongeable de $P$, 
il est immédiat que la définition de \ref{coro-caradev} est la même que celle donnée dans 
\cite[3.2.5]{caro-2006-surcoh-surcv}.

2) Pour le cas général, la définition de \cite[3.2.14]{caro-2006-surcoh-surcv} de la dévissabilité est la suivante:
le complexe $\E ^{(\bullet)}$ est {\og dévissable sur $Y$ en isocristaux surconvergents\fg} s'il existe
un recouvrement fini ouvert $(Y _l) _l$ de $Y$ par des $d$-sous-variétés plongeables dans $P$ tel 
que, pour tout $l$,
$\E ^{(\bullet)}$ soit dévissable sur $Y _l$ en isocristaux surconvergents (au sens de \ref{coro-caradev}).
Cette définition paraît à première vue différente. En fait, grâce à la proposition \ref{Y'cupYdev} ci-dessous, 
ces deux définitions coïncident. Ainsi, la définition \ref{coro-caradev} 
étend celle de \cite[3.2.14]{caro-2006-surcoh-surcv} sans hypothèse de propreté sur $P$. 
Pour valider cette proposition \ref{Y'cupYdev}, nous aurons d'abord besoin d'établir les lemmes ou proposition ci-dessous qui sont souvent 
des analogues d'énoncés de \cite[3.2]{caro-2006-surcoh-surcv} qui deviendront équivalents {\it a posteriori} dans le cas où $P$ est propre. 
\end{rema}

\begin{lemm}
 \label{lemm-3.2.7-coh-cv}
Soient $Y,Y'$ deux sous-variétés lisses $d$-plongeables dans $P$
tels que $Y '\subset Y$.
Soit $\E ^{(\bullet)} \in \underrightarrow{LD}  ^\mathrm{b} _{\Q, \mathrm{qc}}
(\overset{^\mathrm{g}}{} \smash{\widehat{\D}} _{\PP} ^{(\bullet)})$.
Si $\R\underline{\Gamma} ^\dag _{Y }  (\E ^{(\bullet)}) \in 
\smash{\underrightarrow{LD}} ^{\mathrm{b}} _{\Q ,\mathrm{isoc},Y }
(\smash{\widehat{\D}} _{\PP} ^{(\bullet)})$
alors 
$\R\underline{\Gamma} ^\dag _{Y'}  (\E ^{(\bullet)}) \in 
\smash{\underrightarrow{LD}} ^{\mathrm{b}} _{\Q ,\mathrm{isoc},Y '}
(\smash{\widehat{\D}} _{\PP} ^{(\bullet)})$.
\end{lemm}

\begin{proof}
Choisissons deux cadres de la forme $(\PP, T,X,Y)$ et $(\PP, T',X',Y')$.
D'après 
\ref{stabIsoc*inv-i-bullet},
on dispose du foncteur exact
$\R\underline{\Gamma} ^\dag _{Y'} =\R\underline{\Gamma} ^\dag _{X'} \circ (\hdag T')
\colon 
\mathrm{Isoc} ^{(\bullet)} ( \PP, T ,X /K)
\to 
\mathrm{Isoc} ^{(\bullet)} ( \PP, T ',X '/K)$.
D'où le résultat.
\end{proof}

\begin{lemm}
\label{trans-strat}
Soient $Y$ une sous-variété de $P$,
$Y =\sqcup _{i=1,\dots , r} Y _i$ une stratification de $Y$ dans $P$.
Pour tout $i=1,\dots , r$, soit 
$Y _i =\sqcup _{j=1,\dots ,  j _i} Y _{i,j}$
une stratification lisse 
(resp. une $d$-stratification lisse dans $P$).
Alors
$Y= (\sqcup _{j=1,\dots ,  j _1} Y _{1,j}) \sqcup \dots \sqcup 
(\sqcup _{j=1,\dots ,  j _r} Y _{r,j})$
est une stratification lisse 
(resp. une $d$-stratification lisse dans $P$).
On dira que cette stratification est une {\og sous-stratification\fg} 
(resp. une {\og sous-$d$-stratification lisse dans $P$\fg}) de 
$Y =\sqcup _{i=1,\dots , r} Y _i$.
\end{lemm}

\begin{proof}
Soient $1\leq i\leq r$ et  $1\leq j\leq j _i$. 
Comme $Y _i$ est un ouvert 
de $\sqcup _{i'=i,\dots ,r} Y _{i'}$, alors
$\sqcup _{j'=j,\dots ,  j _i}Y _{i,j'}$
est un ouvert de 
$Z _{i,j}:=(\sqcup _{j'=j,\dots ,  j _i}Y _{i,j'}) 
\sqcup
(\sqcup _{j=1,\dots ,  j _{i+1}} Y _{1,j}) \sqcup \dots \sqcup 
(\sqcup _{j=1,\dots ,  j _r} Y _{r,j})$.
Comme   
$Y _{i,j} $ est un ouvert de 
$\sqcup _{j'=j,\dots ,  j _i}Y _{i,j'}$, 
alors$Y _{i,j} $ est un ouvert de $Z _{i,j}$.
\end{proof}

\begin{prop}
 \label{3.2.7-coh-cv}
Soient $\E ^{(\bullet)} \in \underrightarrow{LD}  ^\mathrm{b} _{\Q, \mathrm{qc}}
(\overset{^\mathrm{g}}{} \smash{\widehat{\D}} _{\PP} ^{(\bullet)})$,
$Y$ une sous-variété de $P$
et
$Y =\sqcup _{i=1,\dots , r} Y _i$ une stratification de $Y$ dans $P$.
Si, pour tout $i=1,\dots r$,
le complexe $\E ^{(\bullet)}$
se dévisse en isocristaux surconvergents sur $Y _i$, 
alors le complexe $\E ^{(\bullet)}$ se dévisse en isocristaux surconvergents 
au-dessus d'une sous-$d$-stratification lisse dans $P$ de
$Y =\sqcup _{i=1,\dots , r} Y _i$.
\end{prop}

\begin{proof}
Cela découle aussitôt des deux lemmes \ref{lemm-3.2.7-coh-cv} et \ref{trans-strat}.
\end{proof}

\begin{lemm}
\label{exist-d-strat-lisse}
Soit $Y$ une sous-variété de $P$.
Il existe une $d$-stratification lisse de $Y$ dans $P$.
\end{lemm}

\begin{proof}
Comme $P$ est lisse, on se ramène facilement au cas où
$P$ est intègre. 
Notons $X$ l'adhérence de $Y$ dans $P$.
Le cas où $Y=X$ est évident. 
Sinon, il existe alors des diviseurs $T _1, \dots, T _r$ de $P$ tels que 
$X \setminus Y= \cap _{j=1,\dots, r} T _j$.
Si $y$ est un point générique d'une composante irréductible de dimension 
égale à $\dim Y$, alors il existe un $j$ tel que $y \not \in T _j$.
On pose alors $Y _1:= X\setminus T _j \subset Y$.
En procédant par récurrence lexicographique sur la dimension de $Y$ et 
sur le nombre de composantes irréductibles de degré maximal, on en déduit 
une $d$-stratification lisse de $Y\setminus Y _1$ dans $P$. D'où le résultat.
\end{proof}

\begin{lemm}
\label{Y'subsetYdev}
Soient $Y$ une sous-variété de $P$,
$Y '$ une sous-variété de $Y$.
Soit $\E ^{(\bullet)} \in \underrightarrow{LD}  ^\mathrm{b} _{\Q, \mathrm{qc}}
(\overset{^\mathrm{g}}{} \smash{\widehat{\D}} _{\PP} ^{(\bullet)})$.
Si $\E ^{(\bullet)}$ se dévisse sur $Y$ en isocristaux surconvergents alors il l'est sur $Y'$.
En particulier, 
la réciproque de la proposition \ref{3.2.7-coh-cv} est valable.
\end{lemm}

\begin{proof}
Soit $Y =\sqcup _{i=1,\dots , r} Y _i$
une $d$-stratification lisse de $Y$ dans $P$ au-dessus de laquelle 
$\E ^{(\bullet)}$ se dévisse en isocristaux surconvergents.
On obtient la stratification
$Y '=\sqcup _{i=1,\dots , r} (Y _i \cap Y')$.
Par \ref{3.2.7-coh-cv}, on se ramène alors au cas où $Y$ est lisse, $d$-plongeable dans $P$
et où
$\R\underline{\Gamma} ^\dag _{Y}  (\E ^{(\bullet)}) \in 
\smash{\underrightarrow{LD}} ^{\mathrm{b}} _{\Q ,\mathrm{isoc},Y}
(\smash{\widehat{\D}} _{\PP} ^{(\bullet)})$.
Dans ce cas, d'après le lemme \ref{lemm-3.2.7-coh-cv},
$\E ^{(\bullet)}$ se dévisse en isocristaux surconvergents sur n'importe quelle $d$-stratification lisse de $Y'$ dans $P$. 
Or, par \ref{exist-d-strat-lisse}, il existe une $d$-stratification lisse de $Y'$ dans $P$.
\end{proof}

\begin{prop}
\label{dev-triangulé}
Soit $Y$ une sous-variété de $P$. 

\begin{enumerate}
\item Soient $\E ^{(\bullet)}, \FF ^{(\bullet)}  \in \underrightarrow{LD}  ^\mathrm{b} _{\Q, \mathrm{qc}}
(\overset{^\mathrm{g}}{} \smash{\widehat{\D}} _{\PP} ^{(\bullet)})$
dévissables en isocristaux surconvergents sur $Y$. 
Il existe alors une $d$-stratification lisse de $Y$ dans $P$ telle que 
$\E ^{(\bullet)}$ et 
$ \FF ^{(\bullet)} $
se dévissent simultanément en isocristaux surconvergents au-dessus de celle-ci.  

\item La sous-catégorie pleine de
$\underrightarrow{LD}  ^\mathrm{b} _{\Q, \mathrm{qc}}
(\overset{^\mathrm{g}}{} \smash{\widehat{\D}} _{\PP} ^{(\bullet)})$
des complexes dévissables en isocristaux surconvergents sur $Y$
est triangulée.
\end{enumerate}
\end{prop}

\begin{proof}
La preuve est analogue à \cite[3.2.10]{caro-2006-surcoh-surcv}:
soit $Y =\sqcup _{i=1,\dots , r} Y _i$
une $d$-stratification lisse de $Y$ dans $P$ au-dessus de laquelle 
$\E ^{(\bullet)}$ 
se dévisse en isocristaux surconvergents.
D'après \ref{Y'subsetYdev}, 
$\FF ^{(\bullet)}$ se dévisse en isocristaux surconvergents au-dessus de chacun de $Y _i$.
Par \ref{3.2.7-coh-cv},
le complexe $\FF ^{(\bullet)}$ se dévisse en isocristaux surconvergents 
au-dessus d'une sous-$d$-stratification lisse dans $P$ de
$Y =\sqcup _{i=1,\dots , r} Y _i$.
Par \ref{lemm-3.2.7-coh-cv}, c'est aussi le cas de $\E ^{(\bullet)}$ au-dessus de cette dernière.
Déduisons-en la seconde assertion de la proposition. 
Par dévissage en isocristaux surconvergents de $\E ^{(\bullet)}, \FF ^{(\bullet)}  $
au-dessus de la même $d$-stratification lisse dans $P$, 
on se ramène alors au cas immédiat où $Y$ est lisse et $d$-plongeable dans $P$
et où $\E ^{(\bullet)}$ et $\FF ^{(\bullet)}$ sont des objets de 
$\smash{\underrightarrow{LD}} ^{\mathrm{b}} _{\Q ,\mathrm{isoc},Y}
(\smash{\widehat{\D}} _{\PP} ^{(\bullet)})$, ce qui est immédiat.

\end{proof}

\begin{prop}
\label{Y'cupYdev}
Soient $Y$, $Y'$ deux sous-variétés de $P$.
Soit $\E ^{(\bullet)} \in \underrightarrow{LD}  ^\mathrm{b} _{\Q, \mathrm{qc}}
(\overset{^\mathrm{g}}{} \smash{\widehat{\D}} _{\PP} ^{(\bullet)})$.
Le complexe $\E ^{(\bullet)}$ se dévisse en isocristaux surconvergents  sur $Y\cup Y'$ si et seulement si 
le complexe $\E ^{(\bullet)}$ se dévisse en isocristaux surconvergents  sur $Y$ et sur $Y'$.

\end{prop}

\begin{proof}
La nécessité découle de \ref{Y'subsetYdev}.
Supposons à présent que le complexe $\E ^{(\bullet)}$ se dévisse en isocristaux surconvergents  sur $Y$ et sur $Y'$.
On dispose du triangle distingué de Mayer-Vietoris:
\begin{equation}
\label{Y'cupYdev-trian-dist}
\R\underline{\Gamma} ^\dag _{Y \cap Y '}  (\E ^{(\bullet)})
\to
\R\underline{\Gamma} ^\dag _{Y  }  (\E ^{(\bullet)})
\oplus
\R\underline{\Gamma} ^\dag _{Y ' }  (\E ^{(\bullet)})
\to
\R\underline{\Gamma} ^\dag _{Y  \cup Y '}  (\E ^{(\bullet)})
\to
\R\underline{\Gamma} ^\dag _{Y  \cap Y '}  (\E ^{(\bullet)})[1].
\end{equation}
Posons $\FF ^{(\bullet)}:= \R\underline{\Gamma} ^\dag _{Y  }  (\E ^{(\bullet)})$.
Comme $\R\underline{\Gamma} ^\dag _{Y  }  (\FF ^{(\bullet)})= \FF ^{(\bullet)}$
et
comme $\R\underline{\Gamma} ^\dag _{Y '\setminus Y }  (\FF ^{(\bullet)})=0$, 
par \ref{3.2.7-coh-cv}, on en déduit que 
$\FF ^{(\bullet)}$ se dévisse en isocristaux surconvergents au-dessus d'une sous-$d$-stratification lisse dans $P$
de la stratification
$Y \cup Y'= Y \sqcup (Y '\setminus Y)$.
De même, on vérifie que 
$\R\underline{\Gamma} ^\dag _{Y \cap Y '}  (\E ^{(\bullet)})$
et
$\R\underline{\Gamma} ^\dag _{Y '}  (\E ^{(\bullet)})$ se dévissent en isocristaux surconvergents sur $Y \cup Y'$. 
Par \ref{dev-triangulé}, on conclut alors via le triangle distingué \ref{Y'cupYdev-trian-dist}.

\end{proof}

\begin{vide}
Soit $(\PP, T,X,Y)$ un cadre.
 Les propriétés de 
\cite[3.2]{caro-2006-surcoh-surcv} restent valable en passant du cas où $\PP$ est propre au cas général (i.e. $\PP$ est séparé et lisse). 
Par exemple:
\label{3.2.23-coh-cv}
La catégorie $\underrightarrow{LD}  ^\mathrm{b} _{\Q, \textrm{dév}}
(\overset{^\mathrm{g}}{} \smash{\widehat{\D}} _{\PP} ^{(\bullet)})$
est la plus petite sous-catégorie pleine triangulée de
$\underrightarrow{LD}  ^\mathrm{b} _{\Q, \mathrm{qc}}
(\overset{^\mathrm{g}}{} \smash{\widehat{\D}} _{\PP} ^{(\bullet)})$
contenant $\smash{\underrightarrow{LD}} ^{\mathrm{b}} _{\Q ,\mathrm{isoc},Y'}
(\smash{\widehat{\D}} _{\PP} ^{(\bullet)})$, pour toute
sous-variété $Y'$ lisse $d$-plongeable dans $P$.
\end{vide}

\begin{vide}
\label{4.1.16}
Soit $(\PP, T,X,Y)$ un cadre.
On note 
$\underrightarrow{LD}  ^\mathrm{b} _{\Q, \textrm{dév}} (\PP, T, X/K)$
la sous-catégorie strictement pleine de 
$\underrightarrow{LD}  ^\mathrm{b} _{\Q, \textrm{dév}}
(\overset{^\mathrm{g}}{} \smash{\widehat{\D}} _{\PP} ^{(\bullet)} (T ))$
des complexes $\E ^{(\bullet)}$
tels que le morphisme canonique 
$\R\underline{\Gamma} ^\dag _{X}  (\E ^{(\bullet)})
\to 
\E ^{(\bullet)}$
soit un isomorphisme.

$\bullet$ Cette catégorie $\underrightarrow{LD}  ^\mathrm{b} _{\Q, \textrm{dév}} (\PP, T, X/K)$ 
est égale à la sous-catégorie strictement pleine de 
$\underrightarrow{LD}  ^\mathrm{b} _{\Q, \textrm{qc}}
(\overset{^\mathrm{g}}{} \smash{\widehat{\D}} _{\PP} ^{(\bullet)} )$
des complexes $\E ^{(\bullet)}$ dévissables en isocristaux surconvergents au-dessus de $Y$ et
tels que l'on dispose de l'isomorphisme
$\R\underline{\Gamma} ^\dag _{Y}  (\E ^{(\bullet)})
\riso 
\E ^{(\bullet)}$ (en effet, il suffit de considérer la stratification 
$P= (P \setminus X) \sqcup Y \sqcup (X \cap T)$ et d'invoquer 
\ref{3.2.7-coh-cv}).

$\bullet$ Le foncteur  \ref{eq-coh-lim} induit l'équivalence de catégories
\begin{gather}
\label{equi-devcoh-devqc}
\underrightarrow{\lim} \colon
\underrightarrow{LD}  ^\mathrm{b} _{\Q, \textrm{dév}}
 (\PP, T, X/K)
\cap
\underrightarrow{LD}  ^\mathrm{b} _{\Q, \mathrm{coh}}
(\overset{^\mathrm{g}}{} \smash{\widehat{\D}} _{\PP} ^{(\bullet)} (T ))
\cong
D ^\mathrm{b} _{\textrm{dév}}
 (\PP, T, X/K),
\end{gather}
où la catégorie à droite a été définie dans \cite[6.2.2]{caro-pleine-fidelite}
(voir \ref{rappel-dev-coh}).

$\bullet$ Comme les isocristaux surconvergents munis d'une structure de Frobenius sont surholonomes
(e.g. voir \cite{caro-Tsuzuki}), avec de plus \ref{pl-fid-surcoh}, 
on obtient l'égalité
$F\text{-}\underrightarrow{LD}  ^\mathrm{b} _{\Q, \textrm{dév}}
(\overset{^\mathrm{g}}{} \smash{\widehat{\D}} _{\PP} ^{(\bullet)} (T ))
=
F\text{-}\underrightarrow{LD}  ^\mathrm{b} _{\Q, \textrm{surcoh}}
(\overset{^\mathrm{g}}{} \smash{\widehat{\D}} _{\PP} ^{(\bullet)} (T ))$
(sans structure de Frobenius, on ne bénéficie que de l'inclusion
$\underrightarrow{LD}  ^\mathrm{b} _{\Q, \textrm{surcoh}}
(\overset{^\mathrm{g}}{} \smash{\widehat{\D}} _{\PP} ^{(\bullet)} (T ))
\subset
\underrightarrow{LD}  ^\mathrm{b} _{\Q, \textrm{dév}}
(\overset{^\mathrm{g}}{} \smash{\widehat{\D}} _{\PP} ^{(\bullet)} (T ))$).
De plus, on dispose alors de l'équivalence de catégories:
\begin{gather}
\label{Frob-equi-devcoh-devqc}
\underrightarrow{\lim} \colon
F\text{-}\underrightarrow{LD}  ^\mathrm{b} _{\Q, \textrm{dév}}
 (\PP, T, X/K)
\cong
F\text{-}D ^\mathrm{b} _{\textrm{dév}}
 (\PP, T, X/K).
\end{gather}

\end{vide}

\begin{rema}
On peut définir 
$D ^\mathrm{b} _{\textrm{dév}}
(\smash{\D} ^\dag _{\PP} (\hdag T) _\Q)$
comme étant la sous-catégorie strictement pleine
de 
$D ^\mathrm{b} _{\textrm{coh}}
(\smash{\D} ^\dag _{\PP} (\hdag T) _\Q)$
qui induise l'équivalence de catégories 
\ref{equi-devcoh-devqc} donnée ci-dessus. 

\end{rema}

\subsection{Sur la stabilité par produit tensoriel de la dévissabilité en isocristaux, surholonomie et holonomie}

\begin{lemm}
\label{otimes-iscoY}
Soient $(\PP, T,X,Y)$ et $(\PP', T',X',Y')$ deux cadres lisses en dehors de leur diviseur 
(voir la définition \ref{defi-cadre}) et $\PP '':= \PP \times \PP'$.
\begin{enumerate}
\item
Le bifoncteur 
$\smash{\overset{\L}{\boxtimes}}   ^{\dag} _{\O _\S} $
se factorise sous la forme
\begin{equation}
\label{otimes-iscoY1}
\smash{\overset{\L}{\boxtimes}}   ^{\dag} _{\O _\S} 
\colon 
\smash{\underrightarrow{LD}} ^{\mathrm{b}} _{\Q ,\mathrm{isoc},Y}
(\smash{\widehat{\D}} _{\PP} ^{(\bullet)})
\times 
\smash{\underrightarrow{LD}} ^{\mathrm{b}} _{\Q ,\mathrm{isoc},Y'}
(\smash{\widehat{\D}} _{\PP'} ^{(\bullet)})
\to 
\smash{\underrightarrow{LD}} ^{\mathrm{b}} _{\Q ,\mathrm{isoc},Y\times Y'}
(\smash{\widehat{\D}} _{\PP''} ^{(\bullet)}).
\end{equation}

\item Si $\PP= \PP'$ et si $Y \cap Y'$ est lisse, 
le bifoncteur 
$\smash{\overset{\L}{\otimes}}   ^{\dag}
_{\O  _{\PP,\Q}}
$
se factorise sous la forme
\begin{equation}
\label{otimes-iscoY2}
\smash{\overset{\L}{\otimes}}   ^{\dag}
_{\O  _{\PP,\Q}}
\colon 
\smash{\underrightarrow{LD}} ^{\mathrm{b}} _{\Q ,\mathrm{isoc},Y}
(\smash{\widehat{\D}} _{\PP} ^{(\bullet)})
\times 
\smash{\underrightarrow{LD}} ^{\mathrm{b}} _{\Q ,\mathrm{isoc},Y'}
(\smash{\widehat{\D}} _{\PP} ^{(\bullet)})
\to 
\smash{\underrightarrow{LD}} ^{\mathrm{b}} _{\Q ,\mathrm{isoc},Y\cap Y'}
(\smash{\widehat{\D}} _{\PP} ^{(\bullet)}).
\end{equation}

\end{enumerate}

\end{lemm}

\begin{proof}
Comme la catégorie 
$\smash{\underrightarrow{LD}} ^{\mathrm{b}} _{\Q ,\mathrm{isoc},Y\times Y'}
(\smash{\widehat{\D}} _{\PP''} ^{(\bullet)})$ (resp. $\smash{\underrightarrow{LD}} ^{\mathrm{b}} _{\Q ,\mathrm{isoc},Y\cap Y'}
(\smash{\widehat{\D}} _{\PP'} ^{(\bullet)})$)
est une sous-catégorie triangulée de 
$\smash{\underrightarrow{LD}} ^{\mathrm{b}} _{\Q ,\mathrm{qc}}
(\smash{\widehat{\D}} _{\PP''} ^{(\bullet)})$
(resp. $\smash{\underrightarrow{LD}} ^{\mathrm{b}} _{\Q ,\mathrm{qc}}
(\smash{\widehat{\D}} _{\PP'} ^{(\bullet)})$),
quitte à utiliser des triangles distingués de troncation et à procéder par récurrence
sur le nombre d'espace de cohomologie non nul, on se ramène au cas où les complexes
de 
$\smash{\underrightarrow{LD}} ^{\mathrm{b}} _{\Q ,\mathrm{isoc},Y}
(\smash{\widehat{\D}} _{\PP} ^{(\bullet)})$ 
(resp. $\smash{\underrightarrow{LD}} ^{\mathrm{b}} _{\Q ,\mathrm{isoc},Y'}
(\smash{\widehat{\D}} _{\PP'} ^{(\bullet)})$)
sont 
des objets de 
$ \mathrm{Isoc} ^{(\bullet)}( \PP, T, X/K)$
(resp. $ \mathrm{Isoc} ^{(\bullet)}( \PP', T', X'/K)$),
i.e. 
à la situation déjà traitée en respectivement \ref{boxtimesDT} et
\ref{stbintnlisse-bullet}.
\end{proof}

\begin{theo}
\label{prod-tens-ext-formel-pre}
Soient $(\PP, D,X,Y)$ et $(\PP', D',X',Y')$ deux cadres.
On note $\PP'':=\PP \times \PP '$, $X '' := X \times X'$, 
$p \colon \PP'' \to \PP$
et
$p' \colon \PP'' \to \PP'$ les projections canoniques et
$D''= p ^{-1} (D) \cup p ^{\prime -1} (D')$. 
\begin{enumerate}
\item 
On dispose des factorisations
\begin{align}
\label{stabdév-botimes-pre}
\smash{\overset{\L}{\boxtimes}}   ^{\dag} _{\O _\S} 
  \colon&
\underrightarrow{LD}  ^\mathrm{b} _{\Q, \text{\rm dév}} (\PP, D, X/K) \times
\underrightarrow{LD}  ^\mathrm{b} _{\Q, \text{\rm dév}}(\PP', D', X'/K)
\rightarrow
\underrightarrow{LD}  ^\mathrm{b} _{\Q, \text{\rm dév}} (\PP'', D'', X''/K),
\\
\label{stabdév-otimes-pre}
\smash{\overset{\L}{\otimes}}^{\dag} _{\O _{\PP } ( \hdag D ) _{\Q}}
  \colon&
\underrightarrow{LD}  ^\mathrm{b} _{\Q, \text{\rm dév}}   (\PP, D, X/K) \times
\underrightarrow{LD}  ^\mathrm{b} _{\Q, \text{\rm dév}}  (\PP, D, X/K)
\rightarrow
\underrightarrow{LD}  ^\mathrm{b} _{\Q, \text{\rm dév}}   (\PP, D, X/K).
\end{align}
\item La catégorie 
$F\text{-}\underrightarrow{LD}  ^\mathrm{b} _{\Q, \textrm{surcoh}}
(\overset{^\mathrm{g}}{} \smash{\widehat{\D}} _{\PP} ^{(\bullet)} (D))$
est  stable par produit tensoriel.
\end{enumerate}
\end{theo}

\begin{proof}
Vérifions d'abord \ref{stabdév-botimes-pre}.
 Soient 
 $\E ^{(\bullet)}\in\underrightarrow{LD}  ^\mathrm{b} _{\Q, \text{\rm dév}} (\PP, D, X/K) $
et
 $\E ^{\prime(\bullet)}\in\underrightarrow{LD}  ^\mathrm{b} _{\Q, \text{\rm dév}}(\PP', D', X'/K)$.
Soit 
$Y =\sqcup _{i=1,\dots , r} Y _i$
une $d$-stratification lisse de $Y$ dans $P$ au-dessus de laquelle 
$\E ^{(\bullet)}$ se dévisse en isocristaux surconvergents.
Soit 
$Y '=\sqcup _{j=1,\dots , s} Y '_j$
une $d$-stratification lisse de $Y'$ dans $P'$ au-dessus de laquelle 
$\E ^{\prime (\bullet)}$ se dévisse en isocristaux surconvergents.
Comme $\underrightarrow{LD}  ^\mathrm{b} _{\Q, \text{\rm dév}} (\PP'', D'', X''/K)$
est une sous-catégorie triangulée de 
$\underrightarrow{LD}  ^\mathrm{b} _{\Q, \textrm{qc}}
(\overset{^\mathrm{g}}{} \smash{\widehat{\D}} _{\PP''} ^{(\bullet)})$,
on se ramène par dévissage (voir la remarque \ref{rema-dev-tridist}) au cas où 
$ \E ^{(\bullet)} \in \smash{\underrightarrow{LD}} ^{\mathrm{b}} _{\Q ,\mathrm{isoc},Y_i}
(\smash{\widehat{\D}} _{\PP} ^{(\bullet)}) $
et
$ \E ^{\prime (\bullet)} \in \smash{\underrightarrow{LD}} ^{\mathrm{b}} _{\Q ,\mathrm{isoc},Y' _j}
(\smash{\widehat{\D}} _{\PP'} ^{(\bullet)}) $.
Grâce à 
\ref{otimes-iscoY1},
on en déduit que 
$\E ^{(\bullet)}\smash{\overset{\L}{\boxtimes}}   ^{\dag} _{\O _\S}\E ^{\prime (\bullet)}
\in
\smash{\underrightarrow{LD}} ^{\mathrm{b}} _{\Q ,\mathrm{isoc},Y _i\times Y' _j}
(\smash{\widehat{\D}} _{\PP''} ^{(\bullet)})
\subset
\underrightarrow{LD}  ^\mathrm{b} _{\Q, \text{\rm dév}} (\PP'', D'', X''/K)$ (l'inclusion résulte du premier point de 
\ref{4.1.16}).
D'où le résultat.
Traitons à présent \ref{stabdév-otimes-pre}.
Par \ref{dev-triangulé}, 
il existe $Y =\sqcup _{i=1,\dots , r} Y _i$
une même $d$-stratification lisse de $Y$ dans $P$ au-dessus de laquelle 
$\E ^{(\bullet)}$ et $\E ^{\prime (\bullet)}$ se dévissent en isocristaux surconvergents.
On procède alors de manière analogue à la vérification de \ref{stabdév-botimes-pre} mais
en utilisant \ref{otimes-iscoY2} à la place de \ref{otimes-iscoY1}.
Enfin, la dernière assertion découle de l'égalité
$F\text{-}\underrightarrow{LD}  ^\mathrm{b} _{\Q, \textrm{dév}}
(\overset{^\mathrm{g}}{} \smash{\widehat{\D}} _{\PP} ^{(\bullet)} (D ))
=
F\text{-}\underrightarrow{LD}  ^\mathrm{b} _{\Q, \textrm{surcoh}}
(\overset{^\mathrm{g}}{} \smash{\widehat{\D}} _{\PP} ^{(\bullet)} (D ))$ (voir le dernier point de 
\ref{4.1.16}).
\end{proof}

\begin{theo}
\label{prod-tens-ext-formel}

Soient $(\PP, T,X,Y)$ et $(\PP', T',X',Y')$ deux cadres.
On note $\PP'':=\PP \times \PP '$, $X '' := X \times X'$, 
$p \colon \PP'' \to \PP$
et
$p' \colon \PP'' \to \PP'$ les projections canoniques et
$T''= p ^{-1} (T) \cup p ^{\prime -1} (T')$. 
Le foncteur \ref{prod-tens-extlimcoh} 
induit les factorisations
\begin{align}
\label{stabdév-botimes}
\smash{\overset{\L}{\boxtimes}}   ^{\dag} _{\O _\S, T, T '} 
  \colon&
  D ^\mathrm{b} _{\textrm{\rm dév}}  (\PP, T, X/K) \times
D ^\mathrm{b} _{\mathrm{\text{\rm dév}}}  (\PP', T', X'/K)
\rightarrow
D ^\mathrm{b} _{\mathrm{\text{\rm dév}}}  (\PP'', T'', X''/K),
\\
\label{stabdév-botimesF}
\smash{\overset{\L}{\boxtimes}}   ^{\dag} _{\O _\S, T, T '}  
\colon&
F\text{-}D ^\mathrm{b} _\mathrm{surhol}  (\PP, T, X/K) \times
F\text{-}D ^\mathrm{b} _\mathrm{surhol}  (\PP', T', X'/K)
\rightarrow
F\text{-}D ^\mathrm{b} _\mathrm{surhol}  (\PP'', T'', X''/K).
\end{align}

\end{theo}

\begin{proof}
Vérifions dans un premier temps \ref{stabdév-botimes}.
 Soient 
$\E ^{(\bullet)} \in \underrightarrow{LD}  ^\mathrm{b} _{\Q, \textrm{dév}}
 (\PP, T, X/K)
\cap
\underrightarrow{LD}  ^\mathrm{b} _{\Q, \mathrm{coh}}
(\overset{^\mathrm{g}}{} \smash{\widehat{\D}} _{\PP} ^{(\bullet)} (T ))
$
et
$\E ^{\prime (\bullet)} \in \underrightarrow{LD}  ^\mathrm{b} _{\Q, \textrm{dév}}
 (\PP', T', X'/K)
\cap
\underrightarrow{LD}  ^\mathrm{b} _{\Q, \mathrm{coh}}
(\overset{^\mathrm{g}}{} \smash{\widehat{\D}} _{\PP'} ^{(\bullet)} (T' ))
$,
$\E := \underrightarrow{\lim}  ~\E ^{(\bullet)}$
 et
$ \E' := \underrightarrow{\lim}  ~\E ^{\prime (\bullet)}$.
 On dispose par définition de l'égalité:
$\E\smash{\overset{\L}{\boxtimes}}   ^{\dag} _{\O _\S, D, D '}\E ' = 
\underrightarrow{\lim} \,
\E^{(\bullet)}\smash{\overset{\L}{\boxtimes}}   ^{\dag} _{\O _\S}\E ^{\prime(\bullet)} $.
 Par stabilité de la cohérence et de la dévissabilité (voir \ref{stabdév-botimes-pre} pour cette dernière) 
par produit tensoriel externe et via l'équivalence de catégories \ref{equi-devcoh-devqc},
on obtient alors la  factorisation de \ref{stabdév-botimes}. 
Grâce à \ref{=Fsurholdev}, 
la factorisation \ref{stabdév-botimesF} résulte de celle de 
\ref{stabdév-botimes}.
\end{proof}

\begin{theo}
\label{stabdév}
Soit $(\PP, D,X,Y)$ un cadre. Le foncteur 
\ref{def-otimes-coh1}
induit les factorisations
\begin{align}
\label{stabdév-otimesF}
\smash{\overset{\L}{\otimes}}^{\dag} _{\O _{\PP } ( \hdag D ) _{\Q}}  
\colon&
F\text{-}D ^\mathrm{b} _\mathrm{surhol}  (\PP, D, X/K) \times
F\text{-}D ^\mathrm{b} _\mathrm{surhol}  (\PP, D, X/K)
\rightarrow
F\text{-}D ^\mathrm{b} _\mathrm{surhol}  (\PP, D, X/K).
\end{align}

\end{theo}

\begin{proof}
Grâce à \ref{=Fsurholdev} et \ref{Frob-equi-devcoh-devqc}, 
cela résulte aussitôt de \ref{stabdév-otimes-pre}.
\end{proof}

\begin{rema}
\label{rema-stab-tens-qcdev++}
Les théorèmes \ref{prod-tens-ext-formel} et \ref{stabdév} sont des  conséquences du théorème 
\ref{prod-tens-ext-formel-pre}, ce dernier étant plus fort.
La raison principale est que si l'image par le foncteur 
$\underrightarrow{\lim}$ d'un complexe  quasi-cohérent $\E ^{(\bullet)}$ est (sur)cohérent, alors il n'est pas vrai
que $\E ^{(\bullet)}$ soit lui-même (sur)cohérent. 
\end{rema}

\begin{theo}
\label{theo-stab-prod-surhol}
Soient $(Y,X)$ et $(Y',X')$ deux couples de $k$-variétés proprement $d$-plongeables.
\begin{enumerate}
\item 
Soit $(\PP, T,X,Y)$ un cadre tel qu'il existe un $\V$-schéma formel $\mathcal{Q} $ propre et lisse et une immersion ouverte de la forme 
$\PP \hookrightarrow \mathcal{Q}$. 
Le bifoncteur produit tensoriel 
$ \smash{\overset{\L}{\otimes}} ^\dag _{\O _{\PP,\Q} }   [-d _{P}]$
(voir \ref{stabdév-otimesF} au décalage près)
se factorise en le bifoncteur produit tensoriel que l'on notera:
\begin{equation}
  \smash{\overset{\L}{\otimes}}^{\dag}  _{\O _{(Y,X)/K}} 
  \ : \
F \text{-}D ^\mathrm{b} _\mathrm{surhol}  (\D ^\dag _{(Y,X)/K}) \times
F \text{-}D ^\mathrm{b} _\mathrm{surhol}  (\D ^\dag _{(Y,X)/K})
\rightarrow
F \text{-}D ^\mathrm{b} _\mathrm{surhol} (\D ^\dag _{(Y,X)/K}).
\end{equation}
\item De même, le produit tensoriel externe
$\smash{\overset{\L}{\boxtimes}}   ^{\dag} _{\O _\S}\, $  induit le foncteur noté
\begin{equation}
  \smash{\overset{\L}{\boxtimes}}   ^{\dag} _{\O _\S}\, 
  \ : \
F \text{-}D ^\mathrm{b} _\mathrm{surhol}  (\D ^\dag _{(Y',X')/K}) \times
F \text{-}D ^\mathrm{b} _\mathrm{surhol}  (\D ^\dag _{(Y,X)/K})
\rightarrow
F \text{-}D ^\mathrm{b} _\mathrm{surhol} (\D ^\dag _{(Y'\times Y,X'\times X)/K}).
\end{equation} 
\end{enumerate}

\end{theo}

\begin{proof}
Soit $(\PP, T,X,Y)$ un cadre tel qu'il existe un $\V$-schéma formel $\mathcal{Q} $ propre et lisse et une immersion ouverte de la forme 
$\PP \hookrightarrow \mathcal{Q}$. 
Pour que le bifoncteur $  \smash{\overset{\L}{\otimes}}^{\dag}  _{\O _{(Y,X)/K}} $ soit bien défini,
il s'agit de vérifier que le bifoncteur 
$  \smash{\overset{\L}{\otimes}} ^\dag _{\O _{\PP,\Q} } [-d _{P}]$ ne dépend pas, à isomorphisme canonique près, 
du choix.
Faisons donc un second choix :
soit $(\widetilde{\PP}, \widetilde{T},X,Y)$ un cadre tel qu'il existe un $\V$-schéma formel $\widetilde{\mathcal{Q}} $ propre et lisse et une immersion ouverte de la forme 
$\widetilde{\PP} \hookrightarrow \widetilde{\mathcal{Q}} $. 
  Quitte comme d'habitude à remplacer $\widetilde{\mathcal{Q}}$ par $\mathcal{Q} \times \widetilde{\mathcal{Q}}$
  et $\widetilde{\mathcal{P}}$ par $\mathcal{P} \times \widetilde{\mathcal{P}}$, on peut supposer qu'il existe
  un morphisme $\widetilde{f}\colon \widetilde{\mathcal{Q}}\rightarrow \mathcal{Q}$ induisant 
  $f \colon  \widetilde{\PP} \to \PP$.
Avec \ref{fg!prodtens}, on dispose alors du diagramme commutatif à isomorphisme canonique près :
\begin{equation}
  \notag
  \xymatrix @C=2cm {
{ F\text{-}D ^\mathrm{b} _\mathrm{surhol}  (\PP, T, X/K) \times F\text{-}D ^\mathrm{b} _\mathrm{surhol}  (\PP, T, X/K)}
\ar[r] ^-{\smash{\overset{\L}{\otimes}} ^\dag _{\O _{\PP,\Q} } [-d _{P}]}
\ar[d] _-{\R \underline{\Gamma} ^\dag _{X} f  ^! \times \R \underline{\Gamma} ^\dag _{X} f  ^! }
^-\cong
&
{F\text{-}D ^\mathrm{b} _\mathrm{surhol}  (\PP, T, X/K)}
\ar[d] ^-{\R \underline{\Gamma} ^\dag _{X} f  ^! }
_-\cong
\\
{F\text{-}D ^\mathrm{b} _\mathrm{surhol}  (\widetilde{\PP}, \widetilde{T}, X/K)
\times F\text{-}D ^\mathrm{b} _\mathrm{surhol}  (\widetilde{\PP}, \widetilde{T}, X/K)}
\ar[r] ^-{\smash{\overset{\L}{\otimes}} ^\dag _{\O _{\widetilde{\PP},\Q} } [-d _{\widetilde{P}}]}
&
{F\text{-}D ^\mathrm{b} _\mathrm{surhol}  (\widetilde{\PP}, \widetilde{T}, X/K).}
}
\end{equation}

On procède de même concernant le produit tensoriel externe.
\end{proof}

On dispose du théorème ci-dessous annoncé dans \cite[2.2.4]{caro-stab-holo} (avec des notations légèrement différentes):
\begin{theo}
\label{theo-stab-prod-hol}
Soient $Y$, $Y'$ deux variétés quasi-projectives. 
Avec les notations de \ref{nota-holY}, 
on bénéficie des foncteurs produits tensoriel interne et externe suivants:
\begin{align}
\smash{\overset{\L}{\boxtimes}}   ^{\dag} _{\O _\S}
\colon 
&
F \text{-}D ^\mathrm{b} _\mathrm{hol}  (\D ^\dag _{Y'/K}) \times
F \text{-}D ^\mathrm{b} _\mathrm{hol}  (\D ^\dag _{Y/K})
\rightarrow
F \text{-}D ^\mathrm{b} _\mathrm{hol} (\D ^\dag _{Y'\times Y/K}),
\\
\smash{\overset{\L}{\otimes}}^{\dag}  _{\O _{Y/K}} 
 \colon
 &
 F\text{-}D ^\mathrm{b} _\mathrm{hol}  (\D ^\dag _{Y/K}) 
 \times 
 D ^\mathrm{b} _\mathrm{hol}  (\D ^\dag _{Y/K})  
 \to 
 D ^\mathrm{b} _\mathrm{hol}  (\D ^\dag _{Y/K}) .
\end{align}

\end{theo}

\begin{proof}
Berthelot a vérifié que le produit tensoriel externe préserve l'holonomie (voir l'exemple \cite[5.3.5.(v)]{Beintro2}).
On vérifie de manière analogue à \ref{theo-stab-prod-surhol} l'indépendance par rapport au choix de l'immersion de $Y$ ou $Y'$ dans un $\V$-schéma formel projectif et lisse. La stabilité de l'holonomie par produit tensoriels s'en déduit alors via \ref{rema-otimes-boxtimes} (et le fait que l'holonomie est stable par image inverse extraordinaire dans ce cas). 
Une preuve alternative de ce théorème \ref{theo-stab-prod-hol} est d'invoquer \ref{theo-stab-prod-surhol}.
\end{proof}

\bibliographystyle{smfalpha}
%\bibliography{Bibliotheque}
\providecommand{\bysame}{\leavevmode ---\ }
\providecommand{\og}{``}
\providecommand{\fg}{''}
\providecommand{\smfandname}{et}
\providecommand{\smfedsname}{\'eds.}
\providecommand{\smfedname}{\'ed.}
\providecommand{\smfmastersthesisname}{M\'emoire}
\providecommand{\smfphdthesisname}{Th\`ese}

\bigskip
\noindent Daniel Caro\\
Laboratoire de Mathématiques Nicolas Oresme\\
Université de Caen
Campus 2\\
14032 Caen Cedex\\
France.\\
email: daniel.caro@unicaen.fr

\end{document}